\documentclass[reqno]{amsart}
\usepackage{amsfonts,amsmath,amssymb}

\numberwithin{equation}{section}

\newtheorem{thm}{Theorem}[section]
\newtheorem{lem}{Lemma}[section]

\newtheorem{cor}{Corollary}[section]
\newtheorem{prop}{Proposition}[section]

\newtheorem{Step}{Step}[section]

\begin{document}
\title[Critical Elliptic Systems]
{Critical Elliptic Systems in Potential Form}
\author{Emmanuel Hebey}
\address{Emmanuel Hebey, Universit\'e de Cergy-Pontoise, 
D\'epartement de Math\'ematiques, Site de 
Saint-Martin, 2 avenue Adolphe Chauvin, 
95302 Cergy-Pontoise cedex, 
France}
\email{Emmanuel.Hebey@math.u-cergy.fr}

\date{July 18, 2005}

\begin{abstract} Let $(M,g)$ be a smooth compact Riemannian $n$-manifold, $n \ge 3$. Let also $p \ge 1$ be an integer, 
and $M_p^s({\mathbb R})$ be the vector space of symmetrical $p\times p$ real matrix. 
For $A: M \to M_p^s({\mathbb R})$ 
smooth, $A = (A_{ij})$, we consider vector valued equations, or systems, like
$$\Delta_g^p{\mathcal U} + A(x){\mathcal U} = \frac{1}{2^\star}D_{{\mathcal U}}\vert{\mathcal U}\vert^{2^\star}
\hskip.1cm ,$$
where ${\mathcal U}: M \to {\mathbb R}^p$ is a $p$-map, 
$\Delta_g^p$ is the Laplace-Beltrami operator 
acting on $p$-maps, and $2^\star$ is critical from the Sobolev viewpoint. 
We investigate various questions 
for this equation, like the existence of minimizing solutions, 
the existence of high energy solutions, 
blow-up theory, and compactness. We provide the complete $H_1^2$-theory for blow-up, 
sharp pointwise estimates, and prove compactness when the equations are not trivially coupled 
and of geometric type.
\end{abstract}

\maketitle

\tableofcontents

Let $(M,g)$ be a smooth compact Riemannian $n$-manifold, $n \ge 3$. Let also $p \ge 1$ be an integer, 
and $M_p^s({\mathbb R})$ be the vector space of symmetrical $p\times p$ real matrix. 
Namely the vector space of $p\times p$ real matrix $S = (S_{ij})$ which are such that $S_{ij} = S_{ji}$ 
for all $i, j$. 
For $A: M \to M_p^s({\mathbb R})$ 
smooth, $A = (A_{ij})$, we consider vector valued equations like
\begin{equation}\label{GenericEqtAlgForm}
\Delta_g^p{\mathcal U} + A(x){\mathcal U} = \frac{1}{2^\star}D_{{\mathcal U}}\vert{\mathcal U}\vert^{2^\star}
\hskip.1cm ,
\end{equation}
where ${\mathcal U}: M \to {\mathbb R}^p$ is a map, sometimes referred to 
as a $p$-map to underline the fact that the target space is ${\mathbb R}^p$, 
$\Delta_g^p$ is the Laplace-Beltrami operator 
acting on $p$-maps, $2^\star = 2n/(n-2)$, 
and $D_{{\mathcal U}}$ is the derivation operator with respect to ${\mathcal U}$. 
Writing that ${\mathcal U} = (u_1,\dots,u_p)$, we have that  
$\vert{\mathcal U}\vert^{2^\star} = \sum_{i=1}^p\vert u_i\vert^{2^\star}$, that 
$\frac{1}{2^\star}D_{{\mathcal U}}\vert{\mathcal U}\vert^{2^\star} = 
\left(\vert u_i\vert^{2^\star-2}u_i\right)_i$, and 
that $\Delta_g^p{\mathcal U} = \left(\Delta_gu_i\right)_i$, 
where $i = 1,\dots,p$, and 
$\Delta_g = -\hbox{div}_g\nabla$ is the Laplace-Beltrami 
operator for functions. Another way we can write 
(\ref{GenericEqtAlgForm}) is like in the following elliptic system in potential form:
\begin{equation}\label{GenericEqt}
\Delta_gu_i + \sum_{j=1}^pA_{ij}(x)u_j = \vert u_i\vert^{2^\star-2}u_i\hskip.1cm ,
\end{equation}
where the equations have to be satisfied in $M$, and for all $i = 1,\dots,p$. We say the system is of order $p$, and refer to it 
as a $p$-system. The 
system has a variational structure. It is also critical from the Sobolev viewpoint since, if 
$H_1^2$ is the Sobolev space of functions in $L^2$ with 
one derivative in $L^2$, then $2^\star$ is the critical Sobolev exponent for the embeddings of 
$H_1^2$ into Lebesgue spaces. When $p=1$, (\ref{GenericEqt}) reduces to Yamabe type equations like 
\begin{equation}\label{YamabeTypeEqt}
\Delta_gu + h(x)u = \vert u\vert^{2^\star-2}u\hskip.1cm ,
\end{equation}
where $h: M \to {\mathbb R}$, $h = A_{11}$,  is smooth. In the same order of ideas, if 
$A$ is diagonal, then (\ref{GenericEqt}) consists of $p$ copies of 
Yamabe type equations like (\ref{YamabeTypeEqt}) with $h = A_{ii}$, $i = 1,\dots,p$. When $A$ is diagonal, 
the equations are independent one from another, and the system is said to be trivial, or trivially coupled.

\medskip Equations like (\ref{GenericEqtAlgForm}) generalize the intensively studied 
Yamabe type equations. In this paper we study  
various questions  for (\ref{GenericEqtAlgForm}), like the existence of minimizing solutions  for (\ref{GenericEqtAlgForm}), 
the existence of high energy solutions, 
blow-up theory, and compactness properties of (\ref{GenericEqtAlgForm}). 
Section \ref{Prel} of the paper is devoted to preliminary definitions 
and remarks on our systems. We insist in this section on some of the differences 
we can find between the scalar case (where $p = 1$) and the vector valued 
case addressed in this paper (where $p \ge 2$). A major difference is the lack of maximum principle. 
We exhibit in Section \ref{Prel} examples of nonnegative solutions of systems like 
(\ref{GenericEqtAlgForm}) which are not positive. In other words,  
examples of solutions such that the factors 
of the solutions are nonnegative,  nonzero, but with 
zeros in $M$ -- a phenomenon which 
does not occur when $p = 1$ (where, by the maximum principle, nonnegative solutions 
are either identically zero or everywhere positive). In Section 
\ref{ExSol} of the paper we discuss the variational structure attached to our systems. 
In particular, we prove that if the minimum energy of the system 
is sufficiently small, then the system possesses a minimizing solution (of small energy). The existence of high energy solutions 
is discussed in Section \ref{HighEn}. Blow-up theory for our systems 
and for nonnegative solutions is discussed in Sections 
\ref{SobTheory} to \ref{StdRescaling}. We adopt in these sections the dynamical 
viewpoint (in the sense of the terminology introduced in \cite{HebBBConf}) 
which consists in considering sequences of solutions (or Palais-Smale sequences) of families 
of equations and not only sequences of solutions 
(or Palais-Smale sequences) of one given equation (in our case, 
of  families of systems and not of one given system). The $H_1^2$-theory for blow-up is discussed in Section 
\ref{SobTheory}. We provide in this section a complete description of the asymptotic behaviour of Palais-Smale 
sequences in the Sobolev space $H_{1,p}^2$ consisting 
of $p$-maps whose components are in $L^2$ with one derivative in $L^2$. 
Global pointwise estimates for sequences of solutions are discussed in Section \ref{PointEst}. 
We prove that we can add pointwise estimates to the Sobolev asymptotics of the 
$H_1^2$-theory when passing from Palais-Smale sequences to sequences of solutions. 
The notion of 
$L^2$-concentration, which turns out to be important 
for compactness issues, is addressed in Section \ref{L2Conc}. 
In Section \ref{SharpPtAsympt} 
we prove sharp local pointwise estimates 
for sequences of solutions. Standard rescaling 
is discussed in Section \ref{StdRescaling}. The 
blow-up theory developed in Sections \ref{SobTheory} to \ref{StdRescaling} is applied in Section 
\ref{Cptness} to get compactness results in the conformally flat case. 
Roughly speaking, we prove that, for conformally flat manifolds, systems 
like (\ref{GenericEqtAlgForm}) 
are compact when their coupling does not involve a trivial coupling related to the 
geometric equation.

\section{Preliminary definitions and remarks}\label{Prel}

In what follows we let $(M,g)$ be a smooth compact Riemannian manifold of dimension $n \ge 3$,  $p \ge 1$ integer, 
and $A: M \to M_p^s({\mathbb R})$, $A = (A_{ij})$, be a smooth map from $M$ to 
$M_p^s({\mathbb R})$. We consider systems like
\begin{equation}\label{GenericEqtSec1}
\Delta_gu_i + \sum_{j=1}^pA_{ij}(x)u_j = \Lambda\vert u_i\vert^{2^\star-2}u_i
\end{equation}
in $M$, for all $i = 1,\dots p$, where $\Lambda \in {\mathbb R}$. We let 
$H_{1,p}^2(M)$ be the space consisting of $p$-maps 
${\mathcal U}: M \to {\mathbb R}^p$, ${\mathcal U} = (u_1,\dots,u_p)$, 
such that the $u_i$'s are all in the standard Sobolev space $H_1^2(M)$. Namely,
$$H_{1,p}^2(M) = \Bigl\{{\mathcal U} = (u_1,\dots,u_p)\hskip.1cm ,\hskip.1cm 
u_i \in H_1^2(M)\hskip.1cm\hbox{for all}\hskip.1cm i\Bigr\}\hskip.1cm ,$$
where $H_1^2(M)$ is the standard Sobolev space of functions in $L^2(M)$ with 
one derivative in $L^2$. 
The space $H_{1,p}^2(M)$ is an Hilbert space when equipped 
with the scalar product
$$\langle{\mathcal U},{\mathcal V}\rangle_{H_{1,p}^2} = \sum_i(u_i,v_i)\hskip.1cm ,$$
where 
$(\cdot,\cdot)$ is the usual scalar product in $H_1^2(M)$, the $u_i$'s are the components of 
${\mathcal U}$, and the $v_i$'s are the components of ${\mathcal V}$.
A map ${\mathcal U} = (u_1,\dots,u_p)$ in $H_{1,p}^2(M)$ is said to be 
a weak solution of (\ref{GenericEqtSec1}) if the equations in (\ref{GenericEqtSec1}) are satisfied in the distributional 
sense by the $u_i$'s. By regularity theory (see the proof of 
Theorem \ref{ExistThm} in Section \ref{ExSol}), any weak solution ${\mathcal U} = (u_1,\dots,u_p)$ of (\ref{GenericEqtSec1}), 
${\mathcal U} \in H_{1,p}^2(M)$, is 
in $C^{2,\theta}$ for all $0 < \theta < 1$ (in the sense 
that the $u_i$'s are in $C^{2,\theta}$ for all $i$). 

\medskip In what follows we say that a $p$-map 
${\mathcal U} = (u_1,\dots,u_p)$ is
{\it nonnegative} if the $u_i$'s are all nonnegative functions (i.e $u_i \ge 0$ for all $i$), 
{\it weakly positive} if the $u_i$'s are all positive functions unless they are identically zero 
(i.e, for any $i$, either $u_i > 0$ or $u_i \equiv 0$), and 
{\it strongly positive} if the $u_i$'s are all positive functions (i.e $u_i > 0$ for all $i$). For short, a $p$-map is said 
to be {\it positive} if it is either weakly positive, or strongly positive.
Following standard terminology in the elliptic 
system literature, we say that a matrix $S = (S_{ij})$ is 
{\it cooperative} if $S_{ij} \ge 0$ for all $i \not= j$. If $S: M \to M_p^s({\mathbb R})$ is a map, 
$S$ is said to be 
cooperative in $M$ if $S_{ij}(x) \ge 0$ for all $i \not= j$, and all $x \in M$. At last, still following standard 
terminology in the elliptic system literature, we say that the system 
(\ref{GenericEqtSec1}) is {\it fully coupled} if the index set $\{1,\dots,p\}$ does not split in two disjoint subsets 
$\{i_1,\dots,i_k\}$ and $\{j_1,\dots,j_{k^\prime}\}$, $k + k^\prime = p$, such that $A_{i_\alpha j_\beta} \equiv 0$ 
for all $\alpha = 1,\dots,k$ and $\beta = 1,\dots,k^\prime$. When (\ref{GenericEqtSec1}) is not fully coupled, up to 
permuting the equations, $A$ can be written in diagonal blocks 
$S: M \to M_k^s({\mathbb R})$ and $T: M \to M_{p-k}^s({\mathbb R})$ for some $k < p$, like in (\ref{ExSec1Matrix}) below, 
and the $p$-system splits in two independent systems. Namely a $k$-system with respect to $S$ and a 
$(p-k)$-system with respect to $T$. If ${\mathcal U} = (u_i)_i$ is a solution of (\ref{GenericEqtSec1}), 
and $\sigma$ is a permutation of $\Sigma_p = \left\{1,\dots,p\right\}$, then $\tilde{\mathcal U} = (u_{\sigma(i)})_i$ 
is  a solution of (\ref{GenericEqtSec1}) when we replace the $A_{ij}$'s by $\tilde A_{ij} = A_{\sigma(i)\sigma(j)}$. 
Possible references on elliptic systems are De Figueiredo \cite{DeF}, 
De Figueiredo and Ding \cite{DeFDin}, De Figueiredo and Felmer \cite{DeFFel}, 
El Hamidi \cite{ElH}, 
Hulshof, Mitidieri and Van der Vorst \cite{HulMitVan}, Mancini and Mitidieri \cite{ManMit}, 
Mitidieri and Sweers \cite{MitSwe}, and Sweers \cite{Swe}.

\medskip In Remarks 1.1 and 1.2 below we discuss two important examples which show that nonnegative 
solutions of systems like (\ref{GenericEqtSec1}) 
are not necessarily weakly positive, and, in the same order of ideas, 
that weakly positive solutions 
of systems like (\ref{GenericEqtSec1}) 
are not necessarily strongly positive.

\medskip\noindent{\bf Remark 1.1:} 
Contrary to the case $p = 1$, nonnegative solutions of a system like (\ref{GenericEqtSec1}) 
are not necessarily weakly positive (and thus not necessarily 
strongly positive as well). We may have a nonnegative solution ${\mathcal U} = (u_1,\dots,u_p)$ with one of the 
$u_i$'s non identically zero, but with zeros in $M$. A possible construction 
when $p=2$ is as 
follows. We let $(S^n,g_0)$ be the unit $n$-sphere. The Yamabe equation on $(S^n,g_0)$ reads as
\begin{equation}\label{YamEqtSphereEx}
\Delta_{g_0}u + \frac{n(n-2)}{4}u = u^{2^\star-1}
\end{equation}
and the positive solutions of (\ref{YamEqtSphereEx}) are given by
\begin{equation}\label{YamEqtSphereSol}
u_{x_0}^\lambda(x) = \left(\frac{n(n-2)}{4}(\lambda^2-1)\right)^{\frac{n-2}{4}}
\left(\lambda - \cos d_{g_0}(x_0,x)\right)^{1-\frac{n}{2}}\hskip.1cm ,
\end{equation}
where $\lambda > 1$ and $x_0 \in S^n$. We fix $x_0 \in S^n$ and let $u_\lambda = u_{x_0}^\lambda$. Then 
\begin{equation}\label{ExSphereEqt1}
u_\lambda = \left(\frac{n(n-2)}{4}\right)^{\frac{n-2}{4}} + \varepsilon_\lambda
\end{equation}
where $\varepsilon_\lambda \to 0$ in $C^0(S^n)$ as $\lambda \to +\infty$. Let $m_\lambda = \min_{S^n}u_\lambda$. Let 
also $\tilde\varepsilon_\lambda$ and $\hat\varepsilon_\lambda$ be the functions given by
\begin{equation}\label{ExSphereEpsiloEqt}
\begin{split}
&\tilde\varepsilon_\lambda(x) = \frac{u_\lambda(x) - m_\lambda}{u_\lambda(x)}\hskip.1cm ,\hskip.1cm\hbox{and}\\
&\hat\varepsilon_\lambda(x) = 
\frac{n(n-2)}{4} - u_\lambda(x)^{2^\star-2} 
+ \frac{\left(u_\lambda(x) - m_\lambda\right)^{2^\star-1}}{u_\lambda(x)}\hskip.1cm ,
\end{split}
\end{equation}
where $x \in S^n$. 
By (\ref{ExSphereEqt1}), $\tilde\varepsilon_\lambda \to 0$ and $\hat\varepsilon_\lambda \to 0$ in 
$C^0(S^n)$ as $\lambda \to +\infty$ while, by (\ref{YamEqtSphereEx}), letting 
$u_1^\lambda = u_\lambda - m_\lambda$ and $u_2^\lambda = u_\lambda$, we easily get that 
${\mathcal U}_\lambda = (u_1^\lambda,u_2^\lambda)$ is a solution of the $2$-system
\begin{equation}\label{ExSphereEqt2}
\Delta_{g_0}u_i + \sum_{j=1}^2A_{ij}^\lambda(x)u_j = \vert u_i\vert^{2^\star-2}u_i
\end{equation}
in $S^n$, for all $i = 1,2$, if
$\tilde\varepsilon_\lambda A_{11}^\lambda + A_{12}^\lambda = \hat\varepsilon_\lambda$ 
and $\tilde\varepsilon_\lambda A_{21}^\lambda + A_{22}^\lambda = \frac{n(n-2)}{4}$. 
Let $A_{11}^\lambda = \frac{n(n-2)}{4}$. Then the $2$-map 
${\mathcal U}_\lambda = (u_1^\lambda,u_2^\lambda)$ is a solution of (\ref{ExSphereEqt2}) for 
all $\lambda$, where 
$A^\lambda = (A^\lambda_{ij})$ can be chosen in the form
\begin{equation}\label{ExSphereMatrixFinalEqt}
A^\lambda = \left(
\begin{matrix}
\frac{n(n-2)}{4} & \varepsilon^\prime_\lambda\\
\varepsilon^\prime_\lambda & \frac{n(n-2)}{4} + \varepsilon^{\prime\prime}_\lambda
\end{matrix}
\right)\hskip.1cm ,
\end{equation}
and where $\varepsilon^\prime_\lambda \to 0$ and $\varepsilon^{\prime\prime}_\lambda \to 0$ in 
$C^0(S^n)$ as $\lambda \to +\infty$. Obviously, $u_1^\lambda = u_\lambda - m_\lambda$ 
is nonnegative with one zero (at $-x_0$). In 
particular, ${\mathcal U}_\lambda = (u_1^\lambda,u_2^\lambda)$ 
is not weakly positive. However, ${\mathcal U}_\lambda$ is 
a nonnegative solution of the $2$-system (\ref{ExSphereEqt2}). 
Summarizing, we constructed a nonnegative solution 
of a system like (\ref{GenericEqtSec1}) with one factor in the solution which is non identically zero, but has zeros in $M$. 
This shows that nonnegative solutions of systems like (\ref{GenericEqtSec1}) 
are not necessarily weakly positive. By (\ref{ExSphereMatrixFinalEqt}), 
for $\lambda > 1$ sufficiently large, 
$A^\lambda$ is positive and the operator 
$\Delta_{g_0}^2+A^\lambda$ is coercive 
(coercivity is defined in Section \ref{ExSol}). Of course, (\ref{ExSphereEqt2}) is also 
fully coupled since, if not, the maximum principle for functions would lead to a 
contradiction.

\medskip\noindent{\bf Remark 1.2:}  In a similar way, 
weakly positive solutions of a system like (\ref{GenericEqtSec1}) 
are not necessarily strongly positive. There is of course a trivial construction where we create $0$-factors 
by adding ``artificial'' equations to the system.  For instance, if ${\mathcal U}$ 
is a solution of (\ref{GenericEqtSec1}) with $\Lambda = 1$, 
then we may regard the map $({\mathcal U},0)$, $0 \in {\mathbb R}$, as a solution of the 
system (\ref{GenericEqtSec1})--(\ref{YamabeTypeEqt}) we get by adding to (\ref{GenericEqtSec1}) 
one equation like (\ref{YamabeTypeEqt}). In such constructions, the resulting system is obviousy not fully coupled. However, 
there are easy examples of fully coupled systems with weakly positive solutions 
having zero factors (and thus with weakly 
positive solutions which are not strongly positive). To get such examples we may consider a positive solution $u$ of 
an equation like (\ref{YamabeTypeEqt}), and note that ${\mathcal U} = (u,u,0)$ is a solution of (\ref{GenericEqtSec1}) 
with $p = 3$ (and $\Lambda = 1$) when  the $A_{ij}$'s are such that $A_{ij} = A_{ji}$, 
$A_{11} + A_{12} = h$, $A_{21} + A_{22} = h$, and $A_{13} = - A_{23}$. For instance, 
if $u$ is a positive solution of (\ref{YamabeTypeEqt}), then ${\mathcal U} = (u,u,0)$ is a 
solution of (\ref{GenericEqtSec1}) 
with $p = 3$ and $\Lambda = 1$, when
\begin{equation*}
A = \left(
\begin{matrix}
\frac{1}{2}h & \frac{1}{2}h & \alpha\\
\frac{1}{2}h & \frac{1}{2}h & -\alpha\\
\alpha & -\alpha & \beta
\end{matrix}
\right)\hskip.1cm ,
\end{equation*}
$h$ is as in (\ref{YamabeTypeEqt}), and $\alpha$, $\beta$ are 
arbitrary functions. The system is obviously fully coupled when $\alpha$ is nonzero.

\medskip The following lemma, that we will use several times in the sequel, shows that 
the above examples in Remarks 1.1 and 1.2 
 stop to hold, or reduce to the trivial case, when $-A$ is cooperative, 
 respectively when 
 $-A$ is cooperative and the system is also fully coupled.
 
\begin{lem}\label{PosSolSec1} Let $(M,g)$ be a smooth compact Riemannian manifold of dimension 
$n \ge 3$,  $p \ge 1$ integer, 
and $A: M \to M_p^s({\mathbb R})$, $A = (A_{ij})$, be a smooth map from $M$ to 
$M_p^s({\mathbb R})$. If $-A$ is cooperative, then any nonnegative solution of (\ref{GenericEqtSec1})  
is weakly positive. If $-A$ is cooperative and (\ref{GenericEqtSec1}) is 
fully coupled, then any weakly positive solution of (\ref{GenericEqtSec1}) is strongly positive.
\end{lem}

\begin{proof} If $-A$ is cooperative,  and ${\mathcal U} = (u_1,\dots,u_p)$ is a solution of 
(\ref{GenericEqtSec1}), we can write that
$$\Delta_gu_i + A_{ii}u_i \ge \Lambda u_i^{2^\star-1}$$
for all $i$. Then the classical maximum principle for scalar equations can be applied so that either $u_i  > 0$ 
everywhere in $M$, or $u_i \equiv 0$. This proves the first assertion in the lemma. 
Concerning the second assertion, we prove that if $-A$ (or $A$) is 
cooperative, and ${\mathcal U}$ is a weakly positive solution 
of the system, then $A$ can be factorized in blocs with respect to the zero and 
nonzero components of ${\mathcal U}$. More precisely, if we write ${\mathcal U} = (u_1,\dots,u_k,0,\dots,0)$ 
with $k < p$, and $u_i > 0$ for all $i$, then
\begin{equation}\label{ExSec1Matrix}
A = \left(
\begin{matrix}
S & 0\\
0 & T
\end{matrix}
\right)\hskip.1cm ,
\end{equation}
where $S: M \to M_k^s({\mathbb R})$, $T: M \to M_{p-k}^s({\mathbb R})$, and the $0$'s are null matrix of 
respective order $k\times(p-k)$ and $(p-k)\times k$. This easily follows from the equations $\sum_{j=1}^kA_{ij}u_j = 0$ 
for all $i \ge k+1$, so that we necessarily have that $A_{ij} = 0$ for all $i \ge k+1$ 
and $j \le k$. The $p$-system (\ref{GenericEqtSec1}) splits into two independent systems -- 
a $k$-system where $A$ is replaced by $S$, and a $(p-k)$-system 
where $A$ is replaced by $T$. In particular, if $-A$ is cooperative and (\ref{GenericEqtSec1}) is 
fully coupled, any weakly positive solution of (\ref{GenericEqtSec1}) is strongly positive.
\end{proof}

 Remark 1.3 below is concerned with 
 the negative case where $\Lambda < 0$ in (\ref{GenericEqtSec1}). With 
respect to scalar equations 
(where $p=1$) we loose uniqueness but still have the a priori $L^\infty$-bound. Like when 
$p=1$, the difficult case in equations like (\ref{GenericEqtSec1}), 
with respect to existence or compactness of solutions, is when $\Lambda$ is positive.

\medskip\noindent{\bf Remark 1.3:}  
There are different behaviour of (\ref{GenericEqtSec1}) depending on the sign of 
$\Lambda$. In the negative case, when $\Lambda < 0$, the positive solution is known to be unique 
if $p = 1$. This obviously stops to hold for systems and there are counter examples involving only 
constant maps, and so only algebra. 
Suppose for instance that $p= 2$ and $n = 6$. Let $\lambda > 0$, $u_\lambda = \lambda$, 
$v_\lambda = \lambda +1$, and $w_\lambda = (\lambda^2+(\lambda+1)^2)/(2\lambda + 1)$. Let also 
$A_{11}(\lambda)$, $A_{22}(\lambda)$, and $A_{12}(\lambda) = A_{21}(\lambda)$ 
be given by $A_{11}(\lambda) = A_{22}(\lambda) = - (3\lambda^2+3\lambda+1)/(2\lambda+1)$ 
and $A_{12}(\lambda) = (\lambda^2+\lambda)/(2\lambda+1)$. Then, for any $\lambda > 0$, $(u_\lambda,v_\lambda)$, 
$(v_\lambda,u_\lambda)$, and $(w_\lambda,w_\lambda)$ are three distinct strongly positive solutions of (\ref{GenericEqtSec1}) 
with $\Lambda = -1$ and $A = A(\lambda)$, where 
$A(\lambda)$ is the matrix of components $A_{ij}(\lambda)$, $i,j = 1,2$. 
However, though uniqueness is not anymore true for systems, it is still true, like when $p = 1$, that there is a bound (depending 
only on $A$ and $\Lambda$) for the $L^\infty$-norm of positive solutions of the system.  If ${\mathcal U} = (u_1,\dots,u_p)$ 
is a positive solution of (\ref{GenericEqtSec1}) with $\Lambda < 0$, 
and $x_i$ is a point where $u_i$ is maximum, then $\Delta_gu_i(x_i) \ge 0$ 
and we get  with the equations in (\ref{GenericEqtSec1}) that
$$\Vert u_i\Vert_\infty^{2^\star-1} \le \frac{1}{\vert\Lambda\vert} \Vert A\Vert_\infty \max_{j=1,\dots,p}\Vert u_j\Vert_\infty$$
for all $i$, where $\Vert A\Vert_\infty = \max_{i,j}\Vert A_{ij}\Vert_\infty$. In particular, there exists 
$C = C(\Vert A\Vert_\infty,\Lambda)$, $C > 0$ depending only on $\Vert A\Vert_\infty$ and $\Lambda$, such that 
$\Vert u_i\Vert_\infty \le C$ for all $i$. 

\section{Minimizing solutions and variational structure}\label{ExSol}

By the Sobolev 
inequality in Euclidean space, there exists $K > 0$ such that 
$\Vert u\Vert_{2^\star} \le K\Vert\nabla u\Vert_2$ for all smooth functions $u$ with compact support 
in ${\mathbb R}^n$. We let $K_n$ be the sharp constant $K$ in this inequality, so that
\begin{equation}\label{ValueSharpCst}
K_n = \sqrt{\frac{4}{n(n-2)\omega_n^{2/n}}}\hskip.1cm ,
\end{equation}
where $\omega_n$ is the volume of the unit $n$-sphere. A preliminary claim is that we also have that
\begin{equation}\label{SharpCstSystEqt}
K_n^{-2} = \inf_{{\mathcal U}\in{\mathcal H}^p_{{\mathbb R}^n}}\sum_{i=1}^p\int_{{\mathbb R}^n}\vert\nabla u_i\vert^2dx
\end{equation}
for all $p \ge 1$, where ${\mathcal H}^p_{{\mathbb R}^n}$ consists of the maps 
${\mathcal U} = (u_1,\dots,u_p)$, $u_i: {\mathbb R}^n \to {\mathbb R}$, 
which are such that $u_i \in D_1^2({\mathbb R}^n)$ for all $i$, and 
$\sum_{i=1}^p\int_{{\mathbb R}^n}\vert u_i\vert^{2^\star}dx = 1$, 
where $D_1^2({\mathbb R}^n)$ is the Beppo-Levi space defined as the completion 
of $C^\infty_0({\mathbb R}^n)$, 
the space of smooth function with compact support in ${\mathbb R}^n$, with respect to the norm 
$\Vert u\Vert = \Vert\nabla u\Vert_2$. 
When $p = 1$, (\ref{SharpCstSystEqt}) 
is obvious since it 
reduces to the basic definition of the sharp constant $K$ in the inequality $\Vert u\Vert_{2^\star} \le K\Vert\nabla u\Vert_2$. 
That (\ref{SharpCstSystEqt})  is also true for $p > 1$ 
was first noticed by Amster, De N\'apoli and Mariani \cite{AmsNapMar}. Let $\Lambda_p$ be the right hand side 
in (\ref{SharpCstSystEqt}). 
Taking ${\mathcal U} = (u,0,\dots,0)$, we easily get that $\Lambda_p \le \Lambda_1$. Conversely, 
since $2/2^\star < 1$, so that $(\sum\vert a_i\vert)^{2/2^\star} \le \sum\vert a_i\vert^{2/2^\star}$, 
$a_i \in {\mathbb R}$, we can write that
\begin{equation}\label{IntroEqt1}
\begin{split}
\left(\sum_{i=1}^p\int_{{\mathbb R}^n}\vert u_i\vert^{2^\star}dx\right)^{2/2^\star}
&\le \left(\sum_{i=1}^p\left(\frac{1}{\Lambda_1}
\int_{{\mathbb R}^n}\vert\nabla u_i\vert^2dx\right)^{2^\star/2}\right)^{2/2^\star}\\
&\le \frac{1}{\Lambda_1} \sum_{i=1}^p
\int_{{\mathbb R}^n}\vert\nabla u_i\vert^2dx
\end{split}
\end{equation}
and it follows that we also have that $\Lambda_p \ge \Lambda_1$. This proves 
the above claim that (\ref{SharpCstSystEqt}) is true 
for all $p$. Now, we return to the case of manifolds. For  
$(M,g)$ a smooth compact Riemannian manifold of dimension $n \ge 3$, 
and $A: M \to M_p^s({\mathbb R})$, $A = (A_{ij})$ smooth, we let 
$I_A: H_{1,p}^2(M) \to {\mathbb R}$ be the functional given by
\begin{equation}\label{FonctionalDef}
I_A({\mathcal U}) = \sum_{i=1}^p\int_M\vert\nabla u_i\vert^2dv_g + 
\sum_{i,j = 1}^p\int_MA_{ij}(x)u_iu_jdv_g\hskip.1cm ,
\end{equation}
where $H_{1,p}^2(M)$ is as in Section \ref{Prel}. Let also $\Phi: H_{1,p}^2(M) \to {\mathbb R}$ 
be the functional given by
$\Phi({\mathcal U}) = \sum_{i=1}^p\int_M\vert u_i\vert^{2^\star}dv_g$. 
The definition of $\Phi$ makes sense thanks to the Sobolev embedding of $H_1^2$ into $L^{2^\star}$. 
We define $\mu_{A,g}^p(M)$ by
\begin{equation}\label{MiniDefEqt}
\mu_{A,g}^p(M) = \inf_{{\mathcal U} \in {\mathcal H}^p_M} I_A({\mathcal U})\hskip.1cm ,
\end{equation}
where ${\mathcal H}^p_M$ consists of the maps ${\mathcal U} \in H_{1,p}^2(M)$ which satisfy the constraint 
$\Phi({\mathcal U}) = 1$. Let $S: M \to {\mathbb R}$ be a function, and $Id_p$ be the $p\times p$-identity matrix. 
We let $\mu_{S,g}^1(M)$ be the infimum $\mu^1_{S Id_1,g}(M)$. An easy claim in the spirit of (\ref{SharpCstSystEqt}) 
is that
\begin{equation}\label{TrivialCouplingMin}
\mu_{S Id_p,g}^p(M) = \mu_{S,g}^1(M)
\end{equation}
if $S$ is such that $\mu_{S,g}^1(M) \ge 0$ (for instance $S \ge 0$). 
In order to prove this claim, we first note that, taking ${\mathcal U} = (u,0,\dots,0)$ in (\ref{MiniDefEqt}), 
we get that 
$\mu_{S Id_p,g}^p(M) \le \mu_{S,g}^1(M)$. We assume $\mu_{S,g}^1(M) \ge 0$. It follows that we also have 
that $\mu_{S Id_p,g}^p(M) \ge 0$. In particular, if  $\mu_{S,g}^1(M) = 0$, then 
$\mu_{S Id_p,g}^p(M) = 0$, and (\ref{TrivialCouplingMin}) is true in this case. Let us now assume that 
$\mu_{S,g}^1(M) > 0$. As in 
(\ref{IntroEqt1}), we may write that
\begin{eqnarray*}
&&\left(\sum_{i=1}^p\int_M\vert u_i\vert^{2^\star}dv_g\right)^{2/2^\star}\\
&&\hskip.4cm \le \left(\sum_{i=1}^p\left(\frac{1}{\mu_{S,g}^1(M)}
\int_M\left(\vert\nabla u_i\vert^2 + Su_i^2\right)dv_g\right)^{2^\star/2}\right)^{2/2^\star}\\
&&\hskip.4cm \le \frac{1}{\mu_{S,g}^1(M)} \sum_{i=1}^p
\int_M\left(\vert\nabla u_i\vert^2 + Su_i^2\right)dv_g
\end{eqnarray*}
and we get that 
$I_{S Id_p}({\mathcal U}) \ge \mu_{S,g}^1(M)\Phi({\mathcal U})^{2/2^\star}$ 
for all maps ${\mathcal U} \in H_{1,p}^2(M)$, where $I_A$ is given by (\ref{FonctionalDef}), and 
$\Phi$ is as in the definition of ${\mathcal H}^p_M$ in (\ref{MiniDefEqt}). As a consequence we get that 
$\mu_{S Id_p,g}^p(M) \ge \mu_{S,g}^1(M)$, and 
this proves the above claim that (\ref{TrivialCouplingMin}) is true when $\mu_{S,g}^1(M) \ge 0$. 
In particular, if $A_g = \frac{n-2}{4(n-1)} S_g Id_p$, where $S_g$ is the scalar curvature of $g$, and if the Yamabe invariant 
$Y_{[g]}(M)$ of $(M,g)$ is nonnegative, then $\mu_{A_g,g}^p(M) = Y_{[g]}(M)$. For instance, if 
$(S^n,g_0)$ is the unit $n$-sphere, then, similarly to (\ref{SharpCstSystEqt}), we can 
also write that $K_n^{-2} = \mu_{A_{g_0},g_0}^p(S^n)$ for all $p$. A possible survey reference for the Yamabe 
material we used here is 
Lee and Parker \cite{LeePar}. We refer also to Aubin \cite{Aub1}, Hebey \cite{Heb1}, and Schoen \cite{Sch3}.

\medskip When $\mu_{S,g}^1(M) < 0$ and $p > 1$, 
an equation like $\mu_{S Id_p,g}^p(M) = \mu_{S,g}^1(M)$ stops to hold and we do get that 
$\mu_{S Id_p,g}^p(M) < \mu_{S,g}^1(M)$. 
In order to see this, we let $u > 0$ be a minimizer for $\mu_{S,g}^1(M)$, and let 
${\mathcal U} \in H_{1,p}^2(M)$ be the $p$-map given by
${\mathcal U} = \left(p^{-1/2^\star}u,\dots,p^{-1/2^\star}u\right)$. Then $\Phi({\mathcal U}) = 1$, 
where $\Phi$ is as in the definition of ${\mathcal H}^p_M$ in (\ref{MiniDefEqt}), and it follows from 
(\ref{FonctionalDef}) that $I_{S Id_p}({\mathcal U}) = p^{1-(2/2^\star)}\mu_{S,g}^1(M)$. 
In particular, 
$\mu_{S Id_p,g}^p(M) < \mu_{S,g}^1(M)$ when $p > 1$ and $\mu_{S,g}^1(M) < 0$. 

\medskip As in Section \ref{Prel}, we let $(M,g)$ be a smooth compact Riemannian manifold of dimension 
$n \ge 3$, $A: M \to M_p^s({\mathbb R})$, $A = (A_{ij})$, be a smooth map, and we consider systems like
\begin{equation}\label{GenericEqtThm1}
\Delta_gu_i + \sum_{j=1}^pA_{ij}(x)u_j = \Lambda\vert u_i\vert^{2^\star-2}u_i
\end{equation}
in $M$, for all $i = 1,\dots p$, where $\Lambda \in {\mathbb R}$. 
We claim here that the following theorem holds. 
The case $p = 1$ in the result is well known and 
goes back to Aubin \cite{Aub} and Trudinger \cite{Tru} (see also Yamabe 
\cite{Yam}). Note that by considering test functions ${\mathcal U} \in 
H_{1,p}^2(M)$ for $I_A$ in the form ${\mathcal U} = (u,0,\dots,0)$, we easily get from the $p=1$ case that 
$\mu_{A,g}^p(M) \le K_n^{-2}$ for all $A$ and all $(M,g)$.

\begin{thm}[Existence of Minimizers]\label{ExistThm} Let $(M,g)$ be a smooth compact 
Riemannian manifold of dimension $n \ge 3$, 
$p \ge 1$, and $A: M \to M_p^s({\mathbb R})$ be a smooth map. Assume $\mu_{A,g}^p(M) < K_n^{-2}$, where $\mu_{A,g}^p(M)$ 
is given by (\ref{MiniDefEqt}), and $K_n$ is given by (\ref{ValueSharpCst}). Then there exists a minimizer 
${\mathcal U}$ for $\mu_{A,g}^p(M)$. Namely, 
${\mathcal U} \in {\mathcal H}^p_M$  and $I_A({\mathcal U}) = \mu_{A,g}^p(M)$, where 
$I_A$ is given by (\ref{FonctionalDef}), and ${\mathcal H}^p_M$ 
is as in (\ref{MiniDefEqt}). In particular, ${\mathcal U} = (u_1,\dots,u_p)$ 
is a minimizing solution of the system (\ref{GenericEqtThm1}) with $\Lambda = \mu_{A,g}^p(M)$.
The $u_i$'s are in $C^{2,\theta}$ for all $i$ and $0 < \theta < 1$. 
The solution ${\mathcal U}$ can be chosen weakly positive if $-A$ is cooperative in $M$, and strongly positive 
if the system is also fully coupled.
\end{thm}

\begin{proof} First we prove the existence of a minimizer 
${\mathcal U}$ for $\mu_{A,g}^p(M)$ when we assume that $\mu_{A,g}^p(M) < K_n^{-2}$. Then we prove the regularity 
of weak solutions, not necessarily minimizing, of systems like (\ref{GenericEqtThm1}). Concerning 
the existence of ${\mathcal U}$, we let $({\mathcal U}_\alpha)_\alpha$ be a minimizing sequence for 
$\mu_{A,g}^p(M)$, and write that ${\mathcal U}_\alpha = (u_\alpha^1,\dots,u_\alpha^p)$ for all $\alpha$. 
Since ${\mathcal U}_\alpha \in {\mathcal H}_M^p$ for all $\alpha$, the $u_\alpha^i$'s 
are uniformly bounded 
in $H_1^2(M)$ for all $i$. Up to passing to a subsequence, we may therefore assume that, 
for any $i$, $u_\alpha^i \rightharpoonup u_i$ weakly in $H_1^2(M)$, $u_\alpha^i \to u_i$ strongly 
in $L^2(M)$, and $u_\alpha^i \to u_i$ almost everywhere in $M$. By the weak convergence in $H_1^2$, 
\begin{equation}\label{ProofExistSec2Eqt1}
\int_M\vert\nabla u_\alpha^i\vert^2dv_g 
= \int_M\vert\nabla(u_\alpha^i-u_i)\vert^2dv_g + \int_M\vert\nabla u_i\vert^2dv_g + o(1)
\end{equation}
for all $i$, where $o(1) \to 0$ as $\alpha \to +\infty$. We also have, for instance by Br\'ezis and Lieb \cite{BreLie}, that 
\begin{equation}\label{ProofExistSec2Eqt2}
\int_M\vert u_\alpha^i\vert^{2^\star}dv_g 
= \int_M\vert u_\alpha^i-u_i\vert^{2^\star}dv_g + \int_M\vert u_i\vert^{2^\star}dv_g + o(1)
\end{equation}
for all $i$ and $\alpha$, 
where $o(1) \to 0$ as $\alpha \to +\infty$. Moreover, by the sharp Sobolev inequality in Hebey and Vaugon 
\cite{HebVau1,HebVau2}, we can write that
\begin{equation}\label{ProofExistSec2Eqt3}
\begin{split}
&\left(\int_M\sum_{i=1}^p\vert u_\alpha^i-u_i\vert^{2^\star}dv_g\right)^{2/2^\star}\\
& \le  \sum_{i=1}^p\left(\int_M\vert u_\alpha^i-u_i\vert^{2^\star}dv_g\right)^{2/2^\star}\\
& \le  K_n^2\sum_{i=1}^p\int_M\vert\nabla(u_\alpha^i-u_i)\vert^2dv_g + 
B \sum_{i=1}^p\int_M\vert u_\alpha^i-u_i\vert^2dv_g
\end{split}
\end{equation}
for all $i$ and $\alpha$, 
where $B > 0$ does not depend on $i$ and $\alpha$. Then, since ${\mathcal U}_\alpha \in {\mathcal H}_M^p$ 
for all $\alpha$, we get by combining (\ref{ProofExistSec2Eqt1})--(\ref{ProofExistSec2Eqt3}) that
\begin{equation}\label{ProofExistSec2Eqt4}
\begin{split}
&\left(1-\sum_{i=1}^p\int_M\vert u_i\vert^{2^\star}dv_g\right)^{2/2^\star}\\
&\le K_n^2 \sum_{i=1}^p\left(\int_M\vert\nabla u_\alpha^i\vert^2dv_g 
- \int_M\vert\nabla u_i\vert^2dv_g\right) + o(1)
\end{split}
\end{equation}
for all $i$ and $\alpha$, 
where $o(1) \to 0$ as $\alpha \to +\infty$. Since $I_A({\mathcal U}_\alpha) = \mu_{A,g}^p(M) + o(1)$ 
for all $\alpha$, and $u_\alpha^i \to u_i$ in $L^2(M)$ 
for all $i$, we get with (\ref{ProofExistSec2Eqt4}) that
\begin{equation}\label{ProofExistSec2Eqt5}
\begin{split}
&\left(1-\sum_{i=1}^p\int_M\vert u_i\vert^{2^\star}dv_g\right)^{2/2^\star}\\
&\le K_n^2 \left(\mu_{A,g}^p(M) - I_A({\mathcal U})\right) + o(1)\\
&\le K_n^2 \mu_{A,g}^p(M) \left(1 - \left(\sum_{i=1}^p\int_M\vert u_i\vert^{2^\star}dv_g\right)^{2/2^\star}\right) + o(1)
\end{split}
\end{equation}
for all $\alpha$, where ${\mathcal U}$ is the $p$-map whose components are the $u_i$'s, 
$o(1) \to 0$ as $\alpha \to +\infty$, and we used the equation
$$I_A({\mathcal U}) \ge \mu_{A,g}^p(M)\left(\sum_{i=1}^p\Vert u_i\Vert_{2^\star}^{2^\star}\right)^{2/2^\star}
\hskip.1cm .$$
Noting that
$$0 \le 1 - \left(\sum_{i=1}^p\int_M\vert u_i\vert^{2^\star}dv_g\right)^{2/2^\star} 
\le \left(1-\sum_{i=1}^p\int_M\vert u_i\vert^{2^\star}dv_g\right)^{2/2^\star}$$
for all $\alpha$, we get with (\ref{ProofExistSec2Eqt5}) that
\begin{equation}\label{ProofExistSec2Eqt6}
K_n^2 \mu_{A,g}^p(M) < 1\hskip.2cm\Rightarrow\hskip.2cm 
\sum_{i=1}^p\int_M\vert u_i\vert^{2^\star}dv_g = 1\hskip.1cm .
\end{equation}
In particular, ${\mathcal U} \in {\mathcal H}_M^p$, 
and then, by writing with (\ref{ProofExistSec2Eqt6}) that
\begin{eqnarray*}
&&\sum_{i=1}^p\left(\int_M\vert\nabla u_\alpha^i\vert^2dv_g - \int_M\vert\nabla u_i\vert^2dv_g\right)\\
&&= I_A({\mathcal U}_\alpha) - I_A({\mathcal U}) + o(1)\\
&&\le \mu_{A,g}^p(M) - \mu_{A,g}^p(M)\left(\sum_{i=1}^p\int_M\vert u_i\vert^{2^\star}dv_g\right)^{2/2^\star} 
+ o(1)\\
&& = o(1)\hskip.1cm ,
\end{eqnarray*}
we get with (\ref{ProofExistSec2Eqt1}) that
\begin{equation}\label{ProofExistSec2Eqt7}
\sum_{i=1}^p\int_M\vert\nabla(u_\alpha^i-u_i)\vert^2dv_g \to 0
\end{equation}
as $\alpha \to +\infty$. By (\ref{ProofExistSec2Eqt7}), $u_\alpha^i \to u_i$ strongly 
in $H_1^2(M)$ for all $i$, as $\alpha \to +\infty$. In particular, 
$I_A({\mathcal U}) = \mu_{A,g}^p(M)$, and 
since ${\mathcal U} \in {\mathcal H}_M^p$, ${\mathcal U} = (u_1,\dots,u_p)$ is a minimizer for 
$\mu_{A,g}^p(M)$. If $-A$ is cooperative, $I_A(\vert{\mathcal U}\vert) \le I_A({\mathcal U})$ where 
$\vert{\mathcal U}\vert = (\vert u_1\vert,\dots,\vert u_p\vert)$. In particular, if $-A$ is cooperative, 
up to replacing ${\mathcal U}$ by $\vert{\mathcal U}\vert$, we can choose ${\mathcal U}$ to be nonnegative. 
Independently, we clearly have that ${\mathcal U}$ is a weak solution of 
the system (\ref{GenericEqtThm1}) with $\Lambda = \mu_{A,g}^p(M)$, in the sense that 
for any $\Phi = (\varphi_1,\dots,\varphi_p)$ in $H_{1,p}^2(M)$, and any $i$,
\begin{equation}\label{ProofRegulEqt0}
\int_M(\nabla u_i\nabla\varphi_i)dv_g
+ \sum_{j=1}^p\int_MA_{ij}u_i\varphi_jdv_g = \Lambda \vert u_i\vert^{2^\star-2}u_i\varphi_i\hskip.1cm .
\end{equation}
Now we prove the regularity 
of ${\mathcal U}$, and more generally of any solution of (\ref{ProofRegulEqt0}). 
Here we follow the arguments developed by Trudinger \cite{Tru} for 
the Yamabe equation. By standard elliptic theory, it suffices to prove that the 
$u_i$'s are all in $L^q(M)$ for some $q > 2^\star$. 
For $\lambda > 0$ we let $F_\lambda$ and $G_\lambda$ be real 
functions defined by
\begin{equation}\label{TruProof1}
\left\{
\begin{array}{l}
F_\lambda(x) = \vert x\vert^{p_0}\hskip.2cm\hbox{if}\hskip.1cm \vert x\vert \le \lambda\\
 \\
F_\lambda(x) = p_0\lambda^{p_0-1}\vert x\vert - (p_0-1)\lambda^{p_0}\hskip.2cm\hbox{if not}
\end{array}
\right.
\end{equation}
and
\begin{equation}\label{TruProof2}
\left\{
\begin{array}{l}
G_\lambda(x) = \vert x\vert^{p_1}\hskip.2cm\hbox{if}\hskip.1cm \vert x\vert \le \lambda\\
 \\
G_\lambda(x) = p_0\lambda^{2(p_0-1)}\vert x\vert - (p_0-1)\lambda^{2p_0-1}\hskip.2cm\hbox{if not}\hskip.1cm ,
\end{array}
\right.
\end{equation}
where $2p_0 = 2^\star$ and $p_1 = 2^\star-1$ so that $2p_0 = p_1 + 1$. For 
$i = 1,\dots,p$, we let also
$$u_{i,\lambda}^1 = F_\lambda(u_i)\hskip.3cm\hbox{and}\hskip.3cm u_{i,\lambda}^2 = G_\lambda(u_i)\hskip.1cm .$$
Noting that $F_\lambda$ and $G_\lambda$ in (\ref{TruProof1}) and (\ref{TruProof2}) 
are Lipschitz functions, we get that the $u_{i,\lambda}^1$'s and $u_{i,\lambda}^2$'s are all in $H_1^2(M)$. 
Independently, we also have that for any $x \ge 0$,
\begin{equation}\label{TruProof4}
\begin{split}
&F_\lambda(x) \le x^{p_0}\hskip.1cm ,\hskip.1cm G_\lambda(x) \le x^{p_1}\hskip.1cm ,\hskip.1cm 
F_\lambda(x)^2 \ge xG_\lambda(x)\\
&\hbox{and}\hskip.1cm F_\lambda^\prime(x)^2 \le p_0G_\lambda^\prime(x)\hskip.1cm\hbox{if}\hskip.1cm x \not= \lambda\hskip.1cm .
\end{split}
\end{equation}
Taking $\varphi_i = u_{i,\lambda}^2$ in (\ref{ProofRegulEqt0}), since $u_i \in L^{2^\star}(M)$ 
for all $i$, and by (\ref{TruProof4}), we can write that
\begin{equation}\label{TruProof5}
\int_MG_\lambda^\prime(u_i)\vert\nabla u_i\vert^2dv_g \le C_1 + \vert\Lambda\vert 
\int_M\vert u_i\vert^{2^\star-1}G_\lambda(u_i)dv_g
\end{equation}
for all $i$, 
where $C_1 > 0$ does not depend on $\lambda$. Still by (\ref{TruProof4}) it follows from (\ref{TruProof5}) that
\begin{equation}\label{TruProof6}
\frac{1}{p_0}\int_M\vert\nabla u_{i,\lambda}^1\vert^2dv_g \le C_1 + 
\vert\Lambda\vert \int_M\vert u_i\vert^{2^\star-2}F_\lambda(u_i)^2dv_g\hskip.1cm .
\end{equation}
Given $K > 0$ we let
$$H_{i,K}^1 = \bigl\{x\hskip.1cm \hbox{s.t.}\hskip.1cm \vert u_i(x)\vert \le K\bigr\}
\hskip.3cm\hbox{and}\hskip.3cm H_{i,K}^2 = \bigl\{x\hskip.1cm \hbox{s.t.}\hskip.1cm \vert u_i(x)\vert \ge K\bigr\}\hskip.1cm .$$
By H\"older's inequality, and the Sobolev inequality for the embedding $H_1^2 \subset L^{2^\star}$, 
namely $\Vert\varphi\Vert_{2^\star} \le A \Vert\nabla\varphi\Vert_{H_1^2}$ for all $\varphi \in H_1^2(M)$, 
we write that
\begin{equation}\label{TruProof7}
\begin{split}
&\int_M\vert u_i\vert^{2^\star-2}F_\lambda(u_i)^2dv_g\\
&\le \int_{H_K^1}\vert u_i\vert^{2^\star-2}F_\lambda(u_i)^2dv_g 
+ \int_{H_K^2}\vert u_i\vert^{2^\star-2}F_\lambda(u_i)^2dv_g\\
&\le \int_{H_K^1}\vert u_i\vert^{2^\star-2}F_\lambda(u_i)^2dv_g 
+ \left(\int_{H_K^2}\vert u_i\vert^{2^\star}dv_g\right)^{2/n}\left(\int_MF_\lambda(u_i)^{2^\star}dv_g\right)^{2/2^\star}\\
&\le \int_{H_K^1}\vert u_i\vert^{2^\star-2}F_\lambda(u_i)^2dv_g 
+ \varepsilon_{i,K}^{2/n}A^2\int_M\left(\vert\nabla u_{i,\lambda}^1\vert^2 + (u_{i,\lambda}^1)^2\right)dv_g\hskip.1cm ,
\end{split}
\end{equation}
where, since $u_i \in L^{2^\star}(M)$,
 $\varepsilon_{i,K} = \int_{H_K^2}\vert u_i\vert^{2^\star}dv_g$
is such that
 $$\lim_{K\to+\infty}\varepsilon_{i,K} = 0\hskip.1cm .$$
We let $K > 0$ be such that 
$p_0A^2\varepsilon_{i,K}^{2/n} < 1$. For $\lambda > K$,
\begin{equation}\label{TruProof8}
\int_{H_K^1}\vert u_i\vert^{2^\star-2}F_\lambda(u_i)^2dv_g  \le K^{2(2^\star-1)}V_g\hskip.1cm ,
\end{equation}
where $V_g$ is the volume of $M$ with respect to $g$. 
Similarly, by (\ref{TruProof4}), and since 
$u_i \in L^{2^\star}(M)$, there exists $C_2 > 0$, independent of $\lambda$, such that
\begin{equation}\label{TruProof9}
\int_M(u_{i,\lambda}^1)^2dv_g \le C_2
\end{equation}
for all $i$. 
Inserting (\ref{TruProof7})--(\ref{TruProof9}) into (\ref{TruProof6}), and since $K$ is such that 
$p_0A^2\varepsilon_K^{2/n} < 1$, we get that
\begin{equation}\label{TruProof10}
\int_M\vert\nabla u_{i,\lambda}^1\vert^2dv_g \le C_3\hskip.1cm ,
\end{equation}
where $C_3 > 0$ is independent of $\lambda$. In particular, by (\ref{TruProof9}), (\ref{TruProof10}), and the
Sobolev inequality, there exists $C_4 > 0$ independent of $\lambda$ such that
$$\int_M(u_{i,\lambda}^1)^{2^\star}dv_g \le C_4$$
for all $i$. 
Letting $\lambda \to +\infty$, it clearly follows that the $u_i$'s are all  in $L^{2^\star p_0}(M)$. Since $p_0 > 1$, we then get that the 
$u_i$'s are all in $L^q(M)$ for some $q > 2^\star$. By standard elliptic theory, the $u_i$'s are in $C^{2,\theta}(M)$, 
$0 < \theta < 1$, for all $i$. This holds for ${\mathcal U}$ minimizing or not. 
If we assume that $-A$ is cooperative, see Lemma \ref{PosSolSec1} in Section \ref{Prel}, the 
maximum principle for functions applies. 
In particular, if ${\mathcal U}$ is minimizing, it can be chosen weakly positive since, 
as mentionned above, up to replacing ${\mathcal U}$ by $\vert{\mathcal U}\vert$, 
it can be chosen nonnegative. If the system is also 
fully coupled, see again Lemma \ref{PosSolSec1} in Section \ref{Prel}, ${\mathcal U}$ is necessarily strongly positive when it is weakly positive. 
This ends the proof of Theorem \ref{ExistThm}.
\end{proof}

In what follows we say that the operator $\Delta_g^p + A$ is {\it coercive} if 
there exists $K > 0$ such that 
$$I_A({\mathcal U}) \ge K\Vert{\mathcal U}\Vert_{H_{1,p}^2}^2$$
for all ${\mathcal U} \in H_{1,p}^2(M)$, where the norm in the right hand side 
is the norm of $H_{1,p}^2(M)$. Define $\lambda_A(g)$ to be the infimum of the $I_A({\mathcal U})$'s 
for ${\mathcal U} \in H_{1,p}^2(M)$, ${\mathcal U} = (u_1,\dots,u_p)$, such that $\sum_i\int_Mu_i^2dv_g = 1$. 
If $\lambda_A(g) > 0$, we can 
write that
\begin{eqnarray*} I_A({\mathcal U}) 
& \ge & \varepsilon I_A({\mathcal U}) + (1-\varepsilon)\lambda_A(g)\sum_{i=1}^p\int_Mu_i^2dv_g\\
& \ge & \varepsilon \Vert{\mathcal U}\Vert_{H_{1,p}^2}^2
\end{eqnarray*}
for all ${\mathcal U} \in H_{1,p}^2(M)$, 
where $\varepsilon >0$ is such that $(1-\varepsilon)\lambda_A(g)Id_p + \varepsilon A \ge \varepsilon Id_p$ 
in the sense of bilinear forms, and $Id_p$ is the $p\times p$-identity matrix. In particular, $\Delta_g^p + A$ is coercive. 
Conversely, since the $L^2$-norm is controlled 
by the $H_1^2$-norm, we easily get that $\lambda_A(g) > 0$ if $\Delta_g^p + A$ is coercive. This proves 
that $\Delta_g^p+A$ is coercive if and only if $\lambda_A(g) > 0$. 
Noting that $\mu_{A,g}^p(M) > 0$ if and only if $\lambda_A(g) > 0$, it follows that
$\mu_{A,g}^p(M) > 0$ if and only if $\Delta_g^p+A$ is coercive.

\medskip When $\mu_{A,g}^p(M) \le 0$, the existence of a minimizing solution for (\ref{GenericEqtThm1}) directly 
follows from Theorem \ref{ExistThm}. When $\mu_{A,g}^p(M) > 0$, and thus $\Delta_g^p+A$ is coercive, there 
are several situations where Theorem \ref{ExistThm} can be applied. For instance, we get with 
the developments in Aubin \cite{Aub} 
that if $n \ge 4$ and
$$A_{ii}(x) < \frac{n-2}{4(n-1)} S_g(x)$$
for some $i$ and some $x$, 
where $S_g$ is the scalar curvature of $g$, then $\mu_{A,g}^p(M) < K_n^{-2}$. We also get that 
$\mu_{A,g}^p(M) < K_n^{-2}$ if 
$\int_MA_{ii}dv_g < K_n^{-2}V_g^{2/2^\star}$
for some $i$, where $V_g$ is the volume of $M$ with 
respect to $g$, or, by Schoen \cite{Sch1}, if $A_{ii} \equiv \frac{n-2}{4(n-1)}S_g$ for some $i$, and the manifold is 
conformally distinct to the unit sphere. 
By Theorem \ref{ExistThm}, (\ref{GenericEqtThm1}) with $\Lambda = 1$ possesses a minimizing solution in such cases.

\medskip\noindent{\bf Remark 2.1:} 
The above examples involve test maps with no coupling. Namely 
test maps like ${\mathcal U} = 
(u_1,\dots,u_p)$ where, given some $i = 1,\dots,p$, we ask that $u_j \equiv 0$ for all $j \not= i$. There are of course several 
examples where coupling 
will help decreasing the value of $I_A$. Let for instance $(S^n,g_0)$ be the unit $n$-sphere, $p = 2$, and 
$A_\alpha$ be the matrix
\begin{equation*}
A_\alpha = \left(
\begin{matrix}
\frac{n(n-2)}{4} & \alpha\\
\alpha & \frac{n(n-2)}{4}
\end{matrix}\right)\hskip.1cm ,
\end{equation*}
where $\alpha: S^n \to {\mathbb R}$ is a smooth function. 
If $\alpha\equiv 0$, then, with the notations at the beginning of the section, $A_0 = A_{g_0}$ 
and we get that $\mu_{A_0,g_0}^2(S^n) = K_n^{-2}$. On the other hand, it is easily seen that 
$\mu_{A_\alpha,g_0}^2(S^n) < K_n^{-2}$ if $\alpha\not\equiv 0$. Indeed, if we let $u_0 > 0$ be a minimizer for the 
Yamabe invariant on the sphere with $\Vert u_0\Vert_{2^\star} = 1$, there is an entire family of such minimizers 
(including one constant function), 
and if we let ${\mathcal U}_\varepsilon = (u_0,-\varepsilon\alpha)$, where $\varepsilon > 0$ is small, then, for 
$I_A$ as in (\ref{FonctionalDef}), and $\Phi$ as in the definition of ${\mathcal H}^p_M$ in (\ref{MiniDefEqt}), 
$$I_{A_\alpha}({\mathcal U}_\varepsilon) 
\le \left(K_n^{-2} - 2\varepsilon\int_{S^n}\alpha^2u_0dv_{g_0} + Y_{S^n}(\alpha)\varepsilon^2\right)
\Phi({\mathcal U}_\varepsilon)^{2/2^\star}\hskip.1cm ,$$
where
$$Y_{S^n}(\alpha) = \int_{S^n}\left(\vert\nabla\alpha\vert^2+\frac{n(n-2)}{4}\alpha^2\right)dv_{g_0}
\hskip.1cm .$$
Taking $\varepsilon > 0$ sufficiently small such that 
$Y_{S^n}(\alpha)\varepsilon < 2\int_{S^n}\alpha^2u_0dv_{g_0}$, 
this proves the above claim 
that $\mu_{A_\alpha,g_0}^2(S^n) < K_n^{-2}$  when $\alpha\not\equiv 0$. 

\medskip\noindent{\bf Remark 2.2:} 
When $p = 1$, minimizers are always positive or negative (since if $u$ is a minimizer, then $\vert u\vert$ is also 
a minimizer and we can apply the maximum principle). When $p > 1$, this is not true anymore. However, 
in several situations, we are still dealing with positive solutions of systems like (\ref{GenericEqtThm1}). 
Suppose for instance that $p = 2$, and let $A$, $A^\prime$ be the matrix
\begin{equation}\label{ExPosMinSec2Matrix}
A = \left(
\begin{matrix}
\alpha & \beta\\
\beta & \gamma
\end{matrix}
\right)\hskip.2cm\hbox{and}\hskip.2cm
A^\prime = \left(
\begin{matrix}
\alpha & -\beta\\
-\beta & \gamma
\end{matrix}
\right)\hskip.1cm ,
\end{equation}
where $\alpha, \beta, \gamma$ are smooth functions in $M$, and $A$ is supposed to be such that
$\Delta_g^p+A$ is coercive so that $\mu(A) = \mu^2_{A,g}(M) > 0$. For ${\mathcal U} = (u,v)$ in $H_{1,2}^2(M)$ 
we let ${\mathcal U}^\prime$ be given by ${\mathcal U}^\prime = (u,-v)$. We suppose 
that $\beta \ge 0$, $\beta \not\equiv 0$, is nontrivial and nonnegative. If ${\mathcal U}_0 = (u_0,v_0)$  
is a minimizer for $\mu(A)$, then $I_A({\mathcal U}_0) \le I_A({\mathcal U}_0^\prime)$ and we 
get that $\int_M\beta u_0v_0dv_g \le 0$. In particular, since $\beta \ge 0$ and $\beta \not\equiv 0$, 
${\mathcal U}_0$ cannot be strongly positive. Pushing further the analysis, 
noting that $I_A({\mathcal U}) = I_{A^\prime}({\mathcal U}^\prime)$ for all ${\mathcal U} \in {\mathcal H}_M^p$, we 
easily get that ${\mathcal U}_0$ is a minimizer for $\mu(A)$ if and only if ${\mathcal U}_0^\prime$ is a minimizer 
for $\mu(A^\prime) = \mu_{A^\prime,g}^2(M)$. But $A^\prime$ is such that $-A^\prime$ is cooperative, so that $\vert{\mathcal U}_0^\prime\vert 
= (\vert u_0\vert,\vert v_0\vert)$ is also a minimizer for  $\mu(A^\prime)$. Since we assumed that $\beta\not\equiv 0$, the system 
is fully coupled. It follows that both $\vert u_0\vert$ and $\vert v_0\vert$ are positive functions. Summarizing, 
if $\beta \le 0$, and $\beta\not\equiv 0$, $-A$ is cooperative and the system is fully coupled. 
Minimizers for $\mu_{A,g}^2(M)$ are like ${\mathcal U}_0 = (u_0,v_0)$ or ${\mathcal U}_0 = (-u_0,-v_0)$, 
where $u_0$ and $v_0$ are positive functions in $M$. 
Then, up to a positive constant scale factor, ${\mathcal U} = (u_0,v_0)$ is a strongly positive solution of the original system
\begin{equation*}
\begin{cases}
\Delta_gu_0 + \alpha u_0 + \beta v_0 = u_0^{2^\star-1}\\
\Delta_gv_0 + \beta u_0 + \gamma v_0 = v_0^{2^\star-1}\hskip.1cm .
\end{cases}
\end{equation*}
Conversely, if $\beta \ge 0$, and $\beta\not\equiv 0$, 
minimizers for $\mu_{A,g}^2(M)$ are like ${\mathcal U}_0 = (u_0,-v_0)$ or 
${\mathcal U}_0 = (-u_0,v_0)$, where $u_0$ and $v_0$ are positive functions in $M$. They are not positive (neither negative). 
However  we are still dealing with 
strongly positive solutions of $2$-systems. In this case, up 
to a positive constant scale factor, ${\mathcal U} = (u_0,v_0)$ is a strongly positive solution of the modified system
\begin{equation*}
\begin{cases}
\Delta_gu_0 + \alpha u_0 + \beta^\prime v_0 = u_0^{2^\star-1}\\
\Delta_gv_0 + \beta^\prime u_0 + \gamma v_0 = v_0^{2^\star-1}\hskip.1cm ,
\end{cases}
\end{equation*}
where $\beta^\prime = -\beta$. As one can check, the above discussion, 
and the arguments we developed, extend to higher order 
systems.

\section{High energy solutions}\label{HighEn}

Let $(M,g)$ be a smooth compact  Riemannian manifold, $n \ge 3$. A preliminary remark is that 
there are several ways to construct positive solutions of $p$-systems from positive solutions of $q$-systems 
when $q < p$. Suppose for instance that we have two positive solutions of scalar equations like 
(\ref{YamabeTypeEqt}). Namely,
$\Delta_gu+hu = u^{2^\star-1}$ and 
$\Delta_gv+kv = v^{2^\star-1}$ in $M$, where $h$ and $k$ are two 
smooth functions in $M$. Then ${\mathcal U} = (u,v)$ is a strongly positive 
solution of the $2$-system
\begin{equation}\label{HighEnergEx1}
\begin{cases}
\Delta_gu + \alpha u + \beta v = u^{2^\star-1}\\
\Delta_gv + \beta u + \gamma v = v^{2^\star-1}
\end{cases}
\end{equation}
as soon as $\alpha$, $\beta$, and $\gamma$ are such that
\begin{equation}\label{HighEnergEx1StructCond}
\begin{split}
\bigl(h(x) - \alpha(x)\bigr) \frac{u(x)}{v(x)} 
& = \bigl(k(x) - \gamma(x)\bigr) \frac{v(x)}{u(x)}\\
& = \beta(x)
\end{split}
\end{equation}
for all $x \in M$. In other words, two positive solutions of scalar equations like 
(\ref{YamabeTypeEqt}) provide several examples (e.g like a one parameter family 
w.r.t $\beta$) of $2$-systems with strongly positive solutions. These solutions, see Remark 2.2 in Section \ref{ExSol}, are 
certainly not minimizing if $\beta > 0$.

\medskip In what follows 
 we discuss particular examples where we do get solutions with arbitrarily large energies (Proposition 
\ref{Prop1Sec3}), and multiple solutions with distinct energies (Proposition \ref{Prop2Sec3}). 
For $M$ a smooth compact Riemannian $(n-1)$-manifold, $n \ge 3$, and 
$S^1$ the circle in ${\mathbb R}^2$ of radius $1$ centered at $0$, we let $\hat M = S^1\times M$ 
and $\hat g$ be the standard product metric on $\hat M$. For $p \ge 1$, 
and $A: M \to M_p^s({\mathbb R})$, we consider the system
\begin{equation}\label{GenericEqtSec3}
\Delta_{\hat g}u_i + \sum_{j=1}^pA_{ij}(x)u_j = \vert u_i\vert^{2^\star-2}u_i
\end{equation}
in $\hat M$, for all $i = 1,\dots p$, where $(t,x)$, $t \in S^1$, $x \in M$,  is the variable in $\hat M$. 
For ${\mathcal U} = (u_1,\dots,u_p)$ a solution of (\ref{GenericEqtSec3}), we let 
$E({\mathcal U})$ be the energy of ${\mathcal U}$ given by $E({\mathcal U}) = 
\sum_{i=1}^p\Vert u_i\Vert_{2^\star}^{2^\star}$. 
We let also $\Lambda_{min} = K_n^{-n}$, where $K_n$ is as in 
(\ref{ValueSharpCst}), and let $Id_p$ be the $p\times p$-identity matrix. If ${\mathcal U} = (u_1,\dots,u_p)$ is 
a minimizer for $\mu_{A,g}^p(M)$, and $\mu_{A,g}^p(M) > 0$, where $\mu_{A,g}^p(M)$ is defined in (\ref{MiniDefEqt}), 
then $\tilde{\mathcal U} = (\tilde u_1,\dots,\tilde u_p)$ 
is a solution of (\ref{GenericEqtSec3}) when we let $\tilde u_i = \mu_{A,g}^p(M)^{(n-2)/4}u_i$ for all $i$. 
Its energy $E(\tilde{\mathcal U})$ is such that $E(\tilde{\mathcal U}) \le \Lambda_{min}$. 
We prove here that the following result holds.

\begin{prop}\label{Prop1Sec3} Let $(M,g)$ be a smooth compact Riemannian $(n-1)$-manifold, $n \ge 3$. For any 
$\Lambda > 0$ there exist positive real numbers 
$K_1(\Lambda) < K_2(\Lambda)$ such that if $A: M \to M_p^s({\mathbb R})$ satisfies that 
\begin{equation}\label{AssProp1Sec3}
K_1(\Lambda)Id_p < A(x) < K_2(\Lambda)Id_p
\end{equation}
for all $x$, in the sense of bilinear forms, then 
(\ref{GenericEqtSec3}) in $\hat M = S^1\times M$ possesses a solution ${\mathcal U}$ of energy 
$\Lambda \le E({\mathcal U}) \le \Lambda + \Lambda_{min}$. The solution can be chosen weakly positive 
if $-A$ is cooperative, and strongly positive if the system is also fully coupled.
\end{prop}

\begin{proof} For $m \ge 1$, we let $S^1(1/m)$ be the circle in ${\mathbb R}^2$ of radius $1/m$ centered at $0$. 
We let $M_m = S^1(1/m)\times M$ 
and $g_m$ be the standard product metric on $M_m$. We denote by $(t,x)$ the variable in $S^1(1/m)\times M$, and 
$G_m$ be the subgroup of $O(2)$ generated 
by $z \to e^{2i\pi/m}z$. We regard $G_m$ as acting on $\hat M = S^1\times M$ by $(t,x) \to (\sigma(t),x)$ for $\sigma \in G_m$. Then 
$\hat M/G_m = M_m$. We define $B_m = B_0(g_m)$ to be the smallest $B$ such that
$$\Vert u\Vert_{2^\star}^2 \le K_n^2\Vert\nabla u\Vert_2^2 + B\Vert u\Vert_2^2$$
for all $u \in H_1^2(M_m)$. For $t \in \left(0,B_m\right)$, and $t^\prime = K_n^{-2}t$, we let 
$\lambda_t = \mu_{t^\prime,g_m}^1(M_m)$ 
with the notations in Section \ref{ExSol}. It can be proved (see for instance 
Druet, Hebey and Robert \cite{DruHebRob2}) that $B_m$ exists, that $B_m \ge V_{g_m}^{-2/n}$, where $V_{g_m}$ is the volume 
of $M_m$ with respect to $g_m$, that $\lambda_t$ is nondecreasing in $T$, and that $\lambda_t \to K_n^{-2}$ 
as $t \to B_m$. For $\Lambda > 0$ we let $m = m_\Lambda$ be given by
$m = [\Lambda_{min}^{-1}\Lambda] + 1$, 
where, for $X >0$, $[ X]$ is the greatest integer not exceeding $X$. We let also $T(\Lambda) 
\in \left(0,B_m\right)$ be such that
\begin{equation}\label{ProofProp1Eqt1}
\Lambda < \left([\Lambda_{min}^{-1}\Lambda] + 1\right)\lambda_t^{\frac{n}{2}}
\end{equation}
for all $T(\Lambda) < t < B_m$. Then we define
$$K_1(\Lambda) = K_n^{-2}T(\Lambda)
\hskip.2cm\hbox{and}\hskip.2cm
K_2(\Lambda) = K_n^{-2}B_m\hskip.1cm .$$
Let $A: M \to M_p^s({\mathbb R})$ be such that 
(\ref{AssProp1Sec3}) holds for all $x \in M$, in the sense of bilinear forms. Since $A$ does 
not depend on the variable $t$, we may regard $A$ as defined in $M_m$. Moreover, since $M$ is compact, there exist 
$t_{min} < t_{max}$ in $\left(K_1(\Lambda),K_2(\Lambda)\right)$ such that $t_{min}Id_p \le A(x) \le t_{max}Id_p$ 
for all $x$, 
in the sense of bilinear forms. In particular, by (\ref{TrivialCouplingMin}), we can write that
\begin{equation}\label{ProofProp1Eqt2}
\mu^1_{t_{min},g_m}(M_m) \le \mu^p_{A,g_m}(M_m) \le 
\mu^1_{t_{max},g_m}(M_m)\hskip.1cm .
\end{equation}
Since $t_{max} < K_2(\Lambda)$, coming back to the 
very first definition of $B_m$, we do get that $\mu^1_{t_{max},g_m}(M_m) < K_n^{-2}$. Then, by Theorem 
\ref{ExistThm}, there exist ${\mathcal U} = (u_1,\dots,u_p)$ a solution in $M_m$ of the system
$$\Delta_{g_m}u_i + \sum_{j=1}^pA_{ij}(x)u_j = \mu^p_{A,g_m}(M_m)\vert u_i\vert^{2^\star-2}u_i$$
for all $i$, and such that $\sum_{i=1}^p\Vert u_i\Vert_{2^\star}^{2^\star} = 1$. 
The solution can be chosen weakly positive if $-A$ is cooperative, and 
strongly positive if the system is also fully coupled. We let $\tilde u_i$ be given by
$$\tilde u_i = \mu^p_{A,g_m}(M_m)^{\frac{n-2}{4}} u_i$$
and $\hat u_i$ be the function in $\hat M$ such that $\hat u_i/G_m = \tilde u_i$. Then $\hat{\mathcal U} 
= (\hat u_1,\dots,\hat u_p)$ is a solution of (\ref{GenericEqtSec3}), and by (\ref{ProofProp1Eqt1}) and 
(\ref{ProofProp1Eqt2}), its energy 
$$E(\hat{\mathcal U}) = m\mu_{A,g_m}^p(M_m)^{n/2}$$
is such that 
$\Lambda \le E(\hat{\mathcal U}) \le \Lambda + \Lambda_{min}$. This proves the proposition.
\end{proof}

Following ideas in Hebey and Vaugon \cite{HebVau} we may also prove existence of several solutions of distinct 
energies in particular cases. For $T > 0$, and $n \ge 3$, let $M_T$ be the manifold $M_T = S^1(T)\times S^{n-1}$ and 
$g_T$ be the standard product metric on $M_T$. For $p \ge 1$, 
and $A: S^{n-1} \to M_p^s({\mathbb R})$, we consider the system
\begin{equation}\label{GenericEqtSec3Mult}
\Delta_{g_T}u_i + \sum_{j=1}^pA_{ij}(x)u_j = \vert u_i\vert^{2^\star-2}u_i
\end{equation}
in $M_T$, for all $i = 1,\dots p$, where $(t,x)$, $t \in S^1(T)$, $x \in S^{n-1}$,  is the variable in $M_T$. 
We claim that the following multiplicity result holds. Such a result goes back to Schoen \cite{Sch3} 
(see also Hebey and Vaugon \cite{HebVau}) where it 
was proved when $p=1$ for the Yamabe equation.

\begin{prop}\label{Prop2Sec3} For any integer $k \ge 1$, and any $q > n/2$, there exists $T(k,q) > 0$ with 
the property that if $T \ge T(k,q)$ and $A: S^{n-1} \to M_p^s({\mathbb R})$ is such that
\begin{equation}\label{AssProp2Sec3}
\left(\frac{(n-2)^2}{4} - \frac{1}{T^q}\right)Id_p \le A(x) \le \frac{(n-2)^2}{4} Id_p
\end{equation}
for all $x$, in the sense of bilinear forms, then (\ref{GenericEqtSec3Mult}) possesses $k$ solutions of distinct energies in 
$M_T = S^1(T)\times S^{n-1}$. Moreover, these solutions can be chosen weakly positive 
if $-A$ is cooperative, and strongly positive if the system is also fully coupled.
\end{prop}

\begin{proof} Let $t > 0$. We know from Hebey and Vaugon \cite{HebVau} that
\begin{equation}\label{SobEst}
\frac{1}{K_n^2}\Vert u\Vert_{2^\star}^2 \le \Vert\nabla u\Vert_2^2 + \left(\frac{(n-2)^2}{4} + \frac{1}{4t^2}\right)\Vert u\Vert_2^2
\end{equation}
for all $u \in H_1^2(M_t)$. We fix $k \ge 1$ integer, and let $G_\alpha$, $\alpha = 1,\dots,k$, be $k$ groups of order 
$\alpha$ like the $G_m$'s in the proof of Proposition \ref{Prop1Sec3}. We let $G_\alpha$ act on $M_T$, $T > 0$, by 
$(t,x) \to (\sigma(t),x)$ for $\sigma \in G_\alpha$. Then $M_T/G_\alpha = M_{T/\alpha}$. Noting that $A$ depends only 
on the variable in $S^{n-1}$, we can regard $A$ as defined on all the $M_{T/\alpha}$'s. 
We assume that $A$ satisfies (\ref{AssProp2Sec3}) with $T = t$, $t > 0$, and let $\theta_n = (n-2)^2/4$. 
Then $\theta_n = \frac{n-2}{4(n-1)} S_{g_t}$ for all $t > 0$, where $S_{g_t}$ is the scalar curvature of $g_t$, and it 
follows from the right inequality in (\ref{AssProp2Sec3}), 
Theorem \ref{ExistThm}, and the resolution 
of the Yamabe problem, that there exists a minimizer ${\mathcal U}_t$ for 
$\mu_{A,g_t}^p(M_t)$. Moreover, still by Theorem \ref{ExistThm}, ${\mathcal U}_t$ can be chosen weakly positive 
if $-A$ is cooperative, and strongly positive if the system is also fully coupled.
For $u \in L^{2^\star}(M_t)$, of norm $1$ in $L^{2^\star}(M_t)$, we can write with 
H\"older's inequality that $\Vert u\Vert_2 \le (2\pi t\omega_{n-1})^{1/n}$, where $\omega_{n-1}$ 
is the volume of the unit $(n-1)$-sphere. 
By (\ref{TrivialCouplingMin}), (\ref{AssProp2Sec3}), and (\ref{SobEst}) we then get that
\begin{equation}\label{ProofProp2Sec3Eqt2}
\begin{split}
&\left(1-F_q(t)\right)K_n^{-2} \le \mu_{A,g_t}^p(M_t) \le K_n^{-2}\hskip.1cm,\hskip.1cm\hbox{where}\\
&F_q(t) = (2\pi\omega_{n-1})^{2/n}K_n^2\left(\frac{1}{t^q} + \frac{1}{4t^2}\right) t^{2/n}\hskip.1cm .
\end{split}
\end{equation}
Given $T > 0$, we let ${\mathcal U}_{T,\alpha}$ be the minimizer ${\mathcal U}_t$ for $t = T/\alpha$. Then we let 
$\hat{\mathcal U}_{T,\alpha}$ be the map on $M_T$ such that $\hat{\mathcal U}_{T,\alpha}/G_\alpha = 
{\mathcal U}_{T,\alpha}$, and 
$${\mathcal W}_{T,\alpha} = \mu_{A,g_{T/\alpha}}^p(M_{T/\alpha})^{\frac{n-2}{4}}\hat{\mathcal U}_{T,\alpha}
\hskip.1cm .$$
It is easily checked that ${\mathcal W}_{T,\alpha}$ is a solution of (\ref{GenericEqtSec3Mult}) in $M_T$ for all $\alpha$, 
and that the energy $E({\mathcal W}_{T,\alpha})$ is such that
\begin{equation}\label{ProofProp2Sec3Eqt3}
E({\mathcal W}_{T,\alpha})^{2/n} = \alpha^{2/n} \mu_{A,g_{T/\alpha}}^p(M_{T/\alpha})
\hskip.1cm .
\end{equation}
In particular, it follows from (\ref{ProofProp2Sec3Eqt2}) and (\ref{ProofProp2Sec3Eqt3}) that 
$E({\mathcal W}_{T,\alpha-1})^{2/n} < E({\mathcal W}_{T,\alpha})^{2/n}$ for all $\alpha = 2,\dots,k$ if
$F_q(T/\alpha) < \left(\alpha^{2/n}-(\alpha-1)^{2/n}\right)\alpha^{-2/n}$ for all $\alpha = 2,\dots,k$. 
Since $q > 2/n$, such inequalities are satisfied for $T \ge T(k,q)$ sufficiently large, and we get $k$ solutions 
${\mathcal W}_{T,\alpha}$ of (\ref{GenericEqtSec3Mult}) in $M_T$ with distinct energies. This proves the proposition.
\end{proof}

\section{The $H_1^2$-theory for blow-up}\label{SobTheory}

In what follows we let $(M,g)$ be a smooth compact Riemannian manifold of dimension $n \ge 3$,  $p \ge 1$ integer, 
and $\left(A(\alpha)\right)_\alpha$, $\alpha\in{\mathbb N}$, be a sequence of smooth maps 
$A(\alpha): M \to M_p^s({\mathbb R})$. We let $A(\alpha) = (A_{ij}^\alpha)$, and consider systems like
\begin{equation}\label{GenericEqtSec4}
\Delta_gu_i + \sum_{j=1}^pA_{ij}^\alpha(x)u_j = \vert u_i\vert^{2^\star-2}u_i
\end{equation}
in $M$, for all $i = 1,\dots p$. We assume in what follows that the $A(\alpha)$'s satisfy that 
there exists $C > 0$, independent of $\alpha$, and a continuous map $A: M \to M_p^s({\mathbb R})$, 
$A = (A_{ij})$, such that
\begin{equation}\label{CondSec4}
\begin{split}
&\vert\sum_{i,j}A^\alpha_{ij}(x)X_iX_j\vert \le C \sum_i(X_i)^2\hskip.1cm ,\hskip.1cm\hbox{and}\\
&A^\alpha_{ij} \to A_{ij}\hskip.1cm\hbox{in}\hskip.1cm L^2(M)\hskip.1cm\hbox{as}\hskip.1cm \alpha \to +\infty\hskip.1cm ,
\end{split}
\end{equation}
where the first equation in (\ref{CondSec4}) should hold 
for all $\alpha$, $x\in M$, and $X = (X_1,\dots,X_p)$ in ${\mathbb R}^p$, and the second equation in (\ref{CondSec4}) 
should hold for all $i, j$.
We denote by $I_{\alpha,p}$ the functional defined for ${\mathcal U} = (u_1,\dots,u_p)$ 
in $H_{1,p}^2(M)$ by
\begin{equation}\label{DefFuncSec4}
\begin{split}
I_{\alpha,p}({\mathcal U})
&= \frac{1}{2}\sum_{i=1}^p\int_M\vert\nabla u_i\vert^2dv_g 
+ \frac{1}{2}Ê\sum_{i,j=1}^p\int_MA^\alpha_{ij}u_iu_jdv_g\\
&\hskip.2cm -\frac{1}{2^\star}\sum_{i=1}^p\int_M\vert u_i\vert^{2^\star}dv_g
\end{split}
\end{equation}
We say that a sequence $({\mathcal U}_\alpha)_\alpha$ in $H_{1,p}^2(M)$ is a {\it Palais-Smale 
sequence} for (\ref{GenericEqtSec4}), or 
for $I_{\alpha,p}$, if the $I_{\alpha,p}({\mathcal U}_\alpha)$'s are bounded with respect to $\alpha$, and 
$DI_{\alpha,p}({\mathcal U}_\alpha) \to 0$ in $H_{1,p}^2(M)^\prime$ as $\alpha \to +\infty$. The 
Palais-Smale sequence is said to be nonnegative if the components of ${\mathcal U}_\alpha$ are nonnegative 
functions. For $(x_\alpha)_\alpha$ a converging sequence 
of points in $M$, and $(\mu_\alpha)_\alpha$ a sequence of positive real numbers converging to zero, we define a 
{\it $1$-bubble} as a sequence $(B_\alpha)_\alpha$ of functions in $M$ given by
\begin{equation}\label{Def1BubbleSec4}
B_\alpha(x) = \left(\frac{\mu_\alpha}{\mu_\alpha^2 + \frac{d_g(x_\alpha,x)^2}{n(n-2)}}\right)^{\frac{n-2}{2}}\hskip.1cm .
\end{equation}
The $x_\alpha$'s are referred to as the {\it centers} and the $\mu_\alpha$'s as 
the {\it weights} of the $1$-bubble $(B_\alpha)_\alpha$. We define a {\it $p$-bubble} as a 
sequence $({\mathcal B}_\alpha)_\alpha$ 
of $p$-maps such that, if we write that ${\mathcal B}_\alpha = (B_\alpha^1,\dots,B_\alpha^p)$, then 
$(B_\alpha^i)_\alpha$ is a $1$-bubble for exactly one $i$, and for $j \not= i$, 
$(B_\alpha^j)_\alpha$ is the trivial zero sequence. In other words, a $p$-bubble is a sequence 
of $p$-maps such that one of the components of the sequence is a $1$-bubble, and the other 
components are trivial zero sequences. 
One remark with respect to the definition 
(\ref{Def1BubbleSec4}) is that if $u: {\mathbb R}^n \to {\mathbb R}$ is given by
\begin{equation}\label{PosSolCritEuclEqt}
u(x) = \left(1 + \frac{\vert x\vert^2}{n(n-2)}\right)^{-\frac{n-2}{2}}\hskip.1cm ,
\end{equation}
then $u$ is a positive solution of the critical Euclidean equation $\Delta u = u^{2^\star-1}$, where 
$\Delta = -\sum\partial^2/\partial x_i^2$. More precisely, $u$ is 
the only positive solution of the equation in ${\mathbb R}^n$ 
which is such that $u(0) = 1$ and $u$ is maximum at $0$. All the 
other positive solutions of the equation $\Delta u = u^{2^\star-1}$ in ${\mathbb R}^n$, see 
Caffarelli, Gidas and Spruck \cite{CafGidSpr} and 
Obata \cite{Oba}, 
are then given by
$\tilde u(x) = \lambda^{(n-2)/2}u\left(\lambda(x-a)\right)$, 
where $\lambda > 0$ and $a \in {\mathbb R}^n$. Another remark with respect to the definition 
(\ref{Def1BubbleSec4}) is that the $B_\alpha$'s in (\ref{Def1BubbleSec4}) live on balls centered at $x_\alpha$ 
and radii of the order of $\sqrt{\mu_\alpha}$. 
Indeed, an equation like $B_\alpha(x) \ge C$ is equivalent to 
$d_g(x_\alpha,x) \le C^\prime\sqrt{\mu_\alpha}$, where $C, C^\prime > 0$, 
and we can write 
that for any $R > 0$,
\begin{equation}\label{C0RangeInteractSec4}
\begin{split}
&\inf_{x \in B_{x_\alpha}(R\sqrt{\mu_\alpha})}B_\alpha(x) = \left(\frac{n(n-2)}{R^2}\right)^{\frac{n-2}{2}} 
+ \varepsilon_\alpha\hskip.1cm ,\hskip.1cm\hbox{and}\\
&\sup_{M\backslash B_{x_\alpha}(R\sqrt{\mu_\alpha})}B_\alpha(x) = \left(\frac{n(n-2)}{R^2}\right)^{\frac{n-2}{2}} 
+ \varepsilon_\alpha\hskip.1cm ,
\end{split}
\end{equation}
where the sequence $(\varepsilon_\alpha)$ is such that $\varepsilon_\alpha \to 0$ as $\alpha \to +\infty$. In particular, the $B_\alpha$'s 
converge to $0$ in $M\backslash B_{x_\alpha}(R_\alpha\sqrt{\mu_\alpha})$ if $R_\alpha \to +\infty$ as $\alpha \to +\infty$. 
On the other hand (see, for instance Druet and Hebey \cite{DruHeb2}), the $B_\alpha$'s in (\ref{Def1BubbleSec4}) 
consume their $H_1^2$-norm on balls centered at $x_\alpha$ 
and radii of the order of $\mu_\alpha$. We sometimes 
refer to $\sqrt{\mu_\alpha}$ as the $C^0$-range of interaction 
of $(B_\alpha)_\alpha$, and to $\mu_\alpha$ as the $H_1^2$-range of interaction of $(B_\alpha)_\alpha$.

\medskip Since we assumed that (\ref{CondSec4}) holds, 
there is a limit system for (\ref{GenericEqtSec4}). The limit system writes as
\begin{equation}\label{LimitSystSec4}
\Delta_gu_i + \sum_{j=1}^pA_{ij}(x)u_j = \vert u_i\vert^{2^\star-2}u_i
\end{equation}
in $M$, for all $i = 1,\dots p$, where the $A_{ij}$'s are as in (\ref{CondSec4}). We let $I_{\infty,p}$ be the functional 
defined for ${\mathcal U} = (u_1,\dots,u_p)$ 
in $H_{1,p}^2(M)$ by
\begin{equation}\label{DefFuncBisSec4}
\begin{split}
I_{\infty,p}({\mathcal U})
&= \frac{1}{2}\sum_{i=1}^p\int_M\vert\nabla u_i\vert^2dv_g 
+ \frac{1}{2}Ê\sum_{i,j=1}^p\int_MA_{ij}u_iu_jdv_g\\
&\hskip.2cm -\frac{1}{2^\star}\sum_{i=1}^p\int_M\vert u_i\vert^{2^\star}dv_g\hskip.1cm .
\end{split}
\end{equation}
We prove in this section that the following result holds. The case $p = 1$ goes back to Struwe \cite{Str}, 
with related references by Br\'ezis and Coron \cite{BreCor}, 
Lions \cite{Lio}, Sacks-Uhlenbeck \cite{SacUhl}, Schoen \cite{Sch3}, and Wente \cite{Wen}. 
Theorem \ref{SobTheorThm} fully answers the question of describing the asymptotic behaviour of 
Palais-Smale sequences for (\ref{GenericEqtSec4}) in terms of Sobolev spaces.

\begin{thm}[$H_1^2$-Theory]\label{SobTheorThm} Let $(M,g)$ be a smooth compact Riemannian manifold of dimension 
$n \ge 3$, $p \ge 1$, and $\left(A(\alpha)\right)_\alpha$, $A(\alpha): M \to M_p^s({\mathbb R})$ , a sequence of smooth maps 
satisfying (\ref{CondSec4}). Let also $({\mathcal U}_\alpha)_\alpha$ be a nonnegative 
Palais-Smale sequence for (\ref{GenericEqtSec4}). Then there exists a nonnegative solution ${\mathcal U}^0$ of the 
limit system (\ref{LimitSystSec4}), there exists $k \in {\mathbb N}$, and there exist $p$-bubbles 
$({\mathcal B}_{j,\alpha})_\alpha$, $j = 1,\dots,k$, such that, up to a subsequence,
\begin{equation}\label{EqtThmSobTherSec4}
{\mathcal U}_\alpha = {\mathcal U}^0 + \sum_{j=1}^k{\mathcal B}_{j,\alpha} + {\mathcal R}_\alpha
\end{equation}
for all $\alpha$, where $({\mathcal R}_\alpha)_\alpha$ is a sequence in $H_{1,p}^2(M)$ converging strongly to $0$ 
in $H_{1,p}^2(M)$ as $\alpha \to +\infty$.
\end{thm}

Let ${\mathcal U}_\alpha = (u_\alpha^1,\dots,u_\alpha^p)$, ${\mathcal U}^0 = (u^0_1,\dots,u^0_p)$, and 
${\mathcal R}_\alpha = (R_\alpha^1,\dots,R_\alpha^p)$. 
Another way we can write the $H_1^2$-decomposition (\ref{EqtThmSobTherSec4}) is that for any $i$,
\begin{equation}\label{EqtThmSobTherSec4Bis}
u_\alpha^i = u^0_i + \sum_{j=1}^{k_i}B_{j,\alpha}^i + R_\alpha^i\hskip.1cm ,
\end{equation}
where the $k_i$'s are nonnegative integers such that $\sum_{i=1}^pk_i = k$, 
possibly $k_i = 0$ for some $i=1,\dots,p$, and the $(B_{j,\alpha}^i)_\alpha$'s 
are the $1$-bubbles from which the 
$p$-bubbles $({\mathcal B}_{j,\alpha})_\alpha$'s in (\ref{EqtThmSobTherSec4}) are defined.

\medskip In addition to (\ref{EqtThmSobTherSec4}), we also have that the energy of the ${\mathcal U}_\alpha$'s split. Namely 
we can write that for any $\alpha$,
\begin{equation}\label{EnerSplitTheorSob}
I_{\alpha,p}({\mathcal U}_\alpha) = I_{\infty,p}({\mathcal U}^0) + \frac{k}{n}K_n^{-n} + o(1)\hskip.1cm ,
\end{equation} 
where $K_n$ is as in (\ref{ValueSharpCst}), $I_{\alpha,p}$ is as in (\ref{DefFuncSec4}), $I_{\infty,p}$ 
is as in (\ref{DefFuncBisSec4}), and $o(1) \to 0$ as $\alpha \to +\infty$. Note that for  
$({\mathcal B}_\alpha)_\alpha$ a $p$-bubble, $I_{\alpha,p}({\mathcal B}_\alpha) = n^{-1}K_n^{-n} + o(1)$. An independent 
remark is that if $-A$ is cooperative, then ${\mathcal U}^0$ is weakly positive, and if $-A$ is cooperative and 
the limit system (\ref{LimitSystSec4}) is fully coupled, then either ${\mathcal U}^0$ is zero, or it is 
strongly positive. 

\medskip Let $\eta: {\mathbb R}^n \to {\mathbb R}$, $0 \le \eta \le 1$, be a smooth cutoff function such that 
$\eta = 1$ in $B_0(\delta)$, and $\eta = 0$ in ${\mathbb R}^n\backslash B_0(2\delta)$, for $\delta > 0$ small. 
In what follows we say that a sequence of functions $(\hat B_\alpha)_\alpha$, $\hat B_\alpha: M \to {\mathbb R}$,  is a 
generalized $1$-bubble if we can write that
\begin{equation}\label{DefGenBubbleTheorSob}
\hat B_\alpha(x) = (R_\alpha)^{\frac{n-2}{2}}\eta_\alpha(x)
u\left(R_\alpha\exp_{x_\alpha}^{-1}(x)\right)\hskip.1cm,
\end{equation}
where $u\not\equiv 0$ is a solution in ${\mathbb R}^n$ of the critical Euclidean equation $\Delta u = \vert u\vert^{2^\star-2}u$,  
$(R_\alpha)_\alpha$ is a sequence of positive real numbers such that $R_\alpha \to +\infty$ as $\alpha \to +\infty$, 
$(x_\alpha)_\alpha$ is a converging sequence of points in $M$, and $\eta_\alpha = \eta\circ\exp_{x_\alpha}^{-1}$. 
The $x_\alpha$'s are referred to 
as the centers and the $R_\alpha$'s as the weights of $(\hat B_\alpha)_\alpha$. 
If $u$ is positive, see 
for instance Druet and Hebey \cite{DruHeb2}, generalized $1$-bubbles are like 
$1$-bubbles in the sense that for any $(\hat B_\alpha)_\alpha$ there exists $(B_\alpha)_\alpha$ such that 
$\hat B_\alpha = B_\alpha + R_\alpha$ for all $\alpha$, where $R_\alpha \to 0$ in $H_1^2(M)$ as $\alpha \to +\infty$. 

\begin{proof} We prove Theorem \ref{SobTheorThm} by coming back to the well understood $p = 1$ case. 
Let $({\mathcal U}_\alpha)_\alpha$ be a nonnegative 
Palais-Smale sequence for (\ref{GenericEqtSec4}). 
As above, we write that ${\mathcal U}_\alpha = (u_\alpha^1,\dots,u_\alpha^p)$. 
First we claim that
\begin{equation}\label{ProofTheorSobAssert1}
\hbox{the sequence}\hskip.1cm ({\mathcal U}_\alpha)_\alpha\hskip.1cm\hbox{is bounded in}
\hskip.1cm H_{1,p}^2(M)\hskip.1cm .
\end{equation} 
In order to prove (\ref{ProofTheorSobAssert1}), we start noting that the equation
\begin{eqnarray*} 
DI_{\alpha,p}({\mathcal U}_\alpha).{\mathcal U}_\alpha
& = & \sum_{i=1}^p\int_M\vert\nabla u_\alpha^i\vert^2dv_g 
+ \sum_{i=1}^p\int_MA_{ij}^\alpha u_\alpha^iu_\alpha^jdv_g\\
&&\hskip.4cm - \sum_{i=1}^p\int_M(u_\alpha^i)^{2^\star}dv_g\hskip.1cm ,
\end{eqnarray*}
and the equation $DI_{\alpha,p}({\mathcal U}_\alpha).{\mathcal U}_\alpha = 
o\left(\Vert{\mathcal U}_\alpha\Vert_{H_{1,p}^2}\right)$, give that
$$I_{\alpha,p}({\mathcal U}_\alpha) = \frac{1}{n}\sum_{i=1}^p\int_M(u_\alpha^i)^{2^\star}dv_g  
+ o\left(\Vert{\mathcal U}_\alpha\Vert_{H_{1,p}^2}\right)\hskip.1cm .$$
By the definition of a Palais-Smale sequence, $\vert I_{\alpha,p}({\mathcal U}_\alpha)\vert \le C$, where 
$C > 0$ is independent of $\alpha$. Hence, 
\begin{equation}\label{ProofClaim1SobTheorEqt1}
\sum_{i=1}^p\int_M(u_\alpha^i)^{2^\star}dv_g \le C + 
o\left(\Vert{\mathcal U}_\alpha\Vert_{H_{1,p}^2}\right)\hskip.1cm ,
\end{equation}
where $C > 0$ is independent of $\alpha$, and by H\"older's inequality, we can also write that 
\begin{equation}\label{ProofClaim1SobTheorEqt2}
\sum_{i=1}^p\int_M(u_\alpha^i)^2dv_g \le C + 
o\left(\Vert{\mathcal U}_\alpha\Vert_{H_{1,p}^2}^{2/2^\star}\right)\hskip.1cm ,
\end{equation}
where $C > 0$ is independent of $\alpha$. By (\ref{ProofClaim1SobTheorEqt1}),
\begin{equation}\label{ProofClaim1SobTheorEqt3}
\begin{split}
&\sum_{i=1}^p\int_M\vert\nabla u_\alpha^i\vert^2dv_g 
+ \sum_{i=1}^p\int_MA_{ij}^\alpha u_\alpha^iu_\alpha^jdv_g\\
&= 2I_{\alpha,p}({\mathcal U}_\alpha) + \frac{2}{2^\star} 
\sum_{i=1}^p\int_M(u_\alpha^i)^{2^\star}dv_g\\
&\le C + 
o\left(\Vert{\mathcal U}_\alpha\Vert_{H_{1,p}^2}\right)\hskip.1cm ,
\end{split}
\end{equation}
where $C > 0$ is independent of $\alpha$. By (\ref{CondSec4}) we can write that
\begin{eqnarray*}
\Vert{\mathcal U}_\alpha\Vert_{H_{1,p}^2}^2
&=& \sum_{i=1}^p\int_M\vert\nabla u_\alpha^i\vert^2dv_g 
+ \sum_{i=1}^p\int_M(u_\alpha^i)^2dv_g\\
&\le&\sum_{i=1}^p\int_M\vert\nabla u_\alpha^i\vert^2dv_g 
+ \sum_{i=1}^p\int_MA_{ij}^\alpha u_\alpha^iu_\alpha^jdv_g 
+ C\sum_{i=1}^p\int_M(u_\alpha^i)^2dv_g
\end{eqnarray*}
and we easily get that the ${\mathcal U}_\alpha$'s are bounded in $H_{1,p}^2(M)$ by combining the above  
equation with (\ref{ProofClaim1SobTheorEqt2}) and (\ref{ProofClaim1SobTheorEqt3}). 
This proves (\ref{ProofTheorSobAssert1}). Now we let 
$I_p$ be the decoupled functional defined for ${\mathcal U} = (u_1,\dots,u_p)$ 
in $H_{1,p}^2(M)$ by
\begin{equation}\label{DefDecFuncProofTheorSob}
I_p({\mathcal U})
= \frac{1}{2}\sum_{i=1}^p\int_M\vert\nabla u_i\vert^2dv_g 
-\frac{1}{2^\star}\sum_{i=1}^p\int_M\vert u_i\vert^{2^\star}dv_g\hskip.1cm .
\end{equation}
By (\ref{ProofTheorSobAssert1}), up to passing to a subsequence, we may assume that 
${\mathcal U}_\alpha \rightharpoonup {\mathcal U}^0$ weakly in $H_{1,p}^2(M)$,  
that ${\mathcal U}_\alpha \to {\mathcal U}^0$ strongly in $L^2_p(M)$, and that ${\mathcal U}_\alpha \to {\mathcal U}^0$ 
almost everywhere, where  ${\mathcal U}^0 = (u^0_1,\dots,u^0_p)$ 
is some map in $H_{1,p}^2(M)$, the convergences have to be understood as $\alpha \to +\infty$, and the convergence in 
$L^2_p(M)$ means that $u_\alpha^i \to u^0_i$ in $L^2(M)$ for all $i$. In particular, 
${\mathcal U}^0$ is a nonnegative $p$-map. We let 
$\hat{\mathcal U}_\alpha = {\mathcal U}_\alpha - {\mathcal U}^0$, and $\hat u_\alpha^i = u_\alpha^i - u^0_i$ 
for $i = 1,\dots,p$. Now we claim that
\begin{equation}\label{ProofTheorSobAssert2}
\begin{split}
&{\mathcal U}^0\hskip.1cm\hbox{is a nonnegative solution of the limit system (\ref{LimitSystSec4})}\hskip.1cm ,\\
&(\hat{\mathcal U}_\alpha)_\alpha\hskip.1cm\hbox{is a Palais-Smale sequence for}\hskip.1cm I_p
\hskip.1cm ,\hskip.1cm\hbox{and}\\
&I_p(\hat{\mathcal U}_\alpha) = I_{\alpha,p}({\mathcal U}_\alpha) - I_{\infty,p}({\mathcal U}^0) + o(1)
\hskip.1cm\hbox{for all}\hskip.1cm \alpha\hskip.1cm ,
\end{split}
\end{equation} 
where $o(1) \to 0$ as $\alpha \to +\infty$. In order to prove (\ref{ProofTheorSobAssert2}), we first observe that 
if $\Phi = (\varphi_1,\dots,\varphi_p)$, then
\begin{equation}\label{ProofClaim2SobTheorEqt1}
\begin{split} 
DI_{\alpha,p}({\mathcal U}_\alpha).\Phi
&=   \sum_{i=1}^p\int_M(\nabla u_\alpha^i\nabla\varphi_i)dv_g 
+ \sum_{i,j=1}^p\int_MA^\alpha_{ij}u_\alpha^i\varphi_jdv_g\\
&\hskip.4cm - \sum_{i=1}^p\int_M(u_\alpha^i)^{2^\star-1}\varphi_idv_g\\
&= o(1)\hskip.1cm ,
\end{split}
\end{equation}
where $o(1) \to 0$ as $\alpha \to +\infty$. By (\ref{CondSec4}), $A^\alpha_{ij} \to A_{ij}$ in $L^q(M)$ 
for all $q \ge 1$ and all $i, j$, as $\alpha \to +\infty$. Then, by H\"older's inequality, and (\ref{ProofTheorSobAssert1}), 
we can write that
\begin{equation}\label{ProofClaim2SobTheorEqt2}
\begin{split}
\int_M\left\vert A_{ij}^\alpha - A_{ij}\right\vert u_\alpha^i\vert\varphi_j\vert dv_g 
&\le \Vert u_\alpha^i\Vert_{2^\star}\left\Vert A^\alpha_{ij}-A_{ij}\right\Vert_n\Vert\varphi_j\Vert_2\\
&= o(1)\hskip.1cm ,
\end{split}
\end{equation}
where $o(1) \to 0$ as $\alpha \to +\infty$. Combining (\ref{ProofClaim2SobTheorEqt1}) and 
(\ref{ProofClaim2SobTheorEqt2}), passing to the limit as $\alpha \to +\infty$, 
it easily follows that
$$\sum_{i=1}^p\int_M(\nabla u_i^0\nabla\varphi_i)dv_g 
+ \sum_{i,j=1}^p\int_MA_{ij}u_i^0\varphi_jdv_g
= \sum_{i=1}^p\int_M(u_i^0)^{2^\star-1}\varphi_idv_g\hskip.1cm .$$
Since $\Phi$ is arbitrary, this proves that 
${\mathcal U}^0$ is a nonnegative solution of the limit system (\ref{LimitSystSec4}). Now we 
compute the energy of $\hat{\mathcal U}_\alpha$. We write that
\begin{equation}\label{ProofClaim2SobTheorEqt3}
\begin{split}
&\int_MA^\alpha_{ij}u_\alpha^iu_\alpha^jdv_g 
- \int_MA_{ij}u_i^0u_j^0dv_g\\
&= \int_M\left(A^\alpha_{ij}-A_{ij}\right)u_\alpha^iu_\alpha^jdv_g 
+ \int_MA_{ij}(u_\alpha^i-u^0_i)u_\alpha^jdv_g\\
&\hskip.4cm + \int_MA_{ij}u^0_i(u_\alpha^j-u^0_j)dv_g
\end{split}
\end{equation}
for all $\alpha$, and all $i$, $j$. By (\ref{CondSec4}), as already mentionned, $A^\alpha_{ij} \to A_{ij}$ in $L^q(M)$ 
for all $q \ge 1$ and all $i, j$. In particular, 
$A^\alpha_{ij} \to A_{ij}$ in $L^{n/2}(M)$, and by (\ref{CondSec4}) and (\ref{ProofTheorSobAssert1}), we get with  
(\ref{ProofClaim2SobTheorEqt3}) that
\begin{equation}\label{ProofClaim2SobTheorEqt4}
\sum_{ij=1}^p\int_MA^\alpha_{ij}u_\alpha^iu_\alpha^jdv_g 
= \sum_{ij=1}^p\int_MA_{ij}u_i^0u_j^0dv_g + o(1)\hskip.1cm ,
\end{equation}
where $o(1) \to 0$ as $\alpha \to +\infty$. Since $u_\alpha^i \rightharpoonup u^0_i$ in $H_1^2$, 
\begin{equation}\label{ProofClaim2SobTheorEqt5}
\int_M\vert\nabla u_\alpha^i\vert^2dv_g = 
\int_M\vert\nabla u^0_i\vert^2dv_g + \int_M\vert\nabla\hat u_\alpha^i\vert^2dv_g + o(1)
\end{equation}
for all $\alpha$ and $i$. By (\ref{ProofClaim2SobTheorEqt4}) and (\ref{ProofClaim2SobTheorEqt5}), 
we then get that
\begin{equation}\label{ProofClaim2SobTheorEqt6}
I_{\alpha,p}({\mathcal U}_\alpha) = I_{\infty,p}({\mathcal U}^0) 
+ I_p(\hat{\mathcal U}_\alpha) - \frac{1}{2^\star}\int_M{\mathcal K}_\alpha dv_g + o(1)
\end{equation}
for all $\alpha$, where $o(1) \to 0$ as $\alpha \to +\infty$, and 
$${\mathcal K}_\alpha = \sum_{i=1}^p
\left(\vert\hat u_\alpha^i+u^0_i\vert^{2^\star} - \vert\hat u_\alpha^i\vert^{2^\star} 
- \vert u^0_i\vert^{2^\star}\right)\hskip.1cm .$$
Noting that there exists $C > 0$, independent of $\alpha$, such that
$$\int_M\vert{\mathcal K}_\alpha\vert dv_g
\le C \sum_{i=1}^p\int_M\left(\vert\hat u_\alpha^i\vert^{2^\star-1}\vert u^0_i\vert dv_g 
+ \vert u^0_i\vert^{2^\star-1}\vert\hat u_\alpha^i\vert\right) dv_g$$
for all $\alpha$, and that, by basic integration theory, 
$$\int_M\vert\hat u_\alpha^i\vert^{2^\star-1}\vert u^0_i\vert dv_g = o(1)
\hskip.2cm\hbox{and}\hskip.2cm 
\int_M\vert u^0_i\vert^{2^\star-1}\vert\hat u_\alpha^i\vert dv_g = o(1)$$
for all $\alpha$ and $i$, we get with (\ref{ProofClaim2SobTheorEqt6}) that
$$I_{\alpha,p}({\mathcal U}_\alpha) = I_{\infty,p}({\mathcal U}^0) 
+ I_p(\hat{\mathcal U}_\alpha) + o(1)$$
for all $\alpha$, where $o(1) \to 0$ as $\alpha \to +\infty$. This proves the third assertion in 
(\ref{ProofTheorSobAssert2}). It remains to prove the second assertion, 
namely that $(\hat{\mathcal U}_\alpha)_\alpha$ 
is a Palais-Smale sequence for $I_p$. Let $\Phi = (\varphi_1,\dots,\varphi_p)$ be 
given in $H_{1,p}^2(M)$. With similar arguments to those used above, and thanks to the 
Sobolev inequality, we can write that
$$\sum_{i,j=1}^p\int_MA^\alpha_{ij}u_\alpha^i\varphi_jdv_g 
= \sum_{i,j=1}^p\int_MA_{ij}u^0_i\varphi_j + 
o\left(\Vert\Phi\Vert_{H_{1,p}^2}\right)$$
for all $\alpha$. In particular, since ${\mathcal U}^0$ is a solution of the 
limit system (\ref{LimitSystSec4}), we can write that
\begin{equation}\label{ProofClaim2SobTheorEqt7}
DI_{\alpha,p}({\mathcal U}_\alpha).\Phi = 
DI_p(\hat{\mathcal U}_\alpha).\Phi - \int_M{\mathcal K}_\alpha^\Phi dv_g 
+ o\left(\Vert\Phi\Vert_{H_{1,p}^2}\right)
\end{equation}
for all $\alpha$, where 
$${\mathcal K}_\alpha^\Phi = \sum_{i=1}^p
\left(\vert\hat u_\alpha^i+u^0_i\vert^{2^\star-2}(\hat u_\alpha^i+u^0_i) - \vert\hat u_\alpha^i\vert^{2^\star-2}\hat u_\alpha^i 
- \vert u^0_i\vert^{2^\star-2}u^0_i\right)\varphi_i\hskip.1cm .$$
Noting that there exists $C > 0$, independent of $\alpha$, such that
\begin{eqnarray*}
\int_M\vert{\mathcal K}_\alpha^\Phi\vert dv_g
&\le& C\sum_{i=1}^p\int_M\left(\vert\hat u_\alpha^i\vert^{2^\star-2}\vert u^0_i\vert \vert\varphi_i\vert 
+ \vert u^0_i\vert^{2^\star-2} \vert\hat u_\alpha^i\vert \vert\varphi_i\vert\right)dv_g\\
&\le& C\sum_{i=1}^p\left(\left\Vert\vert\hat u_\alpha^i\vert^{2^\star-2}u^0_i\right\Vert_{2^\star/(2^\star-1)} 
+ \left\Vert\vert\hat u^0_i\vert^{2^\star-2}u_\alpha^i\right\Vert_{2^\star/(2^\star-1)}\right) \Vert\varphi_i
\Vert_{2^\star}
\end{eqnarray*}
for all $\alpha$, and that, by basic integration theory,
$$\left\Vert\vert\hat u_\alpha^i\vert^{2^\star-2}u^0_i\right\Vert_{2^\star/(2^\star-1)} = o(1)
\hskip.2cm\hbox{and}\hskip.2cm 
\left\Vert\vert\hat u^0_i\vert^{2^\star-2}u_\alpha^i\right\Vert_{2^\star/(2^\star-1)} = o(1)$$
for all $\alpha$ and $i$, we get with (\ref{ProofClaim2SobTheorEqt7}) and the Sobolev inequality that
$$DI_{\alpha,p}({\mathcal U}_\alpha).\Phi = 
DI_p(\hat{\mathcal U}_\alpha).\Phi 
+ o\left(\Vert\Phi\Vert_{H_{1,p}^2}\right)\hskip.1cm .$$
Since $({\mathcal U}_\alpha)_\alpha$ is a Palais-Smale sequence for 
$I_{\alpha,p}$, we get that $(\hat{\mathcal U}_\alpha)_\alpha$ 
is a Palais-Smale sequence for $I_p$. 
This proves (\ref{ProofTheorSobAssert2}). For ${\mathcal U} = (u_1,\dots,u_p)$ and 
$\Phi = (\varphi_1,\dots,\varphi_p)$, we clearly have that
$DI_p({\mathcal U}).\Phi = \sum_{i=1}^pDI_1(u_i).\varphi_i$. 
Then an easy remark is that $(\hat{\mathcal U}_\alpha)_\alpha$ 
is a Palais-Smale sequence for $I_p$ if and only if for any $i$, the sequence $(\hat u_\alpha^i)_\alpha$ 
is a Palais-Smale sequence for $I_1$. We may therefore apply 
the result in the $p = 1$ case to the $(\hat u_\alpha^i)_\alpha$'s 
(see, for instance, Druet, Hebey and Robert \cite{DruHebRob2} or Druet and Hebey \cite{DruHeb2} for 
a presentation of the $p=1$ case in the Riemannian setting). 
In particular, we get from this result in the $p=1$ case that for any $i$, there exists $k_i$ integer, 
and generalized $1$-bubbles $(\hat B_{j,\alpha}^i )_\alpha$ 
as in (\ref{DefGenBubbleTheorSob}), $j = 1,\dots, k_i$, such that, up to a subsequence,
\begin{equation}\label{ProofTheorSobFunctCaseEqt}
\hat u_\alpha^i = \sum_{j=1}^{k_i}\hat B_{j,\alpha}^i + R_\alpha
\hskip.1cm \hskip.2cm\hbox{and}\hskip.2cm
I_1(\hat u_\alpha^i) = \sum_{j=1}^{k_i}E_f(u_j^i) + o(1)\hskip.1cm ,
\end{equation}
where 
$o(1) \to 0$ as $\alpha \to +\infty$, $u_j^i$ is the nontrivial solution of $\Delta u = \vert u\vert^{2^\star-2}u$ from 
which the generalized $1$-bubble $(\hat B_{j,\alpha}^i )_\alpha$ is defined, namely
$$\hat B_{j,\alpha}^i(x) = \left(R_{j,\alpha}^i\right)^{\frac{n-2}{2}}\eta_{j,\alpha}^i(x)
u_j^i\left(R_{j,\alpha}^i\exp_{x_{j,\alpha}^i}^{-1}(x)\right)$$
with notations like in (\ref{DefGenBubbleTheorSob}), and
$E_f(u) = \frac{1}{2}\int_{{\mathbb R}^n}\vert\nabla u\vert^2dx - 
\frac{1}{2^\star}\int_{{\mathbb R}^n}\vert u\vert^{2^\star}dx\hskip.1cm .$
We define $\tilde u_{j,\alpha}^i: B_0\left(\delta\mu_{j,\alpha}^i\right) \to {\mathbb R}$, $\delta > 0$ small, by
$$\tilde u_{j,\alpha}^i(x) = (\mu_{j,\alpha}^i )^{\frac{n-2}{2}}\hat u_\alpha^i
\left(\exp_{x_{j,\alpha}^i}\left(\mu_{j,\alpha}^i x\right)\right)\hskip.1cm ,$$
where $\mu_{j,\alpha}^i = (R_{j,\alpha}^i)^{-1}$. 
Then, see, for instance, Druet, Hebey and Robert \cite{DruHebRob2}, 
we also have that, up to a subsequence, $\tilde u_{j,\alpha}^i \to u_j^i$ a.e in ${\mathbb R}^n$ 
as $\alpha \to +\infty$. This holds for all $i = 1,\dots,p$, and all $j = 1,\dots, k_i$. Let 
$\tilde v_\alpha^i: B_0\left(\delta\mu_{j,\alpha}^i\right) \to {\mathbb R}$ be given by
$$\tilde v_{j,\alpha}^i(x) = (\mu_{j,\alpha}^i )^{\frac{n-2}{2}} u_i^0
\left(\exp_{x_{j,\alpha}^i}\left(\mu_{j,\alpha}^i x\right)\right)\hskip.1cm .$$
Noting that $\tilde v_{j,\alpha}^i \to 0$ a.e in ${\mathbb R}^n$ as $\alpha \to +\infty$, 
and since $u_\alpha^i \ge 0$ for all $i$ and all $\alpha$, it follows that $u_j^i \ge 0$ 
for all $i = 1,\dots,p$, and all $j = 1,\dots, k_i$. By the maximum principle for 
scalar equations, we then get that $u_j^i > 0$ 
for all $i = 1,\dots,p$, and all $j = 1,\dots, k_i$. As already mentioned, this is a situation where 
generalized $1$-bubbles are like 
$1$-bubbles. In particular, for any $i = 1,\dots,p$, and any $j = 1,\dots, k_i$, there exists 
a one bubble $(B_{j,\alpha}^i)_\alpha$ such that 
$\hat B_{j,\alpha}^i = B_{j,\alpha}^i + R_{j,\alpha}^i$ 
for all $\alpha$, where $R_{j,\alpha}^i \to 0$ in $H_1^2(M)$ as $\alpha \to +\infty$. 
In other words, we may replace in (\ref{ProofTheorSobFunctCaseEqt}) 
the generalized $1$-bubbles $(\hat B_{j,\alpha}^i )_\alpha$ by $1$-bubbles 
$(B_{j,\alpha}^i )_\alpha$, and letting $k = \sum_{i=1}^pk_i$, we get that (\ref{EqtThmSobTherSec4}) 
follows from (\ref{ProofTheorSobFunctCaseEqt}). By noting that $E_f(u) = K_n^{-n}/n$ when $u$ is a positive solution 
of $\Delta u = \vert u\vert^{2^\star-2}u$, where $K_n$ is as in (\ref{ValueSharpCst}), 
we also get (\ref{EnerSplitTheorSob}) 
with (\ref{ProofTheorSobAssert2}) and (\ref{ProofTheorSobFunctCaseEqt}).
This ends the proof of the theorem and of the remark after the theorem 
concerning the splitting of the energy.
\end{proof}

The weak limit ${\mathcal U}^0$ and $k$ are clearly invariants of the decomposition (\ref{EqtThmSobTherSec4}) 
in Theorem \ref{SobTheorThm}. 
Let $k_i$ be the number of $p$-bubbles $({\mathcal B}_{j,\alpha})_\alpha$'s which are such that the $i^{th}$-component 
of $({\mathcal B}_{j,\alpha})_\alpha$ is 
a $1$-bubble. As we easily get from the proof of Theorem \ref{SobTheorThm}, 
$k = \sum_{i=1}^pk_i$ and the $k_i$'s are also invariants of the 
decomposition (\ref{EqtThmSobTherSec4}). 
Uniqueness conditions for decompositions like (\ref{EqtThmSobTherSec4}) are known in the $p=1$ case. They 
can be found, for instance, in Druet and Hebey \cite{DruHeb2}. Thanks to these conditions, 
an additional result we easily get from the above proof is that if we write two decompositions (\ref{EqtThmSobTherSec4}) 
with respect to two families $({\mathcal B}_{j,\alpha})_\alpha$ and $(\tilde{\mathcal B}_{j,\alpha})_\alpha$ of $p$-bubbles, then, 
up to renumbering, for any $i = 1,\dots,p$, and any $j = 1,\dots,k_i$,
$$\frac{\mu_{j,\alpha}^i}{\tilde\mu_{j,\alpha}^i} \to 1
\hskip.2cm\hbox{and}\hskip.2cm
\frac{d_g(x_{j,\alpha}^i,\tilde x_{j,\alpha}^i)}{\mu_{j,\alpha}^i} \to 0$$
as $\alpha \to +\infty$, where the $x_{j,\alpha}^i$'s and $\mu_{j,\alpha}^i$'s (resp. 
the $\tilde x_{j,\alpha}^i$'s and $\tilde\mu_{j,\alpha}^i$'s) are the centers and weights of the $1$-bubbles 
$(B_{j,\alpha}^i)_\alpha$ (resp. $(\tilde B_{j,\alpha}^i)_\alpha$) from which the $({\mathcal B}_{j,\alpha})_\alpha$'s 
(resp. $(\tilde{\mathcal B}_{j,\alpha})_\alpha$'s) are defined. We also get from the proof of Theorem \ref{SobTheorThm}, 
and structure conditions we know to hold in the $p=1$ case, that for any $i = 1,\dots,p$, and any $j_1, j_2 = 1,\dots,k_i$, 
the structure equations
$$\frac{\mu_{j_1,\alpha}^i}{\mu_{j_2,\alpha}^i} + \frac{\mu_{j_2,\alpha}^i}{\mu_{j_1,\alpha}^i} 
+ \frac{d_g(x_{j_1,\alpha}^i,x_{j_2,\alpha}^i)}{\mu_{j_1,\alpha}^i\mu_{j_2,\alpha}^i} \to +\infty$$
hold as $\alpha \to +\infty$ when $j_1\not= j_2$. An independent result 
we get from the proof of the theorem is that, up to 
replacing $p$-bubbles by generalized $p$-bubbles, Theorem \ref{SobTheorThm} still holds 
if we do not assume that the ${\mathcal U}_\alpha$'s are nonnegative, where, following the definition 
of a $p$-bubble,  
we define a generalized $p$-bubble as a sequence 
of $p$-maps such that one of the components of the sequence is a generalized $1$-bubble, and the other 
components are zero. 

\medskip An easy consequence of Theorem \ref{SobTheorThm} is the $L^{2^\star}$-theory for blow-up 
(Corollary \ref{CorLptheory} below) where 
the blow-up phenomenon is described as a sum of Dirac masses in the $L^{2^\star}$-Lebesgue's space. 
The Dirac masses in Corollary \ref{CorLptheory} are the limits of the $(B_{j,\alpha}^i)^{2^\star}$'s 
as $\alpha \to +\infty$, $i = 1,\dots,p$, $j = 1,\dots,k_i$, 
where the $(B_{j,\alpha}^i)_\alpha$'s are the $1$-bubbles in the $H_1^2$-decomposition 
(\ref{EqtThmSobTherSec4Bis}) following Theorem \ref{SobTheorThm}. 
More direct proofs of Corollary \ref{CorLptheory} can be given.

\begin{cor}\label{CorLptheory} Let $(M,g)$ be a smooth compact Riemannian manifold of dimension 
$n \ge 3$, $p \ge 1$, and $\left(A(\alpha)\right)_\alpha$ a sequence of smooth maps 
$A(\alpha): M \to M_p^s({\mathbb R})$ satisfying (\ref{CondSec4}). 
Let also $({\mathcal U}_\alpha)_\alpha$ be a nonnegative 
Palais-Smale sequence for (\ref{GenericEqtSec4}). For any $i = 1,\dots,p$, up to a subsequence, 
\begin{equation}\label{CorLptheoryEqt}
(u_\alpha^i )^{2^\star} \rightharpoonup (u^0_i)^{2^\star} + \sum_{j=1}^{k^\prime_i}\lambda_j^i\delta_{x_j^i}
\end{equation}
weakly in the sense of measures as $\alpha \to +\infty$, 
where ${\mathcal U}^0 = (u^0_1,\dots,u^0_p)$ is a nonnegative solution of the limit system (\ref{LimitSystSec4}), 
$k^\prime_i$ is an integer, the $x^i_j$'s are points in $M$, and the $\lambda_j^i$'s are positive real numbers, 
$j = 1,\dots,k^\prime_i$.
\end{cor}

\begin{proof} By Theorem \ref{SobTheorThm}, up to a subsequence, we may assume that
(\ref{EqtThmSobTherSec4}) and (\ref{EqtThmSobTherSec4Bis}) hold 
for the ${\mathcal U}_\alpha$'s. For $i = 1,\dots,p$, we let 
$S_i$ be the set consisting of the limits as $\alpha \to +\infty$ of the centers $x_{j,\alpha}^i$ of the $1$-bubbles 
$(B_{j,\alpha}^i)_\alpha$ in (\ref{EqtThmSobTherSec4Bis}). Then we let $k^\prime_i$ be the number 
of points in $S_i$, and let the $x_j^i$'s, $j = 1,\dots,k^\prime_i$, be the points in $S_i$. We have that 
$k^\prime_i \le k_i$, and it might be that $k^\prime_i < k_i$ since distinct sequences may have the same limits. 
It might also be that $k^\prime_i = 0$, and hence that $S_i = \emptyset$ for some $i$. Let $\delta_0 > 0$ be such that $2\delta_0$ 
is less than any distance between two distinct points in $S_i$ (when $k^\prime_i \ge 2$). Since Palais-Smale sequences 
are bounded  in $H_{1,p}^2(M)$, see (\ref{ProofTheorSobAssert1}) in the proof of Theorem \ref{SobTheorThm}, and 
by the Sobolev embedding theorem, the $u_\alpha^i$'s are bounded in $L^{2^\star}(M)$. In particular, up to a subsequence, 
we may assume that for any $i$, and any $j = 1,\dots,k^\prime_i$, there exists $\lambda_j^i$ such that
\begin{equation}\label{CorSec4Eqt1}
\lim_{\alpha\to+\infty}\int_{B_{x_j^i}(\delta_0)}(u_\alpha^i)^{2^\star}dv_g = 
\int_{B_{x_j^i}(\delta_0)}(u^0_i)^{2^\star}dv_g + \lambda_j^i\hskip.1cm .
\end{equation}
We fix $i$ and $j$, and 
let $(B_{j,\alpha}^i)_\alpha$ be a $1$-bubble in (\ref{EqtThmSobTherSec4Bis}) such that its centers 
$x_{j,\alpha}^i$ converge to $x_j^i$ as $\alpha \to +\infty$. With the notations in (\ref{EqtThmSobTherSec4Bis}), 
we can write that
\begin{equation}\label{CorSec4Eqt2}
\begin{split}
\int_{B_{x_j^i}(\delta_0)}(u_\alpha^i)^{2^\star}dv_g
& = \int_{B_{x_j^i}(\delta_0)}\left(u^0_i+\sum_{m=1}^{k_i}B_{m,\alpha}^i\right)^{2^\star}dv_g + o(1)\\
& \ge \int_{B_{x_j^i}(\delta_0)}(u_i^0)^{2^\star}dv_g + 
\int_{B_{x_j^i}(\delta_0)}(B_{j,\alpha}^i)^{2^\star}dv_g + o(1)
\end{split}
\end{equation}
where $o(1) \to 0$ as $\alpha \to +\infty$. Combining (\ref{CorSec4Eqt1}) and (\ref{CorSec4Eqt2}), 
noting that the $L^{2^\star}$-integral of the $1$-bubble $(B_{j,\alpha}^i)_\alpha$ in the right hand 
side of (\ref{CorSec4Eqt2}) goes to $K_n^{-n}$ as $\alpha \to +\infty$, where 
$K_n$ is as in (\ref{ValueSharpCst}), 
it follows that $\lambda_j^i > 0$ for all $i = 1,\dots,p$ and all $j = 1,\dots,k^\prime_i$. 
Let ${\mathcal B}_i(\delta)$ be the union from $j=1$ 
to $k^\prime_i$ of the geodesic balls $B_{x_j^i}(\delta)$, $0 < \delta < \delta_0$. By (\ref{EqtThmSobTherSec4Bis}),  
$u_\alpha^i \to u^0_i$ in 
$L^{2^\star}\left(M\backslash{\mathcal B}_i(\delta)\right)$ for all $\delta > 0$. In particular, 
for any $\delta \in (0,\delta_0)$, any $i = 1,\dots,p$, 
and any $j = 1,\dots,k^\prime_i$,
\begin{equation}\label{CorSec4Eqt3}
\int_{B_{x_j^i}(\delta_0)\backslash B_{x_j^i}(\delta)}(u_\alpha^i)^{2^\star}dv_g 
= \int_{B_{x_j^i}(\delta_0)\backslash B_{x_j^i}(\delta)}(u^0_i)^{2^\star}dv_g + o(1)\hskip.1cm ,
\end{equation}
where $o(1) \to 0$ as $\alpha \to +\infty$.  
Let $f \in C^0(M)$ and $i = 1,\dots,p$. Given $\varepsilon > 0$, 
let also $\delta \in (0,\delta_0)$ be such that $\vert f(x)-f(x_j^i)\vert < \varepsilon$ if $d_g(x_j^i,x) < \delta$, 
$j = 1,\dots,k^\prime_i$. We can write that
\begin{equation}\label{CorSec4Eqt4}
\int_Mf\left((u_\alpha^i)^{2^\star} - (u^0_i)^{2^\star}\right)dv_g 
= \sum_{j=1}^{k^\prime_i} \int_{B_{x_j^i}(\delta)}f\left((u_\alpha^i)^{2^\star} - (u^0_i)^{2^\star}\right)dv_g + o(1)\hskip.1cm ,
\end{equation}
where $o(1) \to 0$ as $\alpha \to +\infty$, and, by (\ref{CorSec4Eqt1}) and (\ref{CorSec4Eqt3}), we 
can also write that
\begin{equation}\label{CorSec4Eqt5}
\lambda_j^if(x_j^i) - C\varepsilon \le 
\int_{B_{x_j^i}(\delta)}f\left((u_\alpha^i)^{2^\star} - (u^0_i)^{2^\star}\right)dv_g
\le \lambda_j^if(x_j^i) + C\varepsilon
\end{equation}
for all $\alpha$ and all $j = 1,\dots,k^\prime_i$, where $C > 0$ does not depend on $\alpha$, $\varepsilon$, $i$, and $j$.
Since $\varepsilon > 0$ is arbitrary, and $f \in C^0(M)$ is arbitrary, 
we get with (\ref{CorSec4Eqt4}) and (\ref{CorSec4Eqt5}) that
(\ref{CorLptheoryEqt}) is true. This proves Corollary \ref{CorLptheory}.
\end{proof}

\section{Pointwise estimates}\label{PointEst}

In what follows we let $(M,g)$ be a smooth compact Riemannian manifold of dimension $n \ge 3$,  $p \ge 1$ integer, 
and $\left(A(\alpha)\right)_\alpha$, $\alpha\in{\mathbb N}$, be a sequence of smooth maps 
$A(\alpha): M \to M_p^s({\mathbb R})$. We let $A(\alpha) = (A_{ij}^\alpha)$, and consider systems like
\begin{equation}\label{GenericEqtSec5}
\Delta_gu_i + \sum_{j=1}^pA_{ij}^\alpha(x)u_j = \vert u_i\vert^{2^\star-2}u_i
\end{equation}
in $M$, for all $i = 1,\dots p$. We assume in this section that the $A(\alpha)$'s satisfy that 
there exists a $C^{0,\theta}$-map $A: M \to M_p^s({\mathbb R})$, 
$A = (A_{ij})$ and $0 < \theta < 1$, such that
\begin{equation}\label{CondSec5}
A^\alpha_{ij} \to A_{ij}\hskip.1cm\hbox{in}\hskip.1cm C^{0,\theta}(M)
\end{equation}
for all $i, j$ as $\alpha \to +\infty$. In particular, (\ref{CondSec4}) is satisfied. 
A sequence $({\mathcal U}_\alpha)_\alpha$ 
is said to be a sequence of nonnegative solutions of (\ref{GenericEqtSec5}) if, for any $\alpha$,  
${\mathcal U}_\alpha$ is a nonnegative solution of (\ref{GenericEqtSec5}). Clearly the sequence is a 
Palais-Smale sequence for (\ref{GenericEqtSec5}) when (and, actually, if and only if) it is also bounded in $H_{1,p}^2(M)$. 
We prove in this section that passing from Palais-Smale sequences to sequences of solutions 
we can add pointwise estimates to the description in Theorem \ref{SobTheorThm}. The main result of 
this section is as follows.

\begin{thm}[Pointwise Estimates]\label{PointEstThm} Let $(M,g)$ be a smooth compact Riemannian 
manifold of dimension 
$n \ge 3$, $p \ge 1$, and $\left(A(\alpha)\right)_\alpha$ a sequence of smooth maps 
$A(\alpha): M \to M_p^s({\mathbb R})$ satisfying (\ref{CondSec5}).
Let $({\mathcal U}_\alpha)_\alpha$ 
be a bounded sequence in $H_{1,p}^2(M)$ of nonnegative solutions of (\ref{GenericEqtSec5}). 
In addition to the decomposition (\ref{EqtThmSobTherSec4}) in Theorem \ref{SobTheorThm}, there exists $C > 0$ such that, 
up to a subsequence, 
\begin{equation}\label{PointEstThmEqt}
\left(\min_{i,j}d_g(x_{j,\alpha}^i,x)\right)^{\frac{n-2}{2}} 
\sqrt{\sum_{i=1}^p\left(u_\alpha^i(x)-u^0_i(x)\right)^2} \le C
\end{equation}
for all $\alpha$ and all $x \in M$, where ${\mathcal U}^0$ is as in 
(\ref{EqtThmSobTherSec4}), ${\mathcal U}_\alpha = (u_\alpha^1,\dots,u_\alpha^p)$ for all $\alpha$, 
the $u_i^0$'s are the components of ${\mathcal U}^0$, and the 
$x_{j,\alpha}^i$'s are the centers of the $1$-bubbles $(B_{j,\alpha}^i)_\alpha$ in 
(\ref{EqtThmSobTherSec4Bis}) from which the 
$p$-bubbles $({\mathcal B}_{j,\alpha})_\alpha$'s in (\ref{EqtThmSobTherSec4}) are defined.
\end{thm}

\begin{proof} Let $\Phi_\alpha$ be the function such that $\Phi_\alpha(x)$ is the minimum over $i, j$ of 
the $d_g(x_{j,\alpha}^i,x)$'s, where $x \in M$, and let $\Psi_\alpha$ be the function
\begin{equation}\label{PsiDefProofThPointEst}
\Psi_\alpha(x) = 
\sqrt{\sum_{i=1}^pu_\alpha^i(x)^2}
\Phi_\alpha(x)^{\frac{n-2}{2}}
\hskip.1cm .
\end{equation}
In order to prove (\ref{PointEstThmEqt}), it suffices to prove that $(\Psi_\alpha)_\alpha$ is bounded 
in $L^\infty(M)$. We proceed by contradiction. We let the $y_\alpha$'s be points in $M$ such that 
the $\Psi_\alpha$'s are maximum at $y_\alpha$ and $\Psi_\alpha(y_\alpha) \to +\infty$ as $\alpha \to +\infty$. 
Up to a subsequence, we may assume that $u_\alpha^{i_0}(y_\alpha) \ge u_\alpha^i(y_\alpha)$ 
for some $i_0=1,\dots,p$, and all $i$. We set $\mu_\alpha = u_\alpha^{i_0}(y_\alpha)^{-2/(n-2)}$. Then 
$\mu_\alpha \to 0$ as $\alpha \to +\infty$, and by (\ref{PsiDefProofThPointEst}) we also have that
\begin{equation}\label{PsiDefProofThPointEstEqt1}
\frac{d_g(x_{j,\alpha}^i,y_\alpha)}{\mu_\alpha} \to +\infty
\end{equation}
for all $i ,j$, as $\alpha \to +\infty$. Let $\delta > 0$ be less than the injectivity radius of 
$(M,g)$. For $i = 1,\dots,p$, we define the function $v_\alpha^i$ in $B_0(\delta\mu_\alpha^{-1})$ by
\begin{equation}\label{PsiDefProofThPointEstEqt2}
v_\alpha^i(x) = \mu_\alpha^{\frac{n-2}{2}}u_\alpha^i\left(\exp_{y_\alpha}(\mu_\alpha x)\right)
\hskip.1cm ,
\end{equation}
where $B_0(\delta\mu_\alpha^{-1})$ is the Euclidean ball of radius $\delta\mu_\alpha^{-1}$ 
centered at $0$, and $\exp_{y_\alpha}$ is the exponential map at $y_\alpha$. Given $R > 0$ 
and $x \in B_0(R)$, the Euclidean ball of radius $R$ centered at $0$, we can write with 
(\ref{PsiDefProofThPointEst}) and (\ref{PsiDefProofThPointEstEqt2}) that
\begin{equation}\label{PsiDefProofThPointEstEqt3}
v_\alpha^i(x) \le \frac{\mu_\alpha^{\frac{n-2}{2}}\Psi_\alpha\left(\exp_{y_\alpha}(\mu_\alpha x)\right)}
{\Phi_\alpha\left(\exp_{y_\alpha}(\mu_\alpha x)\right)^{\frac{n-2}{2}}}
\end{equation}
for all $i$, when $\alpha$ is sufficiently large. For any $i, j$, and $x \in B_0(R)$,
\begin{eqnarray*} d_g\left(x_{j,\alpha}^i,\exp_{y_\alpha}(\mu_\alpha x)\right)
& \ge & d_g\left(x_{j,\alpha}^i,y_\alpha\right) - R\mu_\alpha\\
& \ge & \left(1 - \frac{R\mu_\alpha}{\Phi_\alpha(y_\alpha)}\right)\Phi_\alpha(y_\alpha)
\end{eqnarray*}
when $\alpha$ is sufficiently large so that, by 
(\ref{PsiDefProofThPointEstEqt1}), the right hand side of the last equation is positive. Coming back to 
(\ref{PsiDefProofThPointEstEqt3}), thanks to the definition of the $y_\alpha$'s, we then get that 
for any $i$, and any $x \in B_0(R)$,
\begin{equation}\label{PsiDefProofThPointEstEqt4}
\begin{split}
v_\alpha^i(x)
& \le \frac{\mu_\alpha^{\frac{n-2}{2}}\Psi_\alpha(y_\alpha)}
{\Phi_\alpha\left(\exp_{y_\alpha}(\mu_\alpha x)\right)^{\frac{n-2}{2}}}\\
& \le \sqrt{p} \left(1 - \frac{R\mu_\alpha}{\Phi_\alpha(y_\alpha)}\right)^{-\frac{n-2}{2}}
\end{split}
\end{equation}
when $\alpha$ is sufficiently large. In particular, by (\ref{PsiDefProofThPointEstEqt1}) and 
(\ref{PsiDefProofThPointEstEqt4}), up to passing 
to a subsequence, the $v_\alpha^i$'s are uniformly bounded in any compact subset 
of ${\mathbb R}^n$ for all $i$. Let ${\mathcal V}_\alpha = (v_\alpha^1,\dots,v_\alpha^p)$. The 
${\mathcal V}_\alpha$'s are solutions of the system
\begin{equation}\label{PsiDefProofThPointEstEqt5}
\Delta_{g_\alpha}v_\alpha^i + \sum_{j=1}^p\mu_\alpha^2\tilde A_{ij}^\alpha v_\alpha^j = (v_\alpha^i)^{2^\star-1}
\hskip.1cm ,
\end{equation}
where
\begin{eqnarray*}
&&\tilde A_{ij}^\alpha(x) = A_{ij}^\alpha\left(\exp_{y_\alpha}(\mu_\alpha x)\right)
\hskip.1cm ,\hskip.1cm\hbox{and}\\
&&g_\alpha(x) 
= \left(\exp_{y_\alpha}^\star g\right)(\mu_\alpha x)\hskip.1cm .
\end{eqnarray*}
Let $\xi$ be the Euclidean metric. Clearly, for any compact subset 
$K$ of ${\mathbb R}^n$, $g_\alpha \to \xi$ in $C^2(K)$ as $\alpha \to +\infty$. Then, by standard elliptic theory, 
we get that the $v_\alpha^i$'s are uniformly bounded 
in $C^{2,\theta}_{loc}({\mathbb R}^n)$ for all $i$, where $0 < \theta < 1$. In particular, up to a subsequence, 
we can assume that $v_\alpha^i \to v_i$ in $C^2_{loc}({\mathbb R}^n)$ as $\alpha \to +\infty$ 
for all $i$, where the $v_i$'s 
are nonnegative functions in $C^2({\mathbb R}^n)$. The $v_i$'s are bounded in ${\mathbb R}^n$ by 
(\ref{PsiDefProofThPointEstEqt4}), and such that $v_{i_0}(0) = 1$ by construction. Without loss of generality, we 
may also assume that the $v_i$'s are in ${\mathcal D}_1^2({\mathbb R}^n)$ and in $L^{2^\star}({\mathbb R}^n)$ 
for all $i$, where 
${\mathcal D}_1^2({\mathbb R}^n)$ is the Beppo-Levi space defined as the completion of $C^\infty_0({\mathbb R}^n)$, 
the space of smooth function with compact support in ${\mathbb R}^n$, with respect to the norm 
$\Vert u\Vert = \Vert\nabla u\Vert_2$. We let ${\mathcal V} = (v_1,\dots,v_p)$. According to the above, 
${\mathcal V} \not\equiv 0$ 
is nonnegative and nonzero. For any $i$, and any $R > 0$, 
$$\int_{B_{y_\alpha}(R\mu_\alpha)}(u_\alpha^i)^{2^\star}dv_g 
= \int_{B_0(R)}(v_\alpha^i)^{2^\star}dv_{g_\alpha}\hskip.1cm .$$
It follows that for any $i$, and any $R > 0$,
\begin{equation}\label{PsiDefProofThPointEstEqt6}
\int_{B_{y_\alpha}(R\mu_\alpha)}(u_\alpha^i)^{2^\star}dv_g 
= \int_{{\mathbb R}^n}(v_i)^{2^\star}dx + \varepsilon_R(\alpha)\hskip.1cm ,
\end{equation}
where $\varepsilon_R(\alpha)$ is such that $\lim_R\lim_\alpha\varepsilon_R(\alpha) = 0$, 
and the limits are as $\alpha\to+\infty$ and $R\to+\infty$. 
Thanks to the decomposition (\ref{EqtThmSobTherSec4}) in Theorem \ref{SobTheorThm}, see also 
(\ref{EqtThmSobTherSec4Bis}), we can write 
that for any $i$, and any $R >0$,
$$\int_{B_{y_\alpha}(R\mu_\alpha)}(u_\alpha^i)^{2^\star}dv_g = 
\int_{B_{y_\alpha}(R\mu_\alpha)}\left(u^0_i + \sum_{j=1}^{k_i}B_{j,\alpha}^i + R_\alpha^i\right)^{2^\star}dv_g
\hskip.1cm ,$$
where $R_\alpha^i \to 0$ in 
$H_1^2(M)$ as $\alpha \to +\infty$. In particular, we get that for any $i$, and any $R >0$,
\begin{equation}\label{PsiDefProofThPointEstEqt7}
\int_{B_{y_\alpha}(R\mu_\alpha)}(u_\alpha^i)^{2^\star}dv_g \le 
C\sum_{j=1}^{k_i}\int_{B_{y_\alpha}(R\mu_\alpha)}(B_{j,\alpha}^i)^{2^\star}dv_g + o(1)
\hskip.1cm ,
\end{equation}
where $o(1) \to 0$ as $\alpha \to +\infty$, and $C > 0$ is independent of $\alpha$, $i$, and $R$. 
Now we claim that, thanks to (\ref{PsiDefProofThPointEstEqt1}),
\begin{equation}\label{PsiDefProofThPointEstEqt8}
\lim_{\alpha\to+\infty}\int_{B_{y_\alpha}(R\mu_\alpha)}(B_{j,\alpha}^i)^{2^\star}dv_g = 0
\end{equation}
for all $R > 0$ and all $i, j$. In order to prove (\ref{PsiDefProofThPointEstEqt8}), 
we distinguish two cases. In what follows 
we fix $R > 0$, $i$, and $j$, and let 
the $\mu_{j,\alpha}^i$'s be the weights of the $1$-bubbles $(B_{j,\alpha}^i)_\alpha$. 
In the first case we assume that for any $R^\prime > 0$, up to a subsequence,
$$B_{y_\alpha}(R\mu_\alpha)\bigcap B_{x_{j,\alpha}^i}(R^\prime \mu_{j,\alpha}^i) = \emptyset$$
for all $\alpha$. Then,
$$\int_{B_{y_\alpha}(R\mu_\alpha)}(B_{j,\alpha}^i)^{2^\star}dv_g \le 
\int_{M\backslash B_{x_{j,\alpha}^i}(R^\prime \mu_{j,\alpha}^i)}(B_{j,\alpha}^i)^{2^\star}dv_g\hskip.1cm .$$
Noting that
$$\int_{M\backslash B_{x_{j,\alpha}^i}(R^\prime \mu_{j,\alpha}^i)}(B_{j,\alpha}^i)^{2^\star}dv_g 
= \varepsilon_{R^\prime}(\alpha)\hskip.1cm ,$$
where
$\lim_{R^\prime}\lim_\alpha\varepsilon_{R^\prime}(\alpha) = 0$, 
and the limits are as $\alpha\to+\infty$ and $R^\prime\to+\infty$, 
we get that (\ref{PsiDefProofThPointEstEqt8}) is true in this case. 
In the second case we assume that there exists $R^\prime > 0$ such that, 
up to a subsequence,
$$B_{y_\alpha}(R\mu_\alpha)\bigcap B_{x_{j,\alpha}^i}(R^\prime \mu_{j,\alpha}^i) \not= \emptyset$$
for all $\alpha$. Then
$$d_g(x_{j,\alpha}^i,y_\alpha) \le R\mu_\alpha + R^\prime\mu_{j,\alpha}^i$$
and it follows from (\ref{PsiDefProofThPointEstEqt1}) that
$\mu_\alpha = o(\mu_{j,\alpha}^i)$ and $d_g(x_{j,\alpha}^i,y_\alpha) = O(\mu_{j,\alpha}^i)$.
Writing that
$$B_{y_\alpha}(R\mu_\alpha) \subset 
\exp_{x_{j,\alpha}^i}\left(\mu_{j,\alpha}^iB_{z_\alpha}\left(C\Lambda_\alpha\right)\right)
\hskip.1cm ,$$
where
$z_\alpha = \frac{1}{\mu_{j,\alpha}^i}\exp_{x_{j,\alpha}^i}^{-1}(y_\alpha)$ 
converges in ${\mathbb R}^n$ (up to a subsequence), $C > 1$ is independent of $\alpha$ and $R$, 
and $\Lambda_\alpha = R\mu_\alpha/\mu_{j,\alpha}^i$, we then get that 
$$\int_{B_{y_\alpha}(R\mu_\alpha)}(B_{j,\alpha}^i)^{2^\star}dv_g 
\le \int_{B_{z_\alpha}\left(C\Lambda_\alpha\right)}u^{2^\star}dv_{\tilde g_\alpha}\hskip.1cm ,$$
where $u$ is given by (\ref{PosSolCritEuclEqt}), and $\tilde g_\alpha$ is the metric given by 
$\tilde g_\alpha(x) = \left(\exp_{x_{j,\alpha}^i}^\star g\right)(\mu_{j,\alpha}^ix)$.
Since $\mu_\alpha = o(\mu_{j,\alpha}^i)$, we have that
$$\int_{B_{z_\alpha}\left(C\Lambda_\alpha\right)}u^{2^\star}dv_{\tilde g_\alpha} = o(1)$$
and this proves (\ref{PsiDefProofThPointEstEqt8}) in this case. In particular, (\ref{PsiDefProofThPointEstEqt8}) is true, 
and coming back to (\ref{PsiDefProofThPointEstEqt6}) and (\ref{PsiDefProofThPointEstEqt7}), we get that, for any $i$, and 
any $R > 0$,
\begin{equation}\label{PsiDefProofThPointEstEqt9}
\int_{{\mathbb R}^n}(v_i)^{2^\star}dx = \varepsilon_R(\alpha)\hskip.1cm ,
\end{equation}
where $\varepsilon_R(\alpha)$ is such that $\lim_R\lim_\alpha\varepsilon_R(\alpha) = 0$, 
and the limits are as $\alpha\to+\infty$ and $R\to+\infty$. Letting $\alpha \to +\infty$, and then $R \to +\infty$, this implies that 
$\int_{{\mathbb R}^n}(v_i)^{2^\star}dx = 0$ for all $i$. Since ${\mathcal V} \not\equiv 0$, we get the desired contradiction. 
The theorem is proved.
\end{proof}

Let $({\mathcal U}_\alpha)_\alpha$ 
be a bounded sequence in $H_{1,p}^2(M)$ of nonnegative solutions of (\ref{GenericEqtSec5}). Up to passing to a subsequence, 
the decomposition (\ref{EqtThmSobTherSec4}) in Theorem \ref{SobTheorThm}, and the estimate 
(\ref{PointEstThmEqt}) of Theorem \ref{PointEstThm} are satisfied by the ${\mathcal U}_\alpha$'s. 
We let ${\mathcal S}_{geom}$ be the set consisting of the limits of the $x_{j,\alpha}^i$'s, 
where the $x_{j,\alpha}^i$'s, are the centers of the $1$-bubbles $(B_{j,\alpha}^i)_\alpha$ 
in (\ref{EqtThmSobTherSec4Bis}) from which the 
$p$-bubbles $({\mathcal B}_{j,\alpha})_\alpha$'s in (\ref{EqtThmSobTherSec4}) are defined.  
Let ${\mathcal U}_\alpha = (u_\alpha^1,\dots,u_\alpha^p)$. Easy consequences of 
Theorem \ref{PointEstThm} are that ${\mathcal S}_{geom} \not= \emptyset$ if and only if the sequences 
$(u_\alpha^i)_\alpha$ 
are not all bounded in $L^\infty(M)$, but also that the  
sequences $(u_\alpha^i)_\alpha$ 
are all bounded in $L^\infty_{loc}(M\backslash{\mathcal S}_{geom})$, $i = 1,\dots,p$. From 
(\ref{GenericEqtSec5}), and by standard elliptic theory, it easily follows that
\begin{equation}\label{COConvOutside}
u_\alpha^i \to u^0_i\hskip.2cm\hbox{in}\hskip.1cm C^2_{loc}(M\backslash{\mathcal S}_{geom})
\end{equation}
as $\alpha \to +\infty$, for all $i$, where the $u^0_i$'s are, as in (\ref{PointEstThmEqt}), the components of 
${\mathcal U}^0$ in (\ref{EqtThmSobTherSec4}). 
The points in 
${\mathcal S}_{geom}$ are referred to as the {\it geometrical blow-up points} of the sequence $({\mathcal U}_\alpha)_\alpha$. 
Needless to say, since different sequences may have the same limits, ${\mathcal S}_{geom}$ may consist of any 
number $m \le k$ of points, where $k$ is the number of $p$-bubbles involved in (\ref{EqtThmSobTherSec4}). 
We say that the sequence $({\mathcal U}_\alpha)_\alpha$ {\it blows up} when 
${\mathcal S}_{geom} \not= \emptyset$.

\medskip An important complement to Theorem \ref{PointEstThm} is given by the following lemma.

\begin{lem}\label{LemRefPtEst} Let $(M,g)$ be a smooth compact Riemannian 
manifold of dimension 
$n \ge 3$, $p \ge 1$, and $\left(A(\alpha)\right)_\alpha$ a sequence of smooth maps 
$A(\alpha): M \to M_p^s({\mathbb R})$ satisfying (\ref{CondSec5}).
Let $({\mathcal U}_\alpha)_\alpha$ 
be a bounded sequence in $H_{1,p}^2(M)$ of nonnegative solutions of (\ref{GenericEqtSec5}). 
In addition to the decomposition (\ref{EqtThmSobTherSec4}) in Theorem \ref{SobTheorThm}, 
and to the estimate (\ref{PointEstThmEqt}) in Theorem \ref{PointEstThm}, there also holds that, up to a subsequence, 
$$\lim_{R\to+\infty} \lim_{\alpha\to+\infty}
\left[~\sup_{x \in M\backslash\Omega_\alpha(R)}\left(\min_{i,j}d_g(x_{j,\alpha}^i,x)\right)^{\frac{n-2}{2}} 
\sqrt{\sum_{i=1}^p\left(u_\alpha^i(x)-u^0_i(x)\right)^2}~\right] = 0\hskip.1cm ,$$
where ${\mathcal U}^0$ is as in 
(\ref{EqtThmSobTherSec4}), ${\mathcal U}_\alpha = (u_\alpha^1,\dots,u_\alpha^p)$ for all $\alpha$, 
the $u_i^0$'s are the components of ${\mathcal U}^0$, the 
$x_{j,\alpha}^i$'s and $\mu^i_{j,\alpha}$'s are the centers 
and weights of the $1$-bubbles $(B_{j,\alpha}^i)_\alpha$ in (\ref{EqtThmSobTherSec4Bis}) 
from which the 
$p$-bubbles $({\mathcal B}_{j,\alpha})_\alpha$'s in (\ref{EqtThmSobTherSec4}) are defined,  and where 
$\Omega_\alpha(R)$ for $R > 0$ is given by $\Omega_\alpha(R) = \bigcup_{i,j}B_{x^i_{j,\alpha}}(R\mu^i_{j,\alpha})$.
\end{lem}

\begin{proof} As in the proof of Theorem \ref{PointEstThm}, we let $\Phi_\alpha$ be the function 
such that $\Phi_\alpha(x)$ is the minimum over $i, j$ of 
the $d_g(x_{j,\alpha}^i,x)$'s, where $x \in M$. We proceed by contradiction and assume that there exists 
a sequence $(y_\alpha)_\alpha$ of points in $M$, and that there exists $\delta_0 > 0$ such that for any 
$i, j$,
\begin{equation}\label{ProofLemSec5Eqt1}
\frac{d_g(x_{j,\alpha}^i,y_\alpha)}{\mu^i_{j,\alpha}} \to +\infty
\end{equation}
as $\alpha \to +\infty$, and such that
\begin{equation}\label{ProofLemSec5Eqt2}
\Phi_\alpha(y_\alpha)^{\frac{n-2}{2}}\sum_{i=1}^p\left\vert u_\alpha^i(y_\alpha)-u^0_i(y_\alpha)\right\vert \ge \delta_0
\end{equation}
for all $\alpha$. By (\ref{COConvOutside}), $\Phi_\alpha(y_\alpha) \to 0$ as $\alpha \to +\infty$ since, if not, 
$u_\alpha^i(y_\alpha)-u^0_i(y_\alpha) \to 0$ as $\alpha \to +\infty$ for all $i$.
Up to a subsequence, we may assume that $u_\alpha^{i_0}(y_\alpha) \ge u_\alpha^i(y_\alpha)$ 
for some $i_0=1,\dots,p$, and all $i$. We set $\mu_\alpha = u_\alpha^{i_0}(y_\alpha)^{-2/(n-2)}$. Since 
$\Phi_\alpha(y_\alpha)u^0_i(y_\alpha) \to 0$ as $\alpha \to +\infty$ for all $i$, up to passing to another subsequence, 
we get with 
(\ref{ProofLemSec5Eqt2}) that
\begin{equation}\label{ProofLemSec5Eqt3}
\frac{\Phi_\alpha(y_\alpha)}{\mu_\alpha} \ge \delta_1
\end{equation}
for all $\alpha$, and some $\delta_1 > 0$. In particular, $\mu_\alpha \to 0$ as 
$\alpha \to +\infty$. Let $\delta > 0$ be less than the injectivity radius of 
$(M,g)$. For $i = 1,\dots,p$, we define the function $v_\alpha^i$ in $B_0(\delta\mu_\alpha^{-1})$ by
\begin{equation}\label{ProofLemSec5Eqt4}
v_\alpha^i(x) = \mu_\alpha^{\frac{n-2}{2}}u_\alpha^i\left(\exp_{y_\alpha}(\mu_\alpha x)\right)
\hskip.1cm ,
\end{equation}
where $B_0(\delta\mu_\alpha^{-1})$ is the Euclidean ball of radius $\delta\mu_\alpha^{-1}$ 
centered at $0$, and $\exp_{y_\alpha}$ is the exponential map at $y_\alpha$. By 
(\ref{ProofLemSec5Eqt3}) we can write that for any sequence $(x_\alpha)_\alpha$ of points in 
$B_0(\delta_1/2)$, and for any $i, j$,
\begin{eqnarray*} d_g\left(x_{j,\alpha}^i,\exp_{y_\alpha}(\mu_\alpha x_\alpha)\right) 
& \ge & d_g\left(x_{j,\alpha}^i,y_\alpha\right) - d_g\left(y_\alpha,\exp_{y_\alpha}(\mu_\alpha x_\alpha)\right)\\
& \ge & \delta_1\mu_\alpha - \vert x_\alpha\vert \mu_\alpha\\
& \ge & C \mu_\alpha
\end{eqnarray*}
for all $\alpha$, where $C > 0$ is independent of $\alpha$. With the estimate 
(\ref{PointEstThmEqt}) of Theorem \ref{PointEstThm},
we then get that there exists $C > 0$ such that 
\begin{equation}\label{ProofLemSec5Eqt5}
v_\alpha^i(x) \le C
\end{equation}
for all $x \in B_0(\delta_1/2)$, all $i$, and all $\alpha$. In particular, we may now follow the arguments developed 
in the proof of Theorem \ref{PointEstThm}. We let 
${\mathcal V}_\alpha = (v_\alpha^1,\dots,v_\alpha^p)$. Then the  
${\mathcal V}_\alpha$'s are solutions of the system
\begin{equation}\label{ProofLemSec5Eqt6}
\Delta_{g_\alpha}v_\alpha^i + \sum_{j=1}^p\mu_\alpha^2\tilde A_{ij}^\alpha v_\alpha^j = (v_\alpha^i)^{2^\star-1}
\hskip.1cm ,
\end{equation}
where $\tilde A_{ij}^\alpha(x) = A_{ij}^\alpha\left(\exp_{y_\alpha}(\mu_\alpha x)\right)$, and $g_\alpha(x) 
= \left(\exp_{y_\alpha}^\star g\right)(\mu_\alpha x)$. 
Let $\xi$ be the Euclidean metric. Clearly, for any compact subset 
$K$ of ${\mathbb R}^n$, $g_\alpha \to \xi$ in $C^2(K)$ as $\alpha \to +\infty$. 
By standard elliptic theory, (\ref{ProofLemSec5Eqt5}), and 
(\ref{ProofLemSec5Eqt6}), we then get that
the $v_\alpha^i$'s are uniformly bounded 
in $C^{2,\theta}\left(B_0(\delta_1/4)\right)$ for all $i$, where $0 < \theta < 1$. In particular, 
up to another subsequence, we can assume that the $v_\alpha^i$'s converge in 
$C^2\left(B_0(\delta_1/8)\right)$ for all $i$. If $v_i$ is the limit of the $v_\alpha^i$'s, it 
follows from the definition of $\mu_\alpha$ that $v_{i_0}(0) = 1$. Let $\delta_2 = \delta_1/8$. 
For any $i$, 
\begin{equation}\label{ProofLemSec5Eqt7}
\begin{split}
\int_{B_{y_\alpha}(\delta_2\mu_\alpha)}(u_\alpha^i)^{2^\star}dv_g 
&= \int_{B_0(\delta_2)}(v_\alpha^i)^{2^\star}dv_{g_\alpha}\\
&= \int_{B_0(\delta_2)}(v_i)^{2^\star}dx + o(1)\hskip.1cm ,
\end{split}
\end{equation}
where $o(1) \to 0$ as $\alpha\to+\infty$. Thanks to the decomposition 
(\ref{EqtThmSobTherSec4}) in Theorem \ref{SobTheorThm}, see also 
(\ref{EqtThmSobTherSec4Bis}), we can write 
that for any $i$, 
\begin{equation}\label{ProofLemSec5Eqt8}
\int_{B_{y_\alpha}(\delta_2\mu_\alpha)}(u_\alpha^i)^{2^\star}dv_g \le 
C\sum_{j=1}^{k_i}\int_{B_{y_\alpha}(\delta_2\mu_\alpha)}(B_{j,\alpha}^i)^{2^\star}dv_g + o(1)
\hskip.1cm ,
\end{equation}
where $o(1) \to 0$ as $\alpha\to+\infty$, and $C > 0$ is independent of $\alpha$ and $i$. As in the proof 
of Theorem \ref{PointEstThm}, see (\ref{PsiDefProofThPointEstEqt8}), we can also write that
\begin{equation}\label{ProofLemSec5Eqt9}
\lim_{\alpha\to+\infty}\int_{B_{y_\alpha}(\delta_2\mu_\alpha)}(B_{j,\alpha}^i)^{2^\star}dv_g = 0
\end{equation}
for all $i, j$. We prove (\ref{ProofLemSec5Eqt9}) as we prove (\ref{PsiDefProofThPointEstEqt8}) by 
considering the two cases where $B_{y_\alpha}(\delta_2\mu_\alpha)\cap B_{x^i_{j,\alpha}}(R\mu_{j,\alpha}^i) = \emptyset$ 
for all $R > 0$, and $B_{y_\alpha}(\delta_2\mu_\alpha)\cap B_{x^i_{j,\alpha}}(R\mu_{j,\alpha}^i) \not= \emptyset$ 
for some $R > 0$. In the second case we recover (\ref{PsiDefProofThPointEstEqt1}) thanks to 
(\ref{ProofLemSec5Eqt1}) by noting that (\ref{ProofLemSec5Eqt3}) and the nonempty intersection give that 
$\delta_1\mu_\alpha \le \delta_2\mu_\alpha + R\mu^i_{j,\alpha}$ so that $\mu_\alpha \le C \mu_{j,\alpha}^i$. 
Then, combining (\ref{ProofLemSec5Eqt7})--(\ref{ProofLemSec5Eqt9}), we get that 
$$\int_{B_0(\delta_2)}(v_i)^{2^\star}dx = 0$$
for all $i$, and taking $i = i_0$, we get a contradiction with the equation $v_{i_0}(0) = 1$. 
Lemma \ref{LemRefPtEst} is proved. 
\end{proof}

\section{$L^2$-concentration}\label{L2Conc}

In what follows we let $(M,g)$ be a smooth compact Riemannian manifold of dimension $n \ge 3$,  $p \ge 1$ integer, 
and $\left(A(\alpha)\right)_\alpha$, $\alpha\in{\mathbb N}$, be a sequence of smooth maps 
$A(\alpha): M \to M_p^s({\mathbb R})$. We let also $A(\alpha) = (A_{ij}^\alpha)$, and consider systems like
\begin{equation}\label{GenericEqtSec6}
\Delta_gu_i + \sum_{j=1}^pA_{ij}^\alpha(x)u_j = \vert u_i\vert^{2^\star-2}u_i
\end{equation}
in $M$, for all $i = 1,\dots p$. We assume in this section that the $A(\alpha)$'s satisfy that 
there exists a $C^{0,\theta}$-map $A: M \to M_p^s({\mathbb R})$, 
$A = (A_{ij})$ and $0 < \theta < 1$, such that
\begin{equation}\label{CondSec6}
\begin{split}
&\Delta_g^p+A\hskip.1cm\hbox{is coercive, and}\\
&A^\alpha_{ij} \to A_{ij}\hskip.1cm\hbox{in}\hskip.1cm C^{0,\theta}(M)
\end{split}
\end{equation}
as $\alpha \to +\infty$, where the second equation in (\ref{CondSec6}) 
should hold for all $i, j$. The goal in what follows is to prove 
$L^2$-concentration (Theorem \ref{L2ConcThm} below) for sequences $({\mathcal U}_\alpha)_\alpha$ of 
nonnegative solutions of (\ref{GenericEqtSec6}).

\medskip A remark concerning the coercivity assumption in (\ref{CondSec6}) is that when $-A$ is cooperative, 
the existence of strongly positive solutions to systems like (\ref{GenericEqtAlgForm}) 
implies the coercivity of the operator $\Delta_g^p+A$ in 
(\ref{GenericEqtAlgForm}). For such systems, like when $p = 1$, coercivity follows from the existence 
of positive solutions. 
Let $\lambda_A(g)$ be, as in Section \ref{ExSol}, the infimum of the functional 
$I_A$ in (\ref{FonctionalDef}) over the ${\mathcal U} = (u_1,\dots,u_p)$ which are such that 
$\sum_i\int_Mu_i^2dv_g = 1$. By compactness of 
the embedding of $H_1^2$ in $L^2$, we easily get that there exists a minimizer ${\mathcal U}_A 
\in H_{1,p}^2(M)$ for $\lambda_A(g)$. If $-A$ is cooperative, then ${\mathcal U}_A$ 
can be chosen weakly positive. We let ${\mathcal U}_A = (u^A_1,\dots,u^A_p)$, and let also 
${\mathcal U} = (u_1,\dots,u_p)$ be a solution of (\ref{GenericEqtAlgForm}).
Since $(\Delta_g^p+A){\mathcal U}_A = \lambda_A(g){\mathcal U}_A$, we can write that 
\begin{eqnarray*}
\sum_{i=1}^p\int_Mu^A_iu_i^{2^\star-1}dv_g
& = & \sum_i\int_Mu^A_i\left(\Delta_gu_i + \sum_{j=1}^pA_{ij}u_j\right)dv_g\\
& = & \sum_i\int_Mu_i\left(\Delta_gu^A_i + \sum_{j=1}^pA_{ij}u^A_j\right)dv_g\\
 & =  & \lambda_A(g) \sum_{i=1}^p\int_Mu_iu^A_idv_g\hskip.1cm .
\end{eqnarray*}
In particular, 
we get that if $-A$ is cooperative and ${\mathcal U}$ is 
a strongly positive solution of (\ref{GenericEqtAlgForm}), 
then $\lambda_A(g) > 0$ and (see Section \ref{ExSol})
$\Delta_g^p+A$ is coercive. This proves the above claim that when $-A$ is cooperative, 
the existence of a strongly positive solution to a system like (\ref{GenericEqtAlgForm}) implies that the operator 
$\Delta_g^p+A$ in the left hand side of (\ref{GenericEqtAlgForm}) 
is coercive. In general, when no 
assumption is made on $A$, and contrary to the scalar case $p = 1$, the existence 
of a strongly positive solution to a system like (\ref{GenericEqtAlgForm}) 
when $p \ge 2$ does not imply the 
coercivity of $\Delta_g^p + A$. For instance, with the 
examples (\ref{HighEnergEx1})--(\ref{HighEnergEx1StructCond}) of Section \ref{HighEn}, 
one easily constructs $2$-systems like (\ref{GenericEqtAlgForm}) with strongly positive solutions 
and such that $A_{12} \gg A_{11} +  A_{22}$, and $A_{11}, A_{22} > 0$. In particular, 
the operator $\Delta_g^2+A$ is not coercive.

\medskip Before we discuss $L^2$-concentration, we need to prove a 
De Giorgi-Nash-Moser iterative scheme 
for our systems. Let $A: M \to M_p^s({\mathbb R})$, $A = (A_{ij})$, be a continuous map. Let also ${\mathcal U} = (u_1,\dots,u_p)$ 
be a nonnegative $p$-map in $H_{1,p}^2(M)$. We say that ${\mathcal U}$ satisfies that
\begin{equation}\label{GenericSystDeGNM}
\Delta_gu_i + \sum_{j=1}^pA_{ij}(x)u_j \le u_i^{2^\star-1}
\end{equation}
in the sense of distributions, for all $i = 1,\dots p$, if for any $\Phi=(\varphi_1,\dots,\varphi_p)$ in $H_{1,p}^2(M)$, 
$\Phi$ nonnegative, and any $i$, we can write that
$$\int_M\left(\nabla u_i\nabla\varphi_i\right)dv_g + \sum_{j=1}^p\int_MA_{ij}u_i\varphi_jdv_g
\le \int_M(u_i)^{2^\star-1}\varphi_idv_g\hskip.1cm ,$$
where $\left(\nabla u_i\nabla\varphi_i\right)$ is the pointwise scalar product of 
$\nabla u_i$ and $\nabla\varphi_i$. In what follows, for ${\mathcal U} = (u_1,\dots,u_p)$, and $q > 0$, we 
let $\vert{\mathcal U}\vert^q = \sum_{i=1}^p\vert u_i\vert^q$. In particular, when $q = 1$, we let 
$\vert{\mathcal U}\vert = \sum_{i=1}^p\vert u_i\vert$. For $\Omega \subset M$ 
we let also $\Vert{\mathcal U}\Vert_{L^\theta(\Omega)} = \sum_i\Vert u_i\Vert_{L^\theta(\Omega)}$.
The De Giorgi-Nash-Moser iterative scheme 
for our systems is as follows.

\begin{lem}\label{DeGNMoser} Let $(M,g)$ be a smooth compact Riemannian 
manifold of dimension $n \ge 3$,  $p \ge 1$ integer, 
and $A: M \to M_p^s({\mathbb R})$, $A = (A_{ij})$, be a continuous map. 
If ${\mathcal U} \in H_{1,p}^2(M)$ is  
a nonnegative $p$-map satisfying (\ref{GenericSystDeGNM}) in the sense of distributions, then 
$\vert{\mathcal U}\vert \in L^\infty(M)$. Moreover, for any 
$x \in M$, any $\Lambda > 0$, any $\delta > 0$, any $s \ge 1$, and any $q > 2^\star$, if ${\mathcal U}$ 
satisfies also that $\Vert{\mathcal U}\Vert_{L^q(B_x(2\delta))} \le \Lambda$, then
\begin{equation}\label{DeGNMEqt}
\max_{y\in B_x(\delta)}\left\vert{\mathcal U}(y)\right\vert \le C \Vert{\mathcal U}\Vert_{L^s(B_x(2\delta))}
\hskip.1cm ,
\end{equation}
where  $C > 0$ does not depend on ${\mathcal U}$. 
\end{lem}

\begin{proof} Let ${\mathcal U} = (u_1,\dots,u_p)$, ${\mathcal U} \not\equiv 0$, 
be a nonnegative $p$-map satisfying (\ref{GenericSystDeGNM}) in the sense of distributions. Applying a 
Trudinger type argument like in the proof of Theorem \ref{ExistThm}, we easily get that 
${\mathcal U} \in L^k(M)$ for some 
$k > 2^\star$. In particular, the first claim in Lemma \ref{DeGNMoser} follows from the second 
claim. Summing the different equations 
in (\ref{GenericSystDeGNM}), we can write that 
\begin{equation}\label{ProofDGNMEqt1}
\begin{split}
\Delta_g\vert{\mathcal U}\vert 
&\le C\vert{\mathcal U}\vert + \sum_{i=1}^pu_i^{2^\star-1}\\
&\le C \left(1 + \vert{\mathcal U}\vert^{2^\star-2}\right)\vert{\mathcal U}\vert
\end{split}
\end{equation}
where $C > 0$ does not depend on ${\mathcal U}$, $\vert{\mathcal U}\vert = \sum_{i=1}^pu_i$ since ${\mathcal U}$ 
is nonnegative, 
and the inequality is to be understood in the sense of distributions. 
Namely, (\ref{ProofDGNMEqt1}) holds in the sense that
\begin{equation}\label{ProofDGNMEqt2}
\int_M(\nabla\vert{\mathcal U}\vert\nabla\varphi)dv_g \le 
C \int_M\left(1 + \vert{\mathcal U}\vert^{2^\star-2}\right)\vert{\mathcal U}\vert\varphi dv_g
 \end{equation}
 for all nonnegative $\varphi \in H_1^2(M)$. If $\Vert{\mathcal U}\Vert_{L^q(B_x(2\delta))} \le \Lambda$ for some 
 $q > 2^\star$, the function $f = \vert{\mathcal U}\vert^{2^\star-2}$ is bounded independently of 
 ${\mathcal U}$ in $L^s\left(B_x(2\delta)\right)$ for some $s > n/2$. In particular, 
by (\ref{ProofDGNMEqt2}), we can apply the standard De Giorgi-Nash-Moser 
 iterative scheme for functions. We get that ${\mathcal U} \in L^\infty(M)$ if 
 ${\mathcal U} \in L^k(M)$ for some 
$k > 2^\star$, and that for any 
$x \in M$, any $\Lambda > 0$, any $\delta > 0$, any $s > 0$, and any $q > 2^\star$, if ${\mathcal U}$ 
satisfies that $\Vert{\mathcal U}\Vert_{L^q(B_x(2\delta))} \le \Lambda$, then
(\ref{DeGNMEqt}) is true. This proves the lemma.
\end{proof}

Concerning the dependency of the constant $C$ in (\ref{DeGNMEqt}) with respect to $A$, 
it easily follows from the above proof that $C$ can be chosen to depend only on the $C^0$-norm of 
$K = \sum_{i,j}\vert A_{ij}\vert$. Another easy remark is that if $A$ is cooperative, then
$L_gu_i \le u_i^{2^\star-1}$ for all $i$, where $L_g = \Delta_g+A_{ii}$, 
and the De Giorgi-Nash-Moser iterative scheme for functions can be applied directly to the 
$u_i$'s. 

\medskip Now that we have a De Giorgi-Nash-Moser iterative scheme for our systems, we return 
to (\ref{GenericEqtSec6}). We assume that (\ref{CondSec6}) holds, and let $({\mathcal U}_\alpha)_\alpha$ 
be a sequence of solutions of (\ref{GenericEqtSec6}). Namely, for any $\alpha$,  
${\mathcal U}_\alpha$ is a solution of (\ref{GenericEqtSec6}). By the second equation 
in (\ref{CondSec6}), there exists $C > 0$ such that for $\alpha \ge \alpha_0$ sufficiently large, 
$I_{A(\alpha)}({\mathcal U}) \ge CI_A({\mathcal U})$ for all ${\mathcal U} \in H_{1,p}^2(M)$, 
where $I_{A(\alpha)}$ and $I_A$ are as in (\ref{FonctionalDef}). By (\ref{GenericEqtSec6}) 
and the coercivity assumption 
in (\ref{CondSec6}) we then get that for $\alpha \ge \alpha_0$ sufficiently large,
\begin{eqnarray*} \int_M\vert{\mathcal U}_\alpha\vert^{2^\star}dv_g
& \ge & C_1 I_A({\mathcal U}_\alpha)\\
& \ge & C_2 \left(\int_M\vert{\mathcal U}_\alpha\vert^{2^\star}dv_g\right)^{2/2^\star}\hskip.1cm ,
\end{eqnarray*}
where $C_1, C_2 > 0$ are independent of $\alpha$. 
In particular, when we assume (\ref{CondSec6}), there exists $K > 0$ such that for any 
sequence $({\mathcal U}_\alpha)_\alpha$ of solutions of (\ref{GenericEqtSec6}), 
and any $\alpha \ge \alpha_0$, 
$\int_M\vert{\mathcal U}_\alpha\vert^{2^\star}dv_g \ge K$ 
if ${\mathcal U}_\alpha \not\equiv 0$. Now $L^2$-concentration states 
as follows. For ${\mathcal U} = (u_1,\dots,u_p)$, 
we write $\Vert{\mathcal U}\Vert_2 = \sum_i\Vert u_i\Vert_2$.

\begin{thm}[$L^2$-concentration]\label{L2ConcThm} Let $(M,g)$ be a smooth compact Riemannian 
manifold of dimension $n \ge 4$,  $p \ge 1$ integer, 
and $\left(A(\alpha)\right)_\alpha$, $A(\alpha): M \to M_p^s({\mathbb R})$, be a sequence of smooth maps 
satisfying (\ref{CondSec6}). Let $({\mathcal U}_\alpha)_\alpha$, ${\mathcal U}_\alpha \not\equiv 0$, 
be a bounded sequence in $H_{1,p}^2(M)$ 
of nonnegative solutions of (\ref{GenericEqtSec6}) such that 
$\Vert{\mathcal U}_\alpha\Vert_2 \to 0$ as $\alpha \to +\infty$. Then, up to a subsequence, 
${\mathcal S}_{geom} \not= \emptyset$ and
\begin{equation}\label{L2ConcEqtThm}
\lim_{\alpha\to+\infty} \frac{\int_{{\mathcal B}_\delta}\vert{\mathcal U}_\alpha\vert^2dv_g}
{\int_M\vert{\mathcal U}_\alpha\vert^2dv_g} = 1
\end{equation}
for all $\delta > 0$, 
where ${\mathcal B}_\delta = \bigcup_{i=1}^mB_{x_i}(\delta)$, and the $x_i $'s are the 
geometrical blow-up points of the 
${\mathcal U}_\alpha$'s, namely ${\mathcal S}_{geom} = \left\{x_1,\dots,x_m\right\}$.
\end{thm}

When $n =3$, bubbles do not concentrate in the $L^2$-norm, 
and $L^2$-concentration turns out to be false when $n = 3$. Dimension $4$ is the 
smallest dimension for which we can get $L^2$-concentration. 

\begin{proof} Up to a subsequence, Theorem \ref{SobTheorThm} and Theorem \ref{PointEstThm} 
apply to the 
${\mathcal U}_\alpha$'s. Since $\Vert{\mathcal U}_\alpha\Vert_2 \to 0$ as $\alpha \to +\infty$, we 
have that ${\mathcal U}^0 \equiv 0$ in (\ref{EqtThmSobTherSec4}). If in addition 
${\mathcal S}_{geom} = \emptyset$, we would get by  (\ref{COConvOutside}) that 
$\vert{\mathcal U}_\alpha\vert \to 0$ uniformly in $M$ as $\alpha \to +\infty$, and this is in contradiction 
with the lower bound $\int_M\vert{\mathcal U}_\alpha\vert^{2^\star}dv_g \ge K$ we discussed 
above. In other words, up to a subsequence, there is 
a nonempty finite set ${\mathcal S}_{geom}$ of geometrical blow-up points for the 
${\mathcal U}_\alpha$'s. By the second equation 
in (\ref{CondSec6}), there exists $C > 0$ such that for $\alpha \ge \alpha_0$ sufficiently large, 
$I_{A(\alpha)}({\mathcal U}) \ge CI_A({\mathcal U})$ for all ${\mathcal U} \in H_{1,p}^2(M)$, 
where $I_{A(\alpha)}$ and $I_A$ are as in (\ref{FonctionalDef}).
Let $f = (f_1,\dots,f_p)$ 
be a smooth map such that $f_i > 0$ in $M$ for all $i$. For instance, $f_i \equiv 1$ 
for all $i$. 
Minimizing $I_{A(\alpha)}$ over the constraint $\sum_i\int_Mf_iu_idv_g = 1$, we easily get 
the existence of ${\mathcal W}_\alpha \in H_{1,p}^2(M)$ such that
\begin{equation}\label{CoercConsL2Conc}
\Delta_g^p{\mathcal W}_\alpha + A(\alpha){\mathcal W}_\alpha = f
\end{equation}
for all $\alpha \ge \alpha_0$. By the coercivity assumption and (\ref{CoercConsL2Conc}), 
$\Vert{\mathcal W}_\alpha\Vert_2 \le C$, and by 
standard regularity results, we get that ${\mathcal W}_\alpha$ is smooth 
and that there exists $K > 0$ such that $\vert{\mathcal W}_\alpha\vert \le K$ in $M$ for all $\alpha \ge \alpha_0$. 
Then, by (\ref{CoercConsL2Conc}), since $f_i > 0$ for all $i$, and $\vert{\mathcal W}_\alpha\vert \le K$, 
we can write that for any $\alpha \ge \alpha_0$,
\begin{equation}\label{L1LPControlL2Conc}
\begin{split}
\sum_{i=1}^p\int_Mu_\alpha^idv_g 
&\le C \sum_{i=1}^p\int_Mf_iu_\alpha^idv_g\\
&= C \sum_{i=1}^p\int_M\left(\Delta_gw_\alpha^i + \sum_jA^\alpha_{ij}w_\alpha^j\right)u_\alpha^idv_g\\
&= C \sum_{i=1}^p\int_M\left(\Delta_gu_\alpha^i + \sum_jA^\alpha_{ij}u_\alpha^j\right)w_\alpha^idv_g\\
&\le C \sum_{i=1}^p\int_M\vert w_\alpha^i\vert (u_\alpha^i)^{2^\star-1}dv_g\\
&\le C \sum_{i=1}^p\int_M(u_\alpha^i)^{2^\star-1}dv_g\hskip.1cm ,
\end{split}
\end{equation}
where $C > 0$ is independent of $\alpha$, and 
${\mathcal W}_\alpha = (w_\alpha^1,\dots,w_\alpha^p)$. In particular, by 
(\ref{L1LPControlL2Conc}), there exists $C > 0$ such that
\begin{equation}\label{L1LPContrBalEqt}
\int_M\vert{\mathcal U}_\alpha\vert dv_g \le 
C \int_M\vert{\mathcal U}_\alpha\vert^{2^\star-1}dv_g
\end{equation}
for all $\alpha \ge \alpha_0$. 
We refer to (\ref{L1LPContrBalEqt}) as the $L^1/L^{2^\star-1}$-controlled 
balance property of the system (\ref{GenericEqtSec6}). 
Up to a subsequence, (\ref{L1LPContrBalEqt}) holds for all $\alpha$. 
Now, for $\delta > 0$, we let ${\mathcal B}_\delta$ 
be as in the statement of the theorem, and let $R_\delta(\alpha)$ be the ratio
\begin{equation}\label{RatioL2Conc}
R_\delta(\alpha) = \frac{\int_{M\backslash{\mathcal B}_\delta}\vert{\mathcal U}_\alpha\vert^2dv_g}
{\int_M\vert{\mathcal U}_\alpha\vert^2dv_g}\hskip.1cm .
\end{equation}
Thanks, for instance, to Theorem \ref{PointEstThm}, we can apply 
the De Giorgi-Nash-Moser iterative scheme in 
Lemma \ref{DeGNMoser} to the ${\mathcal U}_\alpha$'s in $M\backslash{\mathcal B}_\delta$ 
with $s = 2$. By the $L^1/L^{2^\star-1}$-controlled balance (\ref{L1LPContrBalEqt}), 
and the De Giorgi-Nash-Moser iterative scheme, we then get that 
\begin{eqnarray*}
\int_{M\backslash{\mathcal B}_\delta}\vert{\mathcal U}_\alpha\vert^2dv_g
& \le & \left(\max_{M\backslash{\mathcal B}_\delta}\vert{\mathcal U}_\alpha\vert\right) 
\int_{M\backslash{\mathcal B}_\delta}\vert{\mathcal U}_\alpha\vert dv_g\\
& \le & C \sqrt{\int_M\vert{\mathcal U}_\alpha\vert^2dv_g} 
\int_M\vert{\mathcal U}_\alpha\vert^{2^\star-1}dv_g\hskip.1cm ,
\end{eqnarray*}
where $C > 0$ is independent of $\alpha$. In particular, 
\begin{equation}\label{RatioL2ConcProof}
R_\delta(\alpha) \le C \frac{\int_M\vert{\mathcal U}_\alpha\vert^{2^\star-1}dv_g}
{\sqrt{\int_M\vert{\mathcal U}_\alpha\vert^2dv_g}}
\end{equation}
for all $\alpha$, where $C > 0$ is independent of $\alpha$. If we assume now that $n \ge 6$, then 
$2^\star-1 \le 2$, and we can write with H\"older's inequality that
$$\int_M(u_\alpha^i)^{2^\star-1}dv_g \le V_g^{\frac{3-2^\star}{2}} 
\left(\int_M(u_\alpha^i)^2dv_g\right)^{\frac{2^\star-1}{2}}$$
for all $i$, where $V_g$ is the volume of $M$ with respect to $g$. In particular, there exists 
$C > 0$ such that
$$\int_M\vert{\mathcal U}_\alpha\vert^{2^\star-1}dv_g \le C 
\left(\int_M\vert{\mathcal U}_\alpha\vert^2dv_g\right)^{\frac{2^\star-1}{2}}\hskip.1cm ,$$
and since $2^\star > 2$ and $\Vert{\mathcal U}_\alpha\Vert_2 \to 0$, we get with 
(\ref{RatioL2ConcProof}) that $R_\delta(\alpha) \to 0$ as 
$\alpha \to +\infty$. This proves (\ref{L2ConcEqtThm}) when $n \ge 6$. If $n = 5$, then 
$2 \le 2^\star-1 \le 2^\star$, and we can write with H\"older's inequality that
$$\left(\int_M(u_\alpha^i)^{2^\star-1}dv_g\right)^{\frac{1}{2^\star-1}}
\le \left(\int_M(u_\alpha^i)^2dv_g\right)^{\frac{\theta}{2}}
\left(\int_M(u_\alpha^i)^{2^\star}dv_g\right)^{\frac{1-\theta}{2^\star}}\hskip.1cm ,$$
where $\theta = \frac{3}{2(2^\star-1)}$. Since the sequence $({\mathcal U}_\alpha)_\alpha$ is bounded 
in $H_{1,p}^2(M)$, there exists $C > 0$ such that $\int_M\vert{\mathcal U}_\alpha\vert^{2^\star}dv_g \le C$ 
for all $\alpha$. Then we get that
$$\int_M\vert{\mathcal U}_\alpha\vert^{2^\star-1}dv_g \le C 
\left(\int_M\vert{\mathcal U}_\alpha\vert^2dv_g\right)^{\frac{3}{4}}\hskip.1cm ,$$
where $C > 0$ does not depend on $\alpha$, and since $\frac{3}{4} > \frac{1}{2}$ and 
$\Vert{\mathcal U}_\alpha\Vert_2 \to 0$, we get with 
(\ref{RatioL2ConcProof}) that $R_\delta(\alpha) \to 0$ as 
$\alpha \to +\infty$. This proves (\ref{L2ConcEqtThm}) when $n = 5$. Now it remains to prove 
(\ref{L2ConcEqtThm}) when $n = 4$. The argument when $n = 4$ is slightly more delicate and requires 
the $H_1^2$-theory in Section \ref{SobTheory}. 
We start writing that
\begin{equation}\label{L2ConcProofDim4Eqt1}
\begin{split}
\frac{\int_M\vert{\mathcal U}_\alpha\vert^{2^\star-1}dv_g}{\sqrt{\int_M\vert{\mathcal U}_\alpha\vert^2dv_g}} 
&= \sum_{i=1}^p
\frac{\int_M(u_\alpha^i)^{2^\star-1}dv_g}{\sqrt{\int_M\vert{\mathcal U}_\alpha\vert^2dv_g}}\\
&\le 
\sum_{i=1}^p\frac{\int_M(u_\alpha^i)^{2^\star-1}dv_g}{\sqrt{\int_M(u_\alpha^i)^2dv_g}}\hskip.1cm .
\end{split}
\end{equation}
For $i = 1,\dots,p$, we let 
the $x_{j,\alpha}^i$'s and $\mu_{j,\alpha}^i$'s be the centers and weights 
of the $1$-bubbles $(B_{j,\alpha}^i)_\alpha$ in (\ref{EqtThmSobTherSec4Bis}) 
from which the $p$-bubbles $({\mathcal B}_{j,\alpha})_\alpha$'s in 
(\ref{EqtThmSobTherSec4}) are defined, where $j = 1,\dots,k_i$, 
and $\sum_i k_i = k$. For $R > 0$, and $i = 1,\dots,p$, we let also 
$\Omega_{i,\alpha}(R)$ be given by
\begin{equation}\label{L2ConcProofDim4Eqt2}
\Omega_{i,\alpha}(R) = \bigcup_{j=1}^{k_i}B_{x_{j,\alpha}^i}(R\mu_{j,\alpha}^i)\hskip.1cm ,
\end{equation}
and $\Omega_{i,\alpha}(R) = \emptyset$ if $u_\alpha^i \to 0$ in $H_1^2(M)$ as $\alpha \to +\infty$. 
We fix $i = 1,\dots,p$. 
Since $2^\star = 4$, we can write, thanks to H\"older's inequalities, that
\begin{eqnarray*}
\int_M(u_\alpha^i)^{2^\star-1}dv_g 
&\le& \int_{\Omega_{i,\alpha}(R)}(u_\alpha^i)^{2^\star-1}dv_g\\
&&\hskip.4cm + \sqrt{\int_{M\backslash\Omega_{i,\alpha}(R)}(u_\alpha^i)^{2^\star}dv_g}
\sqrt{\int_M(u_\alpha^i)^2dv_g}
\end{eqnarray*}
Then we get that
\begin{equation}\label{L2ConcProofDim4Eqt3}
\frac{\int_M(u_\alpha^i)^{2^\star-1}dv_g}{\sqrt{\int_M(u_\alpha^i)^2dv_g}} 
\le \sqrt{\int_{M\backslash\Omega_{i,\alpha}(R)}(u_\alpha^i)^{2^\star}dv_g} 
+ \frac{\int_{\Omega_{i,\alpha}(R)}(u_\alpha^i)^{2^\star-1}dv_g}{\sqrt{\int_M(u_\alpha^i)^2dv_g}}\hskip.1cm ,
\end{equation}
where $\Omega_{i,\alpha}(R)$ is as in (\ref{L2ConcProofDim4Eqt2}). 
For $\varphi \in C^\infty_0({\mathbb R}^n)$, where $C^\infty_0({\mathbb R}^n)$ is the set of smooth functions with 
compact support in ${\mathbb R}^n$, we let $\varphi_{j,\alpha}^i$ be
the function defined by the equation 
\begin{equation}\label{L2ConcProofDim4Eqt4}
\varphi_{j,\alpha}^i(x) = (\mu_{j,\alpha}^i)^{-\frac{n-2}{2}}
\varphi\left((\mu_{j,\alpha}^i)^{-1}\exp_{x_{j,\alpha}^i}^{-1}(x)\right)
\hskip.1cm .
\end{equation}
Straightforward computations give that 
for any $j_1 \not= j_2$, any $i$, and any $\alpha$,\par
\medskip (i) $\displaystyle\int_M(B_{j_2,\alpha}^i)^{2^\star-1}\varphi_{j_1,\alpha}^idv_g = o(1)$,\par
\medskip\noindent where $o(1) \to 0$ as $\alpha \to +\infty$. Similarly, for any $R > 0$, 
any $i$, any $j$, and any $\alpha$,\par
\medskip (ii) $\displaystyle\int_{M\backslash \Omega_{j,\alpha}^i(R)}
(B_{j,\alpha}^i)^{2^\star}dv_g = \varepsilon_R(\alpha)$,\par
\medskip (iii) $\displaystyle\int_{\Omega_{j,\alpha}^i(R)}
(B_{j,\alpha}^i)^{2^\star-1}\varphi_{j,\alpha}^idv_g 
= \int_{B_0(R)}u^{2^\star-1}\varphi dx + o(1)$,\par
\medskip (iv) $\displaystyle\int_{\Omega_{j,\alpha}^i(R)}
(B_{j,\alpha}^i)^2(\varphi_{j,\alpha}^i)^{2^\star-2}dv_g 
= \int_{B_0(R)}u^2\varphi^{2^\star-2}dx + o(1)$\par
\medskip\noindent where 
$u$ is as in (\ref{PosSolCritEuclEqt}), $\Omega_{j,\alpha}^i(R) = B_{x_{j,\alpha}^i}(R\mu_{j,\alpha}^i)$, 
$o(1) \to 0$ as $\alpha \to +\infty$, and the $\varepsilon_R(\alpha)$'s are such that
\begin{equation}\label{EqtRestL2ConcDim4}
\lim_{R \to +\infty}\limsup_{\alpha\to +\infty}\varepsilon_R(\alpha) = 0
\hskip.1cm .
\end{equation}
By (ii) we can write that
\begin{equation}\label{L2ConcProofDim4Eqt5}
\int_{M\backslash\Omega_{i,\alpha}(R)}(u_\alpha^i)^{2^\star}dv_g = \varepsilon_R(\alpha)\hskip.1cm ,
\end{equation}
where $\Omega_{i,\alpha}(R)$ is as in (\ref{L2ConcProofDim4Eqt2}), and 
the $\varepsilon_R(\alpha)$'s are such that (\ref{EqtRestL2ConcDim4}) holds.
From now on, we let $\varphi$ in (\ref{L2ConcProofDim4Eqt4}) be such that $\varphi = 1$ in 
$B_0(R)$. Then,
$$\int_{\Omega_{i,\alpha}(R)}(u_\alpha^i)^{2^\star-1}dv_g 
\le \sum_{j=1}^{k_i}(\mu_{j,\alpha}^i)^{\frac{n-2}{2}}
\int_{\Omega_{j,\alpha}^i(R)}(u_\alpha^i)^{2^\star-1}\varphi_{j,\alpha}^idv_g$$
while, thanks to (i) and (iii), 
\begin{eqnarray*} \int_{\Omega_{j,\alpha}^i(R)}(u_\alpha^i)^{2^\star-1}\varphi_{j,\alpha}^idv_g
& \le & C\int_{\Omega_{j,\alpha}^i(R)}(B_{j,\alpha}^i)^{2^\star-1}\varphi_{j,\alpha}^idv_g + o(1)\\
& \le & C \int_{B_0(R)}u^{2^\star-1}dx + o(1)\hskip.1cm ,
\end{eqnarray*}
where $o(1) \to 0$ as $\alpha \to +\infty$, and $C > 0$ does not depend on $\alpha$ 
and $R$. In particular, we can write that
\begin{equation}\label{L2ConcProofDim4Eqt6}
\int_{\Omega_{i,\alpha}(R)}(u_\alpha^i)^{2^\star-1}dv_g 
\le \left(C \int_{B_0(R)}u^{2^\star-1}dx + o(1)\right)
\sum_{j=1}^{k_i}(\mu_{j,\alpha}^i)^{\frac{n-2}{2}}\hskip.1cm ,
\end{equation}
where $o(1) \to 0$ as $\alpha \to +\infty$, $u$ is as in (\ref{PosSolCritEuclEqt}), 
$k_i$ is as in (\ref{EqtThmSobTherSec4Bis}), 
and $C > 0$ does not depend on $\alpha$ and $R$.
Independently, for any $j$,
\begin{eqnarray*} \int_M(u_\alpha^i)^2dv_g
& \ge & \int_{\Omega_{j,\alpha}^i(R)}(u_\alpha^i)^2dv_g\\
& \ge & (\mu_{j,\alpha}^i)^{n-2}
\int_{\Omega_{j,\alpha}^i(R)}(u_\alpha^i)^2(\varphi_{j,\alpha}^i)^{2^\star-2}dv_g
\end{eqnarray*}
Here, $2^\star-2 = 2$. As is easily checked, we can write 
with (\ref{EqtThmSobTherSec4Bis}) that
\begin{eqnarray*}
&&\int_{\Omega_{j,\alpha}^i(R)}(u_\alpha^i)^2(\varphi_{j,\alpha}^i)^{2^\star-2}dv_g\\
&& =  \int_{\Omega_{j,\alpha}^i(R)}\left(\sum_{m = 1}^{k_i}B_{m,\alpha}^i\right)^2
(\varphi_{j,\alpha}^i)^{2^\star-2}dv_g + o(1)\\
&& \ge \int_{\Omega_{j,\alpha}^i(R)}(B_{j,\alpha}^i)^2(\varphi_{j,\alpha}^i)^{2^\star-2}dv_g + o(1)
\end{eqnarray*}
and thanks to (iv) we get that
$$\int_{\Omega_{j,\alpha}^i(R)}(u_\alpha^i)^2(\varphi_{j,\alpha}^i)^{2^\star-2}dv_g
\ge \int_{B_0(R)}u^2dx + o(1)\hskip.1cm .$$
Hence, for any $j$,
$$\int_M(u_\alpha^i)^2dv_g \ge (\mu_{j,\alpha}^i)^{n-2} 
\left(\int_{B_0(R)}u^2dx + o(1)\right)$$
and we can write that
\begin{equation}\label{L2ConcProofDim4Eqt7}
\int_M(u_\alpha^i)^2dv_g \ge \left(\max_{j=1,\dots,k_i}\mu_{j,\alpha}^i\right)^{n-2} 
\left(\int_{B_0(R)}u^2dx + o(1)\right)\hskip.1cm ,
\end{equation}
where $o(1) \to 0$ as $\alpha \to +\infty$, $k_i$ is as in (\ref{EqtThmSobTherSec4Bis}), 
and $u$ is as in (\ref{PosSolCritEuclEqt}). Then, since $i$ is arbitrary, we get by combining 
(\ref{RatioL2ConcProof}), (\ref{L2ConcProofDim4Eqt1}), (\ref{L2ConcProofDim4Eqt3}), 
(\ref{L2ConcProofDim4Eqt5}), (\ref{L2ConcProofDim4Eqt6}), and 
(\ref{L2ConcProofDim4Eqt7}) that, for any $R > 0$, 
\begin{equation}\label{L2ConcProofDim4Eqt8}
\limsup_{\alpha\to+\infty}R_\delta(\alpha) \le \varepsilon_R + 
C \frac{\int_{B_0(R)}u^{2^\star-1}dx}{\sqrt{\int_{B_0(R)}u^2dx}}\hskip.1cm ,
\end{equation}
where $\varepsilon_R \to 0$ as $R \to +\infty$, and $C > 0$ does not depend 
on $R$. It is easily seen that
\begin{eqnarray*} \lim_{R\to+\infty}\int_{B_0(R)}u^{2^\star-1}dx
& = & \int_{{\mathbb R}^n}u^{2^\star-1}dx\\
& < & +\infty
\end{eqnarray*}
On the other hand, when $n = 4$,
$$\lim_{R\to+\infty}\int_{B_0(R)}u^2dx = +\infty\hskip.1cm .$$
Coming back to (\ref{L2ConcProofDim4Eqt8}), it follows that 
$R_\delta(\alpha) \to 0$ as $\alpha \to +\infty$, and (\ref{L2ConcEqtThm}) when $n = 4$ 
is true.  This ends the proof of Theorem \ref{L2ConcThm}.
\end{proof}

A possible estimate we could have add to 
$L^2$-concentration is the $L^2/\nabla L^2$-balance stating that for any $\delta > 0$, and any $x \in M$, 
there exists $C > 0$ such that 
\begin{equation}\label{LemSec6BosEqt2}
\sum_{i=1}^p\int_{B_x(\delta)}\vert\nabla u_\alpha^i\vert^2dv_g \le C 
\sum_{i=1}^p\int_{B_x(2\delta)}\left(1 + (u_\alpha^i)^{2^\star-2}\right)(u_\alpha^i)^2dv_g
\end{equation}
for all $\alpha$, where the $u_\alpha^i$'s are the components of the ${\mathcal U}_\alpha$'s. 
We prove (\ref{LemSec6BosEqt2}) by mutliplying the $i$th equation 
of (\ref{GenericEqtSec6}) by $\eta^2u_\alpha^i$, and integrating over $M$, 
where $\eta$ is such that 
$\eta = 0$ in $B_x(\delta)$ and $\eta = 1$ in $M\backslash B_x(2\delta)$. 
Using that by (\ref{CondSec6}), there exists $C > 0$ such that 
$-A(\alpha) \le CId_p$ for all $\alpha$ in the sense of bilinear forms, where $Id_p$ is the 
$p\times p$-identity matrix in $M^s_p({\mathbb R})$, we get (\ref{LemSec6BosEqt2}). 
In particular, it follows from  (\ref{LemSec6BosEqt2}), the pointwise estimate in Section 
\ref{PointEst}, and $L^2$-concentration, that
$$\int_{M\backslash {\mathcal B}_{\delta}}\vert\nabla{\mathcal U}_\alpha\vert^2dv_g 
= o\left(\int_M\vert{\mathcal U}_\alpha\vert^2dv_g\right)$$
for all $\alpha$, where  
$\vert\nabla{\mathcal U}_\alpha\vert^2 = \sum_{i=1}^p\vert\nabla u_\alpha^i\vert^2$, 
$\vert{\mathcal U}_\alpha\vert^2 = \sum_{i=1}^p(u_\alpha^i)^2$, 
${\mathcal B}_\delta = \bigcup_{i=1}^mB_{x_i}(\delta)$, and the $x_i $'s are the 
geometrical blow-up points of the 
${\mathcal U}_\alpha$'s.

\section{Sharp pointwise asymptotics}\label{SharpPtAsympt}

In what follows we let $(M,g)$ be a smooth compact Riemannian manifold of dimension $n \ge 3$,  $p \ge 1$ integer, 
and $\left(A(\alpha)\right)_\alpha$, $\alpha\in{\mathbb N}$, be a sequence of smooth maps 
$A(\alpha): M \to M_p^s({\mathbb R})$. We let also $A(\alpha) = (A_{ij}^\alpha)$, and consider systems like
\begin{equation}\label{GenericEqtSec7}
\Delta_gu_i + \sum_{j=1}^pA_{ij}^\alpha(x)u_j = \vert u_i\vert^{2^\star-2}u_i
\end{equation}
in $M$, for all $i = 1,\dots p$. We assume as in Section 
\ref{L2Conc} that the $A(\alpha)$'s satisfy that 
there exists a $C^{0,\theta}$-map $A: M \to M_p^s({\mathbb R})$, 
$A = (A_{ij})$ and $0 < \theta < 1$, such that
\begin{equation}\label{CondSec7}
\begin{split}
&\Delta_g^p+A\hskip.1cm\hbox{is coercive, and}\\
&A^\alpha_{ij} \to A_{ij}\hskip.1cm\hbox{in}\hskip.1cm C^{0,\theta}(M)
\end{split}
\end{equation}
as $\alpha \to +\infty$, where the second equation in (\ref{CondSec7}) 
should hold for all $i, j$. The limit system we get by 
combining (\ref{GenericEqtSec7}) and (\ref{CondSec7}) reads as
\begin{equation}\label{LimitSystSec7}
\Delta_gu_i + \sum_{j=1}^pA_{ij}(x)u_j = \vert u_i\vert^{2^\star-2}u_i
\end{equation}
in $M$, for all $i = 1,\dots p$. The goal in this section is to prove 
sharp pointwise asymptotics for sequences of 
nonnegative solutions of (\ref{GenericEqtSec7}) when standing close to 
one specific bubble of the $H_1^2$-decomposition of Theorem \ref{SobTheorThm}. 

\medskip Let $({\mathcal U}_\alpha)_\alpha$, ${\mathcal U}_\alpha \not\equiv 0$, 
be a bounded sequence in $H_{1,p}^2(M)$ 
of nonnegative solutions of (\ref{GenericEqtSec7}). Then the sequence $({\mathcal U}_\alpha)_\alpha$ is a 
Palais-Smale sequence for (\ref{GenericEqtSec7}) and we can apply 
Theorem \ref{SobTheorThm}.  In particular, up to a subsequence, (\ref{EqtThmSobTherSec4}) 
and (\ref{EqtThmSobTherSec4Bis}) hold. In what follows, we 
assume that the sequence blows up 
and let the $x_{j,\alpha}^i$'s and $\mu_{j,\alpha}^i$'s be the centers and weights 
of the $1$-bubbles $(B_{j,\alpha}^i)_\alpha$ in (\ref{EqtThmSobTherSec4Bis}) 
from which the $p$-bubbles $({\mathcal B}_{j,\alpha})_\alpha$'s in 
(\ref{EqtThmSobTherSec4}) are defined, where $i = 1,\dots,p$, $j = 1,\dots,k_i$, 
and $\sum_i k_i = k$.  Up to renumbering, and up to a subsequence, we may assume that
\begin{equation}\label{HighestWeightSec7}
\mu_{1,\alpha}^1 = \max_{i,j}\mu_{j,\alpha}^i\hskip.1cm .
\end{equation}
Then we let $\mu_\alpha$ stand for $\mu_{1,\alpha}^1$, 
so that $\mu_\alpha$ is the largest weight among all the possible weights $\mu_{j,\alpha}^i$, and we let 
$x_\alpha$ stand for $x_{1,\alpha}^1$, so that $x_\alpha$ is the corresponding 
$x_{j,\alpha}^i$. In other words, 
up to renumbering, and up to a subsequence, we assume (\ref{HighestWeightSec7}), and then let
\begin{equation}\label{HighestWeightBisSec7}
\mu_\alpha = \mu_{1,\alpha}^1 \hskip.2cm\hbox{and}\hskip.2cm x_\alpha = x_{1,\alpha}^1
\end{equation}
for all $\alpha$. 
We let also $\tilde{\mathcal U}_\alpha$ be the $p$-map defined in the Euclidean ball $B_0(1) \subset 
{\mathbb R}^n$ centered at $0$ and of radius $1$ by 
$\tilde {\mathcal U}_\alpha = (\tilde u_\alpha^1,\dots,\tilde u_\alpha^p)$ and
\begin{equation}\label{DefRescaledFctSec7}
\tilde u_\alpha^i(x) = u_\alpha^i\left(\exp_{x_\alpha}\left(\sqrt{\mu_\alpha}x\right)\right)
\end{equation}
for all $\alpha$, where $\exp_{x_\alpha}$ is the exponential map 
at $x_\alpha$, and ${\mathcal U}_\alpha = (u_\alpha^1,\dots,u_\alpha^p)$. The terminology 
harmonic in Theorem \ref{SharpPtwseAsymptThm} below refers to the Euclidean Laplacian. 
Namely $\varphi: \Omega \to {\mathbb R}$ is harmonic, where $\Omega$ 
is an open subset of ${\mathbb R}^n$, 
if $\Delta\varphi = 0$ in $\Omega$, where $\Delta$ is the Euclidean Laplacian in ${\mathbb R}^n$.

\begin{thm}[Sharp Asymptotics]\label{SharpPtwseAsymptThm} Let $(M,g)$ be a smooth compact Riemannian 
manifold of dimension 
$n \ge 3$, $p \ge 1$ an integer, and $\left(A(\alpha)\right)_\alpha$ a sequence of smooth maps 
$A(\alpha): M \to M_p^s({\mathbb R})$ satisfying (\ref{CondSec7}).
Let $({\mathcal U}_\alpha)_\alpha$ 
be a bounded sequence in $H_{1,p}^2(M)$ of nonnegative solutions of (\ref{GenericEqtSec7}) 
which blows up. 
Then there exist $\delta > 0$, nonnegative real numbers $A_i$, and harmonic functions 
$\varphi_i: B_0(\delta) \to {\mathbb R}$, $i = 1,\dots,p$, such that, up to a subsequence, for any $i$,
\begin{equation}\label{SharpPtseEstSec7ThmStat}
\tilde u_\alpha^i(x) \rightarrow \frac{A_i}{\vert x\vert^{n-2}} + \varphi_i(x)
\end{equation}
in $C^1_{loc}\left(B_0(\delta)\backslash\{0\}\right)$ as $\alpha \to +\infty$, where the $\tilde u_\alpha^i$'s 
are given by (\ref{DefRescaledFctSec7}). If the $-A(\alpha)$'s are cooperative for all $\alpha$, then 
at least one of the $A_i$'s is positive. If the $-A(\alpha)$'s are cooperative for all $\alpha$, 
the limit system (\ref{LimitSystSec7}) is fully coupled, and ${\mathcal U}^0 \not\equiv 0$, where 
${\mathcal U}^0$ is as in Theorem \ref{SobTheorThm}, then the $\varphi_i$'s are positive functions 
in $B_0(\delta)$ for all $i$.
\end{thm}

We prove Theorem \ref{SharpPtwseAsymptThm} in several steps (Steps \ref{step1sec7} to 
\ref{step8sec7} below).  
For $p_1$ and $p_2$ such that 
$2^\star/2 < p_2 < 2^\star < p_1$, and $\sigma > 0$, we define the norm $\Vert\cdot\Vert_{p_1,p_2,\sigma}$ 
on $L^\infty(M)$, the space of bounded functions in $M$, by
$$\Vert u\Vert_{p_1,p_2,\sigma} = \inf\Big\{C > 0\hskip.1cm\hbox{s.t.}\hskip.1cm 
(I_{p_1,p_2}^\sigma)\hskip.1cm\hbox{holds for}\hskip.1cm u\Big\}\hskip.1cm ,$$
where $(I_{p_1,p_2}^\sigma)$ holds for $u$ if there exist 
nonnegative functions $u^1, u^2 \in L^\infty(M)$ such that 
$\vert u\vert \le u^1+u^2$,
$\Vert u^1\Vert_{p_1} \le C$, and 
$$\Vert u^2\Vert_{p_2}\le C\sigma^{\frac{n}{2^\star}-\frac{n}{p_2}}\hskip.1cm .$$
For the sake of completeness, we mention the following specific (and easy) result from elliptic theory that 
we will repeatedly use in the proof of Steps \ref{step1sec7} to \ref{step3sec7}. Namely that 
if $u \in H_2^p(M)$ and $f \in L^p(M)$, $p > 1$, are such that 
$\Delta_gu + \lambda u = f$ for some $\lambda > 0$, then $\Vert u\Vert_{H_2^p} \le C \Vert f\Vert_p$ 
where $C > 0$ depends only on the manifold, $p$, and $\lambda$. In order to prove this we write that 
\begin{equation}\label{Sec7EllTheoryEqt}
\Vert u\Vert_{H_2^p} \le C\left(\Vert u\Vert_p + \Vert L_\lambda u\Vert_p\right)\hskip.1cm ,
\end{equation}
where $L_\lambda = \Delta_g + \lambda$, and 
$C > 0$ depends only on the manifold, $p$, and $\lambda$. Equation (\ref{Sec7EllTheoryEqt}) follows from standard 
elliptic theory. Then, by (\ref{Sec7EllTheoryEqt}), it suffices to prove that $\Vert u\Vert_p \le C \Vert L_\lambda u\Vert_p$ 
where $C > 0$ depends only on the manifold, $p$, and $\lambda$. We proceed by contradiction and 
assume that there 
exists a sequence 
$(u_\alpha)$ in $H_2^p(M)$ such that $\Vert u_\alpha\Vert_p = 1$ 
for all $\alpha$ and $\Vert L_\lambda u_\alpha\Vert_p \to 0$ as 
$\alpha \to +\infty$. By (\ref{Sec7EllTheoryEqt}), the sequence $(u_\alpha)$ is bounded in $H_2^p(M)$. 
By compactness of the embedding of $H_2^p$ in $H_1^p$, up to a subsequence, we may 
assume that $u_\alpha \to u$ in $H_1^p$. Then 
$\Vert u\Vert_p = 1$ while $L_\lambda u = 0$. But $L_\lambda u = 0$ implies $u = 0$ since 
$\lambda > 0$. A contradiction and 
this proves that if $u \in H_2^p(M)$ and $f \in L^p(M)$, $p > 1$, are such that 
$\Delta_gu + \lambda u = f$, for some $\lambda > 0$, then $\Vert u\Vert_{H_2^p} \le C \Vert f\Vert_p$ 
where $C > 0$ depends only on the manifold, $p$, and $\lambda$. 

\medskip Step \ref{step1sec7} in the proof of Theorem \ref{SharpPtwseAsymptThm} is as follows. 
In the sequel, once and for all, the $A(\alpha)$'s are assumed to satisfy (\ref{CondSec7}). We will also 
always assume, up to passing to a subsequence, that the ${\mathcal U}_\alpha$'s satisfy (\ref{EqtThmSobTherSec4}) 
and (\ref{EqtThmSobTherSec4Bis}).

\begin{Step}\label{step1sec7} Let $({\mathcal U}_\alpha)_\alpha$ 
be a bounded sequence in $H_{1,p}^2(M)$ of nonnegative solutions of (\ref{GenericEqtSec7}) 
which blows up. 
There exists $p(n) > 2^\star$ with the property that for any $p_1, p_2$ such that
\begin{equation}\label{EqtStatStep1Sec7}
p_0(n) < p_2 < 2^\star < p_1 < p(n)\hskip.1cm ,
\end{equation}
where $p_0(n) = \max(\frac{2n}{n+2},\frac{n}{n-2})$, there exists $C > 0$ such that, up to a subsequence, 
for any $i$, and any $\alpha$,
\begin{equation}\label{EqtStatStep2Sec7}
\Vert u_\alpha^i\Vert_{p_1,p_2,\mu_\alpha^{-1}} \le C\hskip.1cm ,
\end{equation}
where ${\mathcal U}_\alpha = (u_\alpha^1,\dots,u_\alpha^p)$, and $\mu_\alpha$ is as in 
(\ref{HighestWeightSec7})--(\ref{HighestWeightBisSec7}).
\end{Step}

\begin{proof}[Proof of Step \ref{step1sec7}] We fix $\lambda > 0$ and let $G$ be the Green's function of the 
operator $\Delta_g + \lambda$. By (\ref{GenericEqtSec7}) and (\ref{CondSec7}), we may write that
\begin{eqnarray*} u_\alpha^i(x) 
& = & \int_MG(x,y)u_\alpha^i(y)^{2^\star-1}dv_g(y) - 
\sum_{j=1}^p\int_M\left(A_{ij}^\alpha - \lambda\delta_{ij}\right)G(x,y)u_\alpha^j(y)dv_g(y)\\
& \le & \int_MG(x,y)u_\alpha^i(y)^{2^\star-1}dv_g(y) +  
C\sum_{j=1}^p\int_MG(x,y)u_\alpha^j(y)dv_g(y)
\end{eqnarray*}
for all $i = 1,\dots,p$ and $x \in M$, where $C > 0$ is independent of $\alpha$, $x$, and $i$, and where 
the notation $dv_g(y)$ emphasizes that integration is with respect to $y$. By the $H_1^2$-theory, 
Theorem \ref{SobTheorThm}, writing that $u_\alpha^i(y)^{2^\star-1} = u_\alpha^i(y)^{2^\star-2}u_\alpha^i(y)$, 
we then get that
\begin{equation}\label{ProofStep1Sec7Eqt1}
\begin{split} u_\alpha^i(x)
& \le C_1\sum_{j=1}^p\int_MG(x,y)u_\alpha^j(y)dv_g(y)\\
&\hskip.2cm + C_2 \sum_{j=1}^{k_i}\int_MG(x,y)B_{j,\alpha}^i(y)^{2^\star-2}u_\alpha^i(y)dv_g(y)\\
&\hskip.2cm + C_3 \int_MG(x,y)\vert R_\alpha^i(y)\vert^{2^\star-2}u_\alpha^i(y)dv_g(y)
\end{split}
\end{equation}
for all $i = 1,\dots,p$ and $x \in M$, where $C_1, C_2, C_3$ are positive 
constants independent of $\alpha$, $x$, and $i$, and where $k_i$, the $1$-bubbles $B_{j,\alpha}^i$, and 
the remaining terms $R_\alpha^i$ are as in 
(\ref{EqtThmSobTherSec4}) and (\ref{EqtThmSobTherSec4Bis}). Let $w_\alpha$ be the function given by
\begin{equation}\label{ProofStep1Sec71stFct}
w_\alpha(x) = \sum_{j=1}^p\int_MG(x,y)u_\alpha^j(y)dv_g(y)\hskip.1cm .
\end{equation}
Then $\Delta_gw_\alpha + \lambda w_\alpha = \sum_{j=1}^pu_\alpha^j$, 
and since the right hand side in 
this equation  is bounded in $L^{2^\star}(M)$, independently of $\alpha$, it follows from elliptic theory 
that the sequence consisting of the $w_\alpha$'s is bounded in $H_2^{2^\star}(M)$. In particular, by 
the Sobolev embedding theorem, there exists a dimensional 
constant $p(n) > 2^\star$ such that
\begin{equation}\label{ProofStep1Sec7FinalEqt1}
\Vert w_\alpha\Vert_{p_1} \le C
\end{equation}
for all $\alpha$ and all $2^\star \le p_1 \le p(n)$, where $C > 0$ is independent of $\alpha$, 
and $w_\alpha$ is given by (\ref{ProofStep1Sec71stFct}). Similarly, if we let
$w_{j,\alpha}^i$ be the function given by
\begin{equation}\label{ProofStep1Sec72ndFct}
w_{j,\alpha}^i(x) = \int_MG(x,y)B_{j,\alpha}^i(y)^{2^\star-2}u_\alpha^i(y)dv_g(y)\hskip.1cm ,
\end{equation}
then
\begin{equation}\label{ProofStep1Sec7Eqt2ndFcts}
\Delta_gw_{j,\alpha}^i + \lambda w_{j,\alpha}^i  = (B_{j,\alpha}^i)^{2^\star-2}u_\alpha^i\hskip.1cm .
\end{equation}
Let $1 < q < r$ be such that $1/q = (1/r) + (1/2^\star)$. We can write that
\begin{equation}\label{ProofStep1Sec7Eqt2ndFcts1}
\begin{split} \Vert(B_{j,\alpha}^i)^{2^\star-2}u_\alpha^i\Vert_q
&\le \Vert(B_{j,\alpha}^i)^{2^\star-2}\Vert_r\Vert u_\alpha^i\Vert_{2^\star}\\
&\le C \Vert(B_{j,\alpha}^i)^{2^\star-2}\Vert_r\hskip.1cm ,
\end{split}
\end{equation}
where $C > 0$ is independent of $\alpha$, $i$, and $j$. Moreover, by the equation 
(\ref{Def1BubbleSec4}) of a $1$-bubble, we can write that for $\delta > 0$ sufficiently small,
$$\int_M(B_{j,\alpha}^i)^{(2^\star-2)r}dv_g = 
\int_{B_\delta}(B_{j,\alpha}^i)^{(2^\star-2)r}dv_g + 
O\left((\mu_{j,\alpha}^i)^{2r}\right)\hskip.1cm ,$$
where $B_\delta = B_{x_{j,\alpha}^i}(\delta)$, 
and the $x_{j,\alpha}^i$'s and $\mu_{j,\alpha}^i$'s are the centers and weights of the $1$-bubble 
$(B_{j,\alpha}^i)$, while we also have that
$$\int_{B_\delta}(B_{j,\alpha}^i)^{(2^\star-2)r}dv_g \le C (\mu_{j,\alpha}^i)^{n-2r}$$
if $r > n/4$, where $C > 0$ is independent of $\alpha$, $i$, and $j$. In particular, 
if $r \in \left(\frac{n}{4},\frac{n}{2}\right)$, then
\begin{equation}\label{ProofStep1Sec7Eqt2ndFcts2}
\begin{split} \int_M(B_{j,\alpha}^i)^{(2^\star-2)r}dv_g 
&\le C (\mu_{j,\alpha}^i)^{n-2r}\\
&\le C (\mu_\alpha)^{n-2r}\hskip.1cm ,
\end{split}
\end{equation}
where $C > 0$ is independent of $\alpha$, $i$, and $j$, and where $\mu_\alpha$ is as in 
(\ref{HighestWeightSec7})--(\ref{HighestWeightBisSec7}). By elliptic theory, 
(\ref{ProofStep1Sec7Eqt2ndFcts}), (\ref{ProofStep1Sec7Eqt2ndFcts1}), 
and (\ref{ProofStep1Sec7Eqt2ndFcts2}), we then get that the sequence consisting of the 
$w_{j,\alpha}^i$'s is bounded in $H_2^q(M)$ where $1/q = (1/r) + (1/2^\star)$. By the Sobolev embedding 
theorem, since $r$ is arbitrary in $\left(\frac{n}{4},\frac{n}{2}\right)$, it easily follows that
\begin{equation}\label{ProofStep1Sec7FinalEqt2}
\Vert w_{j,\alpha}^i\Vert_{p_2} \le C (\mu_\alpha^{-1})^{\frac{n}{2^\star}-\frac{n}{p_2}}
\end{equation}
for all $p_2 \in \left(\frac{2n}{n+2},2^\star\right)$, $\alpha$, $i$, and $j$, where $C > 0$ 
is independent of $\alpha$, $i$, and $j$, and 
$w_{j,\alpha}^i$ is given by (\ref{ProofStep1Sec72ndFct}). At last, we let $w_\alpha^i$ be the function given by
\begin{equation}\label{ProofStep1Sec73dFct}
w_\alpha^i(x) = \int_MG(x,y)\vert R_\alpha^i(y)\vert^{2^\star-2}u_\alpha^i(y)dv_g(y)\hskip.1cm .
\end{equation}
Then
$$\Delta_gw_\alpha^i + \lambda w_\alpha^i = \vert R_\alpha^i\vert^{2^\star-2}u_\alpha^i$$
and by Step \ref{step2sec7} below, since $R_\alpha^i \to 0$ in $H_1^2$ as $\alpha \to +\infty$, 
we can write that if $2^\star/2 < p_2 < 2^\star < p_1$, then
\begin{equation}\label{ProofStep1Sec7FinalEqt3}
\begin{split} \Vert w_\alpha^i\Vert_{p_1,p_2,\mu_\alpha^{-1}}
& \le C \Vert(R_\alpha^i)^{2^\star-2}\Vert_{n/2}\Vert u_\alpha^i\Vert_{p_1,p_2,\mu_\alpha^{-1}}\\
& = o\left(\Vert u_\alpha^i\Vert_{p_1,p_2,\mu_\alpha^{-1}}\right)
\end{split}
\end{equation}
for all $\alpha$ and $i$, where $w_\alpha^i$ is given 
by (\ref{ProofStep1Sec73dFct}). 
Combining (\ref{ProofStep1Sec7Eqt1}), 
(\ref{ProofStep1Sec7FinalEqt1}), 
(\ref{ProofStep1Sec7FinalEqt2}), and (\ref{ProofStep1Sec7FinalEqt3}), it follows that there exists 
a dimensional constant $p(n) > 2^\star$ such that 
$$\Vert u_\alpha^i\Vert_{p_1,p_2,\mu_\alpha^{-1}} \le C$$
for all $\alpha$ and $i$, and all $p_1$, $p_2$ such that $p_0(n) < p_2 < 2^\star < p_1 < p(n)$, 
where $p_0(n)$ is as in Step \ref{step1sec7}, and $C > 0$ is independent of $\alpha$, and $i$. 
In particular, (\ref{EqtStatStep1Sec7}) and (\ref{EqtStatStep2Sec7}) are true. 
This proves Step \ref{step1sec7}.
\end{proof}

Step \ref{step2sec7} that we used in Step \ref{step1sec7} is as follows.

\begin{Step}\label{step2sec7} Let $u, v \in H_1^2(M)\cap L^\infty(M)$ and $K \in L^\infty(M)$ be nonnegative functions 
such that
$$\Delta_gu + \lambda u \le Kv$$
in $M$, where $\lambda > 0$. 
Let $p_1, p_2$ be such that $2^\star/2 < p_2 < 2^\star < p_1$. Then
$$\Vert u\Vert_{p_1,p_2,\sigma} \le C\Vert K\Vert_{n/2}\Vert v\Vert_{p_1,p_2,\sigma}$$
for all $\sigma > 0$, where $C > 0$ is independent of $u$, $v$, and $\sigma$.
\end{Step}

\begin{proof}[Proof of Step \ref{step2sec7}] Let $\Lambda > \Vert v\Vert_{p_1,p_2,\sigma}$. Then there exist 
$v_1, v_2 \ge 0$ in $L^\infty(M)$ such that $v \le v_1 + v_2$, $\Vert v_1\Vert_{p_1} \le \Lambda$, and 
$\Vert v\Vert_{p_2} \le \Lambda\sigma^{\frac{n}{2^\star}-\frac{n}{p_2}}$. Let $u_1, u_2$ be such that 
\begin{equation}\label{Proofstep2sec7Eqt1}
\begin{split}
&\Delta_gu_1 + \lambda u_1 = K v_1\hskip.1cm,\hskip.1cm\hbox{and}\\
&\Delta_gu_2 + \lambda u_2 = Kv_2\hskip.1cm .
\end{split}
\end{equation}
Then $u_1$ and $u_2$ are nonnegative functions in $H_2^p(M)$ for all $p > 1$, in particular 
$u_1, u_2 \in L^\infty(M)$, and since
$\Delta_g\left(u-\sum u_i\right) + \lambda\left(u-\sum u_i\right) \le 0$, 
we also have that $u \le u_1 + u_2$. Let $q_1$ and $q_2$ be such that $1/q_i = (1/p_i) + (2/n)$, 
$i = 1,2$. Since $p_i > 2^\star/2$, we have that $q_i > 1$, and we can also write that
\begin{equation}\label{Proofstep2sec7Eqt2}
\Vert Kv_i\Vert_{q_i} \le \Vert K\Vert_{n/2}\Vert v_i\Vert_{p_i}\hskip.1cm ,
\end{equation}
for $i = 1,2$. By elliptic theory, (\ref{Proofstep2sec7Eqt1}), 
(\ref{Proofstep2sec7Eqt2}), and the Sobolev embedding theorem 
we then get that
\begin{eqnarray*} \Vert u_i\Vert_{p_i}
& \le & C\Vert u_i\Vert_{H_2^{q_i}}\\
& \le & C\Vert K\Vert_{n/2}\Vert v_i\Vert_{p_i}
\end{eqnarray*}
for $i = 1,2$, where $C > 0$ depends only on the manifold, 
$\lambda$, and the $p_i$'s, $i = 1,2$. Since $u \le u_1 + u_2$, it 
follows that $\Vert u\Vert_{p_1,p_2,\sigma} \le C\Vert K\Vert_{n/2}\Lambda$, and since 
$\Lambda > \Vert v\Vert_{p_1,p_2,\sigma}$ is arbitrary, Step \ref{step2sec7} is proved.
\end{proof}

The next step in the proof of Theorem \ref{SharpPtwseAsymptThm} is as follows. 
In its statement, we let $\theta(n)$ be given by 
$\theta(n) = \frac{n(n+2)}{2(n-2)}$.

\begin{Step}\label{step3sec7} Let $u, v, w \in H_1^2(M)\cap L^\infty(M)$, $u, v \ge 0$, be 
such that
\begin{equation}\label{Step3Sec7Eqt1}
\Delta_gu + \lambda u = v^{2^\star-1} + w
\end{equation}
in $M$, where $\lambda > 0$. 
Let $p_1, p_2$ be such that $2^\star - 1 < p_2 < 2^\star < p_1 < \theta(n)$, and 
$q_1, q_2$ be such that
$\frac{1}{q_i} = \frac{2^\star-1}{p_i} - \frac{2}{n}$,
for $i = 1,2$. Then
\begin{equation}\label{Step3Sec7Eqt2}
\Vert u\Vert_{q_1,q_2,\sigma} \le C
\left(1 + \max\left(\Vert v\Vert_{p_1,p_2,\sigma},\Vert w\Vert_{p_1,p_2,\sigma}\right)^{2^\star-1}\right)
\end{equation}
for all $\sigma > 0$, where $C > 0$ is independent of $u$, $v$, $w$, and $\sigma$.
\end{Step}

\begin{proof}[Proof of Step \ref{step3sec7}] Let 
$\Lambda >  \max\left(\Vert v\Vert_{p_1,p_2,\sigma},\Vert w\Vert_{p_1,p_2,\sigma}\right)$. Then there 
exist nonnegative functions $f_1, f_2, f^\prime_1$, and $f^\prime_2$ in $L^\infty(M)$ such that
\begin{equation}\label{ProofStep3Sec7Eqt1}
\begin{split}
& v \le f_1+f_2\hskip.1cm ,\hskip.1cm \Vert f_1\Vert_{p_1} \le \Lambda\hskip.1cm ,\hskip.1cm 
\Vert f_2\Vert_{p_2} \le \Lambda\sigma^{\frac{n}{2^\star}-\frac{n}{p_2}}\hskip.1cm ,\hskip.1cm\hbox{and}\\
& \vert w\vert \le f^\prime_1+f^\prime_2\hskip.1cm ,\hskip.1cm 
\Vert f^\prime_1\Vert_{p_1} \le \Lambda\hskip.1cm ,\hskip.1cm 
\Vert f^\prime_2\Vert_{p_2} \le \Lambda\sigma^{\frac{n}{2^\star}-\frac{n}{p_2}}\hskip.1cm .
\end{split}
\end{equation}
We let $D = D(n)$, $D > 0$ depending only on $n$, be such that 
\begin{eqnarray*} v^{2^\star-1} + w
&\le& v^{2^\star-1} + \vert w\vert\\
&\le& D\left(\left(f_1+f^\prime_1\right)^{2^\star-1}+1\right) + D\left(f_2+f^\prime_2\right)^{2^\star-1}
\end{eqnarray*}
and let $H_1, H_2$ be given by $H_1 = D\left(\left(f_1+f^\prime_1\right)^{2^\star-1}+1\right)$ and 
$H_2 = D\left(f_2+f^\prime_2\right)^{2^\star-1}$. We also define $u_1$ and $u_2$ by
\begin{equation}\label{Proofstep3sec7Eqt2}
\begin{split}
&\Delta_gu_1 + \lambda u_1 = H_1\hskip.1cm,\hskip.1cm\hbox{and}\\
&\Delta_gu_2 + \lambda u_2 = H_2\hskip.1cm .
\end{split}
\end{equation}
Then $u_1, u_2$ are nonnegative functions in $H_2^p$ for all $p > 1$. 
In particular, we have that $u_1, u_2 \in L^\infty(M)$, and by elliptic theory and 
the Sobolev embedding theorem, we can write that
\begin{equation}\label{Proofstep3sec7Eqt3}
\Vert u_i\Vert_{q_i}
\le C \Vert u_i\Vert_{H_2^{\tilde p_i}} \le C \Vert H_i\Vert_{\tilde p_i}
\end{equation}
for $i = 1,2$, where $1/q_i = (1/\tilde p_i) - (2/n)$, 
$\tilde p_i = p_i/(2^\star-1)$, and 
$C > 0$ depends only on the manifold, $\lambda$ and the $p_i$'s, $i = 1,2$. By (\ref{Step3Sec7Eqt1}),
$$\Delta_gu + \lambda u \le H_1 + H_2\hskip.1cm ,$$
and it follows from (\ref{Proofstep3sec7Eqt2}) and the maximum principle that $u \le u_1+u_2$. As is easily checked,
$$(2^\star-1)\left(\frac{n}{2^\star}-\frac{n}{p_2}\right) = \frac{n}{2^\star} - \frac{n}{q_2}\hskip.1cm .$$
It follows that
$\Vert H_1\Vert_{\tilde p_1} \le C\left(1 + \Lambda^{2^\star-1}\right)$ and 
$\Vert H_2\Vert_{\tilde p_2} \le C\Lambda^{2^\star-1}\sigma^{\frac{n}{2^\star}-\frac{n}{q_2}}$, 
where $C > 0$ is independent of $u$, $v$, $w$, and $\sigma$. Since 
$\Lambda >  \max\left(\Vert v\Vert_{p_1,p_2,\sigma},\Vert w\Vert_{p_1,p_2,\sigma}\right)$ is arbitrary, 
and since $u \le u_1 + u_2$, we get with (\ref{Proofstep3sec7Eqt3}) that 
(\ref{Step3Sec7Eqt2}) is true. This proves Step \ref{step3sec7}.
\end{proof}

The next step in the proof of Theorem \ref{SharpPtwseAsymptThm} is as follows. Estimates 
like in Step \ref{step4sec7}, where the norm $\Vert\cdot\Vert_{p_1,p_2,\sigma}$ is involved, 
have been introduced in Devillanova and Solimini \cite{DevSol}. 

\begin{Step}\label{step4sec7} Let $({\mathcal U}_\alpha)_\alpha$ 
be a bounded sequence in $H_{1,p}^2(M)$ of nonnegative solutions of (\ref{GenericEqtSec7}) 
which blows up. Let $p_1, p_2$ be such that $2^\star/2 < p_2 < 2^\star < p_1$. There exist 
$C > 0$ and sequences $(v_{1,\alpha}^i)_\alpha$, $(v_{2,\alpha}^i)_\alpha$ of nonnegative functions 
in $L^\infty(M)$ such that, up to a subsequence, $u_\alpha^i \le v_{1,\alpha}^i + v_{2,\alpha}^i$, 
$\Vert v_{1,\alpha}^i\Vert_{p_1} \le C$, and
\begin{equation}\label{Step4Sec7StateEqt}
\Vert v_{2,\alpha}^i\Vert_{p_2} \le C \mu_\alpha^{\frac{n}{p_2}-\frac{n}{2^\star}}
\end{equation}
for all $i$ and $\alpha$, 
where ${\mathcal U}_\alpha = (u_\alpha^1,\dots,u_\alpha^p)$, and $\mu_\alpha$ is as in 
(\ref{HighestWeightSec7})--(\ref{HighestWeightBisSec7}).
\end{Step}

\begin{proof}[Proof of Step \ref{step4sec7}] We prove Step \ref{step4sec7} by induction, starting from 
Step \ref{step1sec7}, using Step \ref{step3sec7}. 
An easy remark is that
\begin{equation}\label{ProofStep4Sec7Eqt1}
\Vert u\Vert_{\tilde p_1,p_2,\sigma} \le C \Vert u\Vert_{p_1,p_2,\sigma}
\end{equation}
if $\tilde p_1 \le p_1$, where $C > 0$ 
depends only on the manifold. We fix $p_1, p_2$ such that $\frac{2^\star}{2} < p_2 < 2^\star < p_1$. We let 
$p_1^0 > 2^\star$ be close to $2^\star$, and let $k_0 \ge 1$ be such that the increasing sequence 
$(p_1^k)_k$ given by
$$\frac{1}{p_1^{k+1}} = \frac{2^\star-1}{p_1^k} - \frac{2}{n}$$
satisfies $p_1^k < \theta(n)$ for all $k \le k_0$, and $p_1^{k_0+1} \ge \theta(n)$, where $\theta(n)$ is 
as in Step \ref{step3sec7}. Similarly, for $p_2^0 < 2^\star$ we construct the decreasing sequence $(p_2^k)_k$ 
by
$$\frac{1}{p_2^{k+1}} = \frac{2^\star-1}{p_2^k} - \frac{2}{n}\hskip.1cm .$$
We choose $p_2^0$ such that $p_2^{k_0+2} = p_2$. Then, since $p_2 > 2^\star/2$, we get that 
$p_2^k > 2^\star-1$ for all $k \le k_0+1$. 
The more $p_1^0 > 2^\star$ is close to $2^\star$, the more $k_0$ is large, and the more $k_0$ is large, the 
more $p_2^0 < 2^\star$ has to be close to $2^\star$. In particular, we can assume that $p_2^0 > 2^\star/(2^\star-1)$. 
We fix $\lambda > 0$, and, by (\ref{GenericEqtSec7}), we write that
$$\Delta_gu_\alpha^i + \lambda u_\alpha^i = (u_\alpha^i)^{2^\star-1} 
+ \sum_{j=1}^p\tilde A^\alpha_{ij}u_\alpha^j$$
for all $i$ and $\alpha$, where $\tilde A^\alpha_{ij} = \lambda \delta_{ij}-A^\alpha_{ij}$. In particular, 
by  (\ref{CondSec7}), there exists $C > 0$ such that 
$\Vert\tilde A^\alpha_{ij}\Vert_{C^0} \le C$ for all $i$, $j$, and $\alpha$. 
Then, by Steps \ref{step1sec7} and \ref{step3sec7}, we can write that there exists $C > 0$ such that, up to a subsequence, 
$$\Vert u_\alpha^i\Vert_{p_1^{k_0+1},p_2^{k_0+1},\mu_\alpha^{-1}} \le C$$
for all $i$ and $\alpha$.
In particular, by (\ref{ProofStep4Sec7Eqt1}),
$\Vert u_\alpha^i\Vert_{\tilde p_1,p_2^{k_0+1},\mu_\alpha^{-1}} \le C$
for $\tilde p_1 < \theta(n)$ as close as we want to $\theta(n)$. We then apply Step \ref{step3sec7} once more and get that 
$$\Vert u_\alpha^i\Vert_{\hat p_1,p_2^{k_0+2},\mu_\alpha^{-1}} \le C\hskip.1cm ,$$
where $\hat p_1 \to +\infty$ as 
$\tilde p_1 \to \theta(n)$. Choosing $\tilde p_1$ sufficiently close to $\theta(n)$, we can assume that $\hat p_1 \ge p_1$, 
and, thanks to (\ref{ProofStep4Sec7Eqt1}), we get that there exist 
$C > 0$ and sequences $(v_{1,\alpha}^i)_\alpha$, $(v_{2,\alpha}^i)_\alpha$ of nonnegative functions 
in $L^\infty(M)$ such that, up to a subsequence, $u_\alpha^i \le v_{1,\alpha}^i + v_{2,\alpha}^i$, 
$\Vert v_{1,\alpha}^i\Vert_{p_1} \le C$, and (\ref{Step4Sec7StateEqt}) holds. 
This proves Step \ref{step4sec7}.
\end{proof}

Going on with the proof of Theorem \ref{SharpPtwseAsymptThm}, we now 
prove integral estimates for sequences 
of nonnegative solutions of (\ref{GenericEqtSec7}). Step \ref{step5sec7} is as follows.

\begin{Step}\label{step5sec7} Let $({\mathcal U}_\alpha)_\alpha$ 
be a bounded sequence in $H_{1,p}^2(M)$ of nonnegative solutions of (\ref{GenericEqtSec7}) 
which blows up. There exist $C_1, C_2 > 0$ such that, up to a subsequence,
\begin{equation}\label{Step5Sec7StatemEqt}
\frac{1}{r^{n-1}} \int_{\partial B_{x_\alpha}(r)} u_\alpha^id\sigma_g \le C_1 
+ C_2 \frac{\mu_\alpha^{\frac{n-2}{2}}}{r^{n-2}}
\end{equation}
for all $i$, all $\alpha$, and all $r > 0$ sufficiently small, 
where ${\mathcal U}_\alpha = (u_\alpha^1,\dots,u_\alpha^p)$, the $x_\alpha$'s and $\mu_\alpha$'s are as in 
(\ref{HighestWeightSec7})--(\ref{HighestWeightBisSec7}), and $d\sigma_g$ is the measure induced on 
$\partial B_{x_\alpha}(r)$ by $g$.
\end{Step}

\begin{proof}[Proof of Step \ref{step5sec7}] Thanks, for instance, to the Bishop-Gromov comparison 
theorem (see Chavel \cite{Cha}), there exists $C > 0$ such that 
$\hbox{Vol}_g\left(B_{x_\alpha}(r)\right) \le Cr^n$ for all $\alpha$ 
and $r > 0$ small, where $\hbox{Vol}_g\left(B_{x_\alpha}(r)\right)$ is the volume 
of the ball $B_{x_\alpha}(r)$. By 
(\ref{GenericEqtSec7}) and (\ref{CondSec7}) we may then write that 
\begin{equation}\label{Proofstep5sec7Eqt1}
\begin{split} \int_{B_{x_\alpha}(r)}(\Delta_gu_\alpha^i)dv_g 
&= \int_{B_{x_\alpha}(r)}(u_\alpha^i)^{2^\star-1}dv_g - 
\sum_{j=1}^p\int_{B_{x_\alpha}(r)}A^\alpha_{ij}u_\alpha^jdv_g\\
&\le  C\sum_{i=1}^p\int_{B_{x_\alpha}(r)}\left((u_\alpha^i)^{2^\star-1}+C\right)dv_g\\
&\le Cr^n + C\sum_{i=1}^p\int_{B_{x_\alpha}(r)}(u_\alpha^i)^{2^\star-1}dv_g
\end{split}
\end{equation}
for all $\alpha$, $i$, and all $r > 0$ small, where 
the above constants $C > 0$ are independent of $\alpha$, $i$, and $r$. By Step \ref{step4sec7} that we apply to the 
${\mathcal U}_\alpha$'s with $p_1 = (2^\star-1)n$ and $p_2 = 2^\star-1$, there exist  
$C > 0$ and sequences $(v_{1,\alpha}^i)_\alpha$, $(v_{2,\alpha}^i)_\alpha$ of nonnegative functions 
in $L^\infty(M)$ such that, up to a subsequence, $u_\alpha^i \le v_{1,\alpha}^i + v_{2,\alpha}^i$, 
$\Vert v_{1,\alpha}^i\Vert_{p_1} \le C$, and 
$\Vert v_{2,\alpha}^i\Vert_{p_2} \le C \mu_\alpha^{(n/p_2)-(n/2^\star)}$
for all $i$ and $\alpha$. In particular, by H\"older's inequality, 
\begin{equation}\label{Proofstep5sec7Eqt2}
\begin{split} \int_{B_{x_\alpha}(r)}(u_\alpha^i)^{2^\star-1}dv_g 
&\le C \int_{B_{x_\alpha}(r)}(v_{1,\alpha}^i)^{2^\star-1}dv_g + 
C \int_{B_{x_\alpha}(r)}(v_{2,\alpha}^i)^{2^\star-1}dv_g\\
&\le C \Vert v_{1,\alpha}^i\Vert_{p_1}^{2^\star-1}
\hbox{Vol}_g\left(B_{x_\alpha}(r)\right)^{\frac{n-1}{n}} 
+ C \left(\mu_\alpha^{\frac{n}{p_2}-\frac{n}{2^\star}}\right)^{2^\star-1}
\end{split}
\end{equation}
for all $\alpha$, $i$, and all $r > 0$ small, where 
the above constants $C > 0$ are, here again, independent of $\alpha$, $i$, and $r$. Combining 
(\ref{Proofstep5sec7Eqt1}) and (\ref{Proofstep5sec7Eqt2}), we get that there exists 
$C_1, C_2 > 0$ such that
\begin{equation}\label{Proofstep5sec7Eqt3}
\int_{B_{x_\alpha}(r)}(\Delta_gu_\alpha^i)dv_g \le C_1r^{n-1} + C_2\mu_\alpha^{\frac{n-2}{2}}
\end{equation}
for all $\alpha$, $i$, and all $r > 0$ small. Given $x_0 \in M$, there exists $\beta_{x_0}$ a smooth function around $x_0$ 
such that for $u$ smooth in $M$, and $r > 0$ small (less than the injectivity radius of the manifold),
\begin{equation}\label{Proofstep5sec7Eqt4}
\begin{split}
&\frac{d}{dr} \left(\frac{1}{r^{n-1}}\int_{\partial B_{x_0}(r)}ud\sigma_g\right)\\
&= \frac{1}{r^{n-1}}\int_{\partial B_{x_0}(r)}\left(\frac{\partial u}{\partial\nu}\right) d\sigma_g 
+ \frac{1}{r^{n-1}}\int_{\partial B_{x_0}(r)}\beta_{x_0}ud\sigma_g\hskip.1cm ,
\end{split}
\end{equation}
where $\partial B_{x_0}(r)$ is the boundary of the geodesic ball $B_{x_0}(r)$, where $d\sigma_g$ is the volume 
element on $\partial B_{x_0}(r)$ induced by $g$, and $\frac{\partial}{\partial\nu}$ is the normal derivative 
with respect to the outward unit normal vector $\nu$. As is well known, see again Chavel \cite{Cha}, 
$\beta_{x_0}(x) = O^\prime\left(d_g(x_0,x)\right)$ where 
the notation in the right hand side of this equation 
stands for a $C^1$-function such that the $k$th derivatives of this function, $k = 0,1$, are bounded 
by $Cd_g(x_0,x)^{1-k}$ where $C > 0$ does not depend on $x_0$ and $x$. 
For $i = 1,\dots,p$, we define $\Phi_\alpha^i: (0,\delta) \to {\mathbb R}$, 
$\delta > 0$ small, by
$$\Phi_\alpha^i(r) = \frac{1}{r^{n-1}} 
\int_{\partial B_{x_\alpha}(r)}u_\alpha^idv_g\hskip.1cm .$$
Then, by (\ref{Proofstep5sec7Eqt4}), we can write that
\begin{equation}\label{Proofstep5sec7Eqt5}
\frac{d\Phi_\alpha^i}{dr}(r) = 
-\frac{1}{r^{n-1}} \int_{B_{x_\alpha}(r)}(\Delta_gu_\alpha^i)dv_g + H_\alpha^i(r)\Phi_\alpha^i(r)
\end{equation}
for all $\alpha$, $i$, and all $r > 0$ small, where the $H_\alpha^i$'s are uniformly bounded with respect to 
$\alpha$ and $r$. By (\ref{COConvOutside}), there exists $\delta > 0$, $\delta$ arbitrarily small, with the 
property that, for any $i$, the $\Phi_\alpha^i(\delta)$'s are bounded uniformly with respect to $\alpha$. 
Integrating (\ref{Proofstep5sec7Eqt5}) between $r$ and $\delta$, $0 < r  < \delta$, we get that 
$$-\left[e^{-\int_0^tH_\alpha^i(s)ds}\Phi_\alpha^i(t)\right]_r^\delta
= \int_r^\delta\left(\frac{e^{-\int_0^tH_\alpha^i(s)ds}}{t^{n-1}}
\int_{B_{x_\alpha}(t)}(\Delta_gu_\alpha^i)dv_g\right)dt\hskip.1cm ,$$
and then, it follows from (\ref{Proofstep5sec7Eqt3}) that
\begin{eqnarray*} \Phi_\alpha^i(r)
&\le& C e^{-\int_0^rH_\alpha^i(s)ds}\Phi_\alpha^i(r)\\
&\le& Ce^{-\int_0^\delta H_\alpha^i(s)ds}\Phi_\alpha^i(\delta) 
+ C \int_r^\delta\left(C_1 + C_2\mu_\alpha^{\frac{n-2}{2}}t^{1-n}\right)dt\\
&\le& C_3 + C_4\frac{1}{r^{n-2}} \mu_\alpha^{\frac{n-2}{2}}
\end{eqnarray*}
for all $\alpha$, $i$, and all $0 < r < \delta$, where $C_1, C_2$ are as in (\ref{Proofstep5sec7Eqt3}), 
and the constants $C$, $C_3$, and $C_4$ are positive and independent of $\alpha$, $i$, and $r$. 
This proves (\ref{Step5Sec7StatemEqt}) and Step \ref{step5sec7}.
\end{proof}

Now that we have Step \ref{step5sec7}, we can prove the first part of Theorem \ref{SharpPtwseAsymptThm}. 
This is the subject of the following step.

\begin{Step}\label{step6sec7} Let $({\mathcal U}_\alpha)_\alpha$ 
be a bounded sequence in $H_{1,p}^2(M)$ of nonnegative solutions of (\ref{GenericEqtSec7}) 
which blows up. There exist $\delta > 0$, nonnegative real numbers $A_i$, and harmonic functions 
$\varphi_i: B_0(\delta) \to {\mathbb R}$, $i = 1,\dots,p$, such that, up to a subsequence, for any $i$,
\begin{equation}\label{SharpPtseEstSec7}
\tilde u_\alpha^i(x) \rightarrow \frac{A_i}{\vert x\vert^{n-2}} + \varphi_i(x)
\end{equation}
in $C^1_{loc}\left(B_0(\delta)\backslash\{0\}\right)$ as $\alpha \to +\infty$, where the $\tilde u_\alpha^i$'s 
are given by (\ref{DefRescaledFctSec7}).
\end{Step}

\begin{proof}[Proof of Step \ref{step6sec7}] We let $g_\alpha$ be given by 
$g_\alpha(x) = \left(\exp_{x_\alpha}^\star g\right)(\sqrt{\mu_\alpha}x)$, where $x \in {\mathbb R}^n$, and 
the $x_\alpha$'s and $\mu_\alpha$'s are as in 
(\ref{HighestWeightSec7})--(\ref{HighestWeightBisSec7}). Clearly, 
$g_\alpha \to \xi$ in $C^2(K)$ for any compact subset $K$ of 
${\mathbb R}^n$, where $\xi$ is the Euclidean metric. We also have that
\begin{equation}\label{ProofStep6Sec7Eqt1}
\Delta_{g_\alpha}\tilde u_\alpha^i + \mu_\alpha\sum_{j=1}^p\tilde A^\alpha_{ij}\tilde u_\alpha^j 
= \mu_\alpha(\tilde u_\alpha^i)^{2^\star-1}
\end{equation}
for all $i$, where $\tilde A^\alpha_{ij}(x) = A^\alpha_{ij}\left(\exp_{x_\alpha}(\sqrt{\mu_\alpha}x)\right)$. 
If the $x_{j,\alpha}^i$'s are the centers 
of the $1$-bubbles $(B_{j,\alpha}^i)_\alpha$ in (\ref{EqtThmSobTherSec4Bis}) 
from which the 
$p$-bubbles $({\mathcal B}_{j,\alpha})_\alpha$'s in (\ref{EqtThmSobTherSec4}) are defined, up to 
passing to a subsequence, we may assume that there exists $C > 0$ such that 
for any $i, j$, either $d_g(x^i_{j,\alpha},x_\alpha) \le C\sqrt{\mu_\alpha}$ for all $\alpha$, or 
$d_g(x^i_{j,\alpha},x_\alpha)/\sqrt{\mu_\alpha} \to +\infty$ as $\alpha \to +\infty$. If ${\mathcal A}$ 
is the finite set consisting of the $i, j$'s such that  $d_g(x^i_{j,\alpha},x_\alpha) \le C\sqrt{\mu_\alpha}$ for all $\alpha$, 
we let $\hat{\mathcal S}$ be given by
\begin{equation}\label{ProofStep6Sec7Eqt2}
\hat{\mathcal S} = \left\{\lim_{\alpha \to +\infty}\frac{1}{\sqrt{\mu_\alpha}} \exp_{x_\alpha}^{-1}(x^i_{j,\alpha})
\hskip.1cm ,\hskip.1cm i, j \in {\mathcal A}\right\}
\end{equation}
where, up to passing to another subsequence, we assume that the limits in $\hat{\mathcal S}$ 
exist. Clearly, $0 \in \hat{\mathcal S}$. 
Given $0 < \delta < R$, we let $K = B_0^\prime(R)\backslash\bigcup_{x \in \hat{\mathcal S}}B_x(\delta)$, 
where $B^\prime_0(R)$ stands for the closed Euclidean ball of radius $R$ centered at $0$. We let also 
$h_\alpha^i$ be given by
\begin{equation}\label{ProofStep6Sec7Eqt3}
h_\alpha^i = \mu_\alpha (\tilde u_\alpha^i)^{2^\star-2}\hskip.1cm .
\end{equation}
Then, by (\ref{ProofStep6Sec7Eqt1}), we can write that
\begin{equation}\label{ProofStep6Sec7Eqt4}
\Delta_{g_\alpha}\tilde u_\alpha^i + \mu_\alpha\sum_{j=1}^p\tilde A^\alpha_{ij}\tilde u_\alpha^j 
= h_\alpha^i \tilde u_\alpha^i
\end{equation}
for all $i$, where the $\tilde A^\alpha_{ij}$'s are as in (\ref{ProofStep6Sec7Eqt1}), and the 
$h_\alpha^i$'s are given by (\ref{ProofStep6Sec7Eqt3}). By the estimate 
(\ref{PointEstThmEqt}) in Theorem \ref{PointEstThm}, that we apply to 
$\tilde x_\alpha = \exp_{x_\alpha}(\sqrt{\mu_\alpha}x)$, we easily get that there exists $C > 0$ such that 
for any $i$, and any $x \in K$,
$\vert h_\alpha^i(x)\vert \le C$. We claim that by Lemma \ref{LemRefPtEst}, we actually have that
\begin{equation}\label{ProofStep6Sec7Eqt5}
h_\alpha^i \to 0\hskip.2cm\hbox{in}\hskip.1cm L^\infty(K)
\end{equation}
for all $i$, as $\alpha \to +\infty$. In order to prove (\ref{ProofStep6Sec7Eqt5}), we first note that 
there exists $C > 0$, $C = C(K)$, such that 
\begin{equation}\label{ProofStep6Sec7Eqt6}
d_g\left(x^i_{j,\alpha},\exp_{x_\alpha}(\sqrt{\mu_\alpha}x)\right) \ge 
C \sqrt{\mu_\alpha}
\end{equation}
for all $\alpha$, $i, j$, and $x \in K$. Then, since $\mu_\alpha$ is the largest 
weight among all the possible weights $\mu_{j,\alpha}^i$, we get with 
(\ref{ProofStep6Sec7Eqt6}) that
\begin{equation}\label{ProofStep6Sec7Eqt7}
\frac{1}{\mu^i_{j,\alpha}}d_g\left(x^i_{j,\alpha},\exp_{x_\alpha}(\sqrt{\mu_\alpha}x)\right) \to +\infty 
\end{equation}
in $L^\infty(K)$ for all $i, j$, as $\alpha \to +\infty$. We may therefore apply Lemma \ref{LemRefPtEst} to the 
$\tilde x_\alpha = \exp_{x_\alpha}(\sqrt{\mu_\alpha}x)$, and we get with this lemma and 
(\ref{ProofStep6Sec7Eqt6}) that $h_\alpha^i \to 0$ in $L^\infty(K)$ 
for all $i$, as $\alpha \to +\infty$. This proves  the above claim that (\ref{ProofStep6Sec7Eqt5}) 
is true. Now that we proved (\ref{ProofStep6Sec7Eqt5}), we claim that for any 
$0 < \delta_1 < \delta_2$, and any $p \in \bigl(1,\frac{n}{n-2}\bigr)$, there exists $C = C(\delta_1,\delta_2,p)$, 
$C > 0$, such that
\begin{equation}\label{ProofStep6Sec7Eqt8}
\int_{R(\delta_1,\delta_2)}(\tilde u_\alpha^i)^pdv_{g_\alpha} \le C
\end{equation}
for all $\alpha$ and $i$, where $R(\delta_1,\delta_2)$ is the Euclidean annulus 
centered at $0$ and of radii $\delta_1, \delta_2$. In order to prove 
(\ref{ProofStep6Sec7Eqt8}), we fix $0 < \delta_1 < \delta_2$, $p \in \bigl(1,\frac{n}{n-2}\bigr)$, 
and let $A(\alpha,\delta_1,\delta_2)$ be the annulus in $M$ centered at 
$x_\alpha$ and of radii $\delta_1\sqrt{\mu_\alpha}, \delta_2\sqrt{\mu_\alpha}$. 
By the integral estimate (\ref{Step5Sec7StatemEqt}) in Step \ref{step5sec7}, we can write that
\begin{equation}\label{ProofStep6Sec7Eqt9}
\frac{1}{\hbox{Vol}_g\left(A_n(\alpha)\right)}\int_{A_n(\alpha)}u_\alpha^idv_g \le C
\end{equation}
for all $\alpha$ and $i$, where $A_n(\alpha) = A(\alpha,\delta_1,\delta_2)$, 
$\hbox{Vol}_g\left(A_n(\alpha)\right)$ is the volume of $A_n(\alpha)$ with respect to $g$, and
$C > 0$ is independent of $\alpha$ and $i$. 
By (\ref{GenericEqtSec7}) and (\ref{CondSec7}), like in 
(\ref{Proofstep5sec7Eqt1}), we can also write that
\begin{equation}\label{ProofStep6Sec7Eqt10}
\int_{A_n(\alpha)}\vert\Delta_gu_\alpha^i\vert dv_g 
\le C\hbox{Vol}_g\left(A_n(\alpha)\right) + 
C\sum_{k=1}^p\int_{A_n(\alpha)}(u_\alpha^k)^{2^\star-1}dv_g
\end{equation}
for all $\alpha$ and $i$, where $C > 0$ is independent of $\alpha$ and $i$.  
By Step \ref{step4sec7} that we apply to the 
${\mathcal U}_\alpha$'s with $p_1 = (2^\star-1)n$ and $p_2 = 2^\star-1$, there exist  
$C > 0$ and sequences $(v_{1,\alpha}^i)_\alpha$, $(v_{2,\alpha}^i)_\alpha$ of nonnegative functions 
in $L^\infty(M)$ such that, up to a subsequence, $u_\alpha^i \le v_{1,\alpha}^i + v_{2,\alpha}^i$, 
$\Vert v_{1,\alpha}^i\Vert_{p_1} \le C$, and 
$\Vert v_{2,\alpha}^i\Vert_{p_2} \le C \mu_\alpha^{(n/p_2)-(n/2^\star)}$
for all $i$ and $\alpha$. In particular, by H\"older's inequality, 
\begin{equation}\label{ProofStep6Sec7Eqt11}
\begin{split} \int_{A_n(\alpha)}(u_\alpha^i)^{2^\star-1}dv_g 
&\le C \int_{A_n(\alpha)}(v_{1,\alpha}^i)^{2^\star-1}dv_g + 
C \int_{A_n(\alpha)}(v_{2,\alpha}^i)^{2^\star-1}dv_g\\
&\le C \Vert v_{1,\alpha}^i\Vert_{p_1}^{2^\star-1}
\hbox{Vol}_g\left(A_n(\alpha)\right)^{\frac{n-1}{n}} 
+ C \left(\mu_\alpha^{\frac{n}{p_2}-\frac{n}{2^\star}}\right)^{2^\star-1}\\
&\le C \hbox{Vol}_g\left(A_n(\alpha)\right)^{\frac{n-1}{n}} 
+ C \mu_\alpha^{\frac{n-2}{2}}
\end{split}
\end{equation}
for all $\alpha$ and $i$, where 
the above constants $C > 0$ are, here again, independent of $\alpha$ and $i$. Combining 
(\ref{ProofStep6Sec7Eqt10}) and (\ref{ProofStep6Sec7Eqt11}), we get that there exists 
$C > 0$ such that
\begin{equation}\label{ProofStep6Sec7Eqt12}
\frac{1}{\hbox{Vol}_g\left(A_n(\alpha)\right)}\int_{A_n(\alpha)}\vert\Delta_gu_\alpha^i\vert dv_g \le 
\frac{C}{\mu_\alpha}
\end{equation}
for all $\alpha$ and $i$. Then, by (\ref{ProofStep6Sec7Eqt9}) and (\ref{ProofStep6Sec7Eqt12}), 
we can write that 
\begin{equation}\label{ProofStep6Sec7Eqt13}
\int_{R(\delta_1,\delta_2)}\tilde u_\alpha^i dv_{g_\alpha} \le C \hskip.2cm\hbox{and}\hskip.2cm
\int_{R(\delta_1,\delta_2)}\vert\Delta_{g_\alpha}\tilde u_\alpha^i\vert dv_{g_\alpha} \le C
\end{equation}
for all $\alpha$ and $i$, where $C > 0$ is independent of $\alpha$ and $i$. 
We let $F_\alpha^i$ be such that 
$F_\alpha^i = \Delta_{g_\alpha}\tilde u_\alpha^i$ in $R(\delta_1,\delta_2)$ and $F_\alpha^i = 0$ 
outside $R(\delta_1,\delta_2)$. 
Given $\delta > \delta_2$ we let also $G_\alpha$ be the Green's function of $\Delta_{g_\alpha}$ in $B_0(\delta)$ with 
zero Dirichlet boundary condition, and set
$$v_\alpha^i(x) = \int_{B_0(\delta)}G_\alpha(x,y)F_\alpha^i(y)dv_{g_\alpha}(y)\hskip.1cm .$$
By standard properties of the Green's function, there exists $C > 0$ such that
\begin{equation}\label{ProofStep6Sec7Eqt14}
G_\alpha(x,y) \le \frac{C}{\left\vert y-x\right\vert^{n-2}}
\end{equation}
for all $x \in R(\delta_1,\delta_2)$, all $y \in B_0(\delta)$, 
and all $\alpha$. 
Given $p \in \bigl(1,\frac{n}{n-2}\bigr)$, we let 
$q$ be such that $\frac{1}{p} + \frac{1}{q} = 1$. 
For $\varphi \in L^q\left(R(\delta_1,\delta_2)\right)$, by 
(\ref{ProofStep6Sec7Eqt14}), we can write that
$$\left\vert\int_{R(\delta_1,\delta_2)}v_\alpha^i\varphi dx\right\vert 
\le C \int_{R(\delta_1,\delta_2)}\left(\int_{R(\delta_1,\delta_2)}\frac{\varphi(x)}{\vert y-x\vert^{n-2}} dx\right)
\left\vert F_\alpha^i(y)\right\vert dy\hskip.1cm .$$
By H\"older's inequality, since $p < \frac{n}{n-2}$, we then get that
\begin{eqnarray*}
&&\left\vert\int_{R(\delta_1,\delta_2)}v_\alpha^i\varphi dx\right\vert\\
&&\hskip.4cm \le  C \Vert\varphi\Vert_{L^q\left(R(\delta_1,\delta_2)\right)} 
\int_{R(\delta_1,\delta_2)}\left(\int_{R(\delta_1,\delta_2)}\frac{dx}{\vert y-x\vert^{p(n-2)}}\right)^{1/p}
\left\vert F_\alpha^i(y)\right\vert dy\\
&&\hskip.4cm \le C \Vert\varphi\Vert_{L^q\left(R(\delta_1,\delta_2)\right)} 
\left\Vert F_\alpha^i\right\Vert_{L^1\left(R(\delta_1,\delta_2)\right)}
\end{eqnarray*}
and, by (\ref{ProofStep6Sec7Eqt13}), it follows that
$$\left\vert\int_{R(\delta_1,\delta_2)}v_\alpha^i\varphi dx\right\vert 
\le C\Vert\varphi\Vert_{L^q\left(R(\delta_1,\delta_2)\right)}$$
for all $\alpha$ and $i$, 
where $C > 0$ does not depend on $\alpha$, $i$, and $\varphi$. 
Thus, by duality, taking $\varphi = (v_\alpha^i)^{p-1}$, we get that 
\begin{equation}\label{ProofStep6Sec7Eqt15}
\int_{R(\delta_1,\delta_2)}(v_\alpha^i)^pdv_{g_\alpha} \le C,
\end{equation}
for all $\alpha$ and $i$, where $C > 0$ is independent of $\alpha$ and $i$. 
Since $\Delta_{g_\alpha}\left(v_\alpha^i - \tilde u_\alpha^i\right) = 0$
in $R(\delta_1,\delta_2)$, it follows from 
the De Giorgi-Nash-Moser 
iterative scheme 
that if $\Omega \subset\subset R(\delta_1,\delta_2)$, then
$$\sup_\Omega\left\vert v_\alpha^i - \tilde u_\alpha^i\right\vert \le 
C \left\Vert v_\alpha^i - \tilde u_\alpha^i\right\Vert_{L^1\left(R(\delta_1,\delta_2)\right)}\hskip.1cm ,$$
where $C > 0$ is independent of $\alpha$ and $i$. 
By (\ref{ProofStep6Sec7Eqt13}) and (\ref{ProofStep6Sec7Eqt15}), 
and since $\delta_1 < \delta_2$ are arbitrary, this implies 
the above claim that (\ref{ProofStep6Sec7Eqt8}) is true. In particular, combining 
(\ref{ProofStep6Sec7Eqt4}), (\ref{ProofStep6Sec7Eqt5}), and (\ref{ProofStep6Sec7Eqt8}), we get 
with standard elliptic theory (see, for instance, Theorem 9.11 
in Gilbarg-Trudinger \cite{GilTru}) that the $\tilde u_\alpha^i$'s are 
uniformly bounded in $H_2^p(\Omega)$ for all $p \in \bigl(1,\frac{n}{n-2}\bigr)$, all $i$, 
and all $\Omega \subset\subset {\mathbb R}^n\backslash\hat{\mathcal S}$. By standard bootstrap 
arguments, it follows that the $\tilde u_\alpha^i$'s are 
uniformly bounded in $H_2^p(\Omega)$ for all $p > 1$, all $i$, 
and all $\Omega \subset\subset {\mathbb R}^n\backslash\hat{\mathcal S}$. By the Sobolev embedding theorem 
we may then assume that, up to a subsequence, and for all $i$, the $\tilde u_\alpha^i$'s converge 
in $C^1_{loc}({\mathbb R}^n\backslash\hat{\mathcal S})$ to some nonnegative 
$\tilde u_i$ as $\alpha \to +\infty$. By (\ref{ProofStep6Sec7Eqt4}), and (\ref{ProofStep6Sec7Eqt5}), 
see also (\ref{CondSec7}), the $\tilde u_i$'s are harmonic in 
${\mathbb R}^n\backslash\hat{\mathcal S}$. Let us now write that $\hat{\mathcal S} = \left\{x_1,\dots,x_m\right\}$, 
where $x_1 = 0$. 
By classical results in harmonic analysis, see, for instance, V\'eron \cite{Ver}, since 
the $\tilde u_i$'s are both harmonic and nonnegative in ${\mathbb R}^n\backslash\hat{\mathcal S}$, 
we can write that for any $i$, there exists 
a harmonic function $\psi_i: {\mathbb R}^n \to {\mathbb R}$, and real numbers $A_j^i \in {\mathbb R}$, such that 
\begin{equation}\label{ProofStep6Sec7Eqt16}
\tilde u_i(x) = \sum_{j=1}^m\frac{A_j^i}{\vert x-x_j\vert^{n-2}} + \psi_i(x)
\end{equation}
for all $x \in {\mathbb R}^n\backslash\hat{\mathcal S}$. Since $\tilde u_i \ge 0$, the harmonic function $\psi_i$
is bounded from below. By Liouville's theorem we then get that $\psi_i$ is constant. Clearly, since 
$\tilde u_i \ge 0$, we also have that the $A_j^i$'s are nonnegative. Letting $\varphi_i$ be given by
\begin{equation}\label{ProofStep6Sec7Eqt17}
\varphi_i(x) = \sum_{j=2}^m\frac{A_j^i}{\vert x-x_j\vert^{n-2}} + \psi_i\hskip.1cm ,
\end{equation}
the $\varphi_i$'s are harmonic in $B_0(\delta)$ for some 
$\delta > 0$ small, and nonnegative. Then, since 
the $\tilde u_\alpha^i$'s converge to $\tilde u_i$ 
in $C^1_{loc}({\mathbb R}^n\backslash\hat{\mathcal S})$, we get with 
(\ref{ProofStep6Sec7Eqt16}) and (\ref{ProofStep6Sec7Eqt17}) that
$$\tilde u_\alpha^i(x) \rightarrow \frac{A_i}{\vert x\vert^{n-2}} + \varphi_i(x)$$
in $C^1_{loc}(B_0(\delta)\backslash\{0\})$ for all $i$, where $A_i = A_1^i$ is nonnegative. 
In particular, (\ref{SharpPtseEstSec7}) holds true, and 
this proves Step \ref{step6sec7}.
\end{proof}

As a remark, note we proved a slightly more general result than the one in 
Step \ref{step6sec7}. Namely that there exists a finite subset $\hat{\mathcal S}$ of ${\mathbb R}^n$, 
given by (\ref{ProofStep6Sec7Eqt2}), such that, up to a subsequence, and for any $i$,
\begin{equation}\label{Sec7RefinStep6Eqt}
\tilde u_\alpha^i(x) \to \sum_{j=1}^m\frac{A_j^i}{\vert x-x_j\vert^{n-2}} + K_i
\end{equation}
in $C^1_{loc}({\mathbb R}^n\backslash\hat{\mathcal S})$ as $\alpha \to +\infty$, where the $A^i_j$'s and $K_i$'s are nonnegative 
constants, and where the $x_j$'s, $j = 1,\dots,m$, are the points in $\hat{\mathcal S}$. 
Now, going on with the proof of Theorem \ref{SharpPtwseAsymptThm}, we claim that if 
the $-A(\alpha)$'s are cooperative for all $\alpha$, then 
at least one of the $A_i$'s in (\ref{SharpPtseEstSec7}) is positive. This is the best we can prove in the sense that 
it might be that one and only one of the $A_i$'s is positive. Easy such examples (for instance, on the sphere) 
can be exhibited. Step \ref{step7sec7} states as follows.

\begin{Step}\label{step7sec7} Let $({\mathcal U}_\alpha)_\alpha$ 
be a bounded sequence in $H_{1,p}^2(M)$ of nonnegative solutions of (\ref{GenericEqtSec7}) 
which blows up. If the $-A(\alpha)$'s are cooperative for all $\alpha$, then 
at least one of the $A_i$'s in (\ref{SharpPtseEstSec7}) is positive.
\end{Step}

\begin{proof}[Proof of Step \ref{step7sec7}] We prove that if the $-A(\alpha)$'s are cooperative for all $\alpha$, then 
we have that $A_i > 0$, where $i$ is such that $x_\alpha = x^i_{j,\alpha}$ for some $j$, and 
the $x^i_{j,\alpha}$'s are the centers of the $1$-bubbles $(B_{j,\alpha}^i)_\alpha$ in (\ref{EqtThmSobTherSec4Bis}) 
from which the $p$-bubbles $({\mathcal B}_{j,\alpha})_\alpha$'s in 
(\ref{EqtThmSobTherSec4}) are defined. By (\ref{HighestWeightBisSec7}), 
$i = 1$. We let $L_g^\alpha$ be the operator acting on functions given by 
$L_g^\alpha u = \Delta_gu + A_{11}^\alpha u$.
By (\ref{CondSec7}), 
the operators $\Delta_g + A(\alpha)$ are coercive for $\alpha$ sufficiently large. 
Up to passing to a subsequence, we can 
assume they are coercive for all $\alpha$. Then, as is easily checked by 
considering $p$-maps like ${\mathcal U} = (u,0,\dots,0)$, $L_g^\alpha$ is also coercive for all $\alpha$. 
We let $G_\alpha$ be the Green's function of $L_g^\alpha$. By (\ref{CondSec7}), and standard properties 
of Green's functions (see, for instance, appendix A in Druet, Hebey and Robert 
\cite{DruHebRob2}), we can write that there exists $C > 0$ such that
\begin{equation}\label{ProofStep7Sec7Eqt1}
G_\alpha(x,y) \ge \frac{C}{d_g(x,y)^{n-2}}
\end{equation}
for all $\alpha$, and all $x \not= y$. Up to passing to another subsequence, we also assume in what follows that 
the ratios $d_g(x_\alpha,x^i_{j,\alpha})/\mu_\alpha$ have limits (possibly $+\infty$) 
as $\alpha \to +\infty$, for all $i, j$, where the $\mu_\alpha$'s are given by (\ref{HighestWeightSec7}) 
and (\ref{HighestWeightBisSec7}) . 
We let $\delta_1 < \delta_2$, $\delta_1, \delta_2 > 0$, be such that
$\left[\delta_1,\delta_2\right]$ does not contain any such limits. Then, there exists $C > 0$ 
such that
\begin{equation}\label{ProofStep7Sec7Eqt2}
d_g\left(x^i_{j,\alpha},\exp_{x_\alpha}(\mu_\alpha x)\right) \ge C \mu_\alpha
\end{equation}
for all $x \in B_0(\delta_2)\backslash B_0(\delta_1)$, all $i, j$, and all $\alpha$ (up to passing to a subsequence so 
that we only have to consider large $\alpha$'s). Let $v_\alpha^1$ be the function given by
\begin{equation}\label{ProofStep7Sec7Eqt3}
v_\alpha^1(x) = \mu_\alpha^{\frac{n-2}{2}} u_\alpha^1\left(\exp_{x_\alpha}(\mu_\alpha x)\right)\hskip.1cm .
\end{equation}
By Theorem \ref{PointEstThm}, and (\ref{ProofStep7Sec7Eqt2}), there exists $C > 0$ such that $v_\alpha^1 \le C$ 
in $B_0(\delta_2)\backslash B_0(\delta_1)$ for all $\alpha$. In particular, we can write that
\begin{equation}\label{ProofStep7Sec7Eqt4}
\int_{B_0(\delta_2)\backslash B_0(\delta_1)}(v_\alpha^1)^{2^\star}dv_{g_\alpha} 
\le C \int_{B_0(\delta_2)\backslash B_0(\delta_1)}(v_\alpha^1)^{2^\star-1}dv_{g_\alpha}
\end{equation}
where $g_\alpha$ is the metric given by $g_\alpha(x) = \left(\exp_{x_\alpha}^\star g\right)(\mu_\alpha x)$. 
Clearly, $g_\alpha \to \xi$ in $C^2(K)$  as $\alpha \to +\infty$ for all $K \subset\subset {\mathbb R}^n$. 
By Theorem \ref{SobTheorThm}, and more precisely (\ref{EqtThmSobTherSec4Bis}) that we consider for $i = 1$, 
we have that
\begin{equation}\label{ProofStep7Sec7Eqt5}
\begin{split} 
\int_{B_0(\delta_2)\backslash B_0(\delta_1)}(v_\alpha^1)^{2^\star}dv_{g_\alpha}
& = \int_{B_{x_\alpha}(\delta_2\mu_\alpha)\backslash B_{x_\alpha}(\delta_1\mu_\alpha)}
(u_\alpha^1)^{2^\star}dv_g\\
& \ge \int_{B_{x_\alpha}(\delta_2\mu_\alpha)\backslash B_{x_\alpha}(\delta_1\mu_\alpha)}B_\alpha^{2^\star}dv_g + o(1)
\end{split}
\end{equation}
for all $\alpha$, where $(B_\alpha)_\alpha$ is the $1$-bubble of center the $x_\alpha$'s and 
weights the $\mu_\alpha$'s, and where $o(1) \to 0$ as $\alpha \to +\infty$. We also have that 
\begin{equation}\label{ProofStep7Sec7Eqt6}
\int_{B_{x_\alpha}(\delta_2\mu_\alpha)\backslash B_{x_\alpha}(\delta_1\mu_\alpha)}B_\alpha^{2^\star}dv_g 
= \int_{B_0(\delta_2)\backslash B_0(\delta_1)}u^{2^\star}dv_{g_\alpha}\hskip.1cm ,
\end{equation}
where $u$ is given by (\ref{PosSolCritEuclEqt}). Combining (\ref{ProofStep7Sec7Eqt5}) and 
(\ref{ProofStep7Sec7Eqt6}), it follows that
\begin{equation}\label{ProofStep7Sec7Eqt7}
\int_{B_0(\delta_2)\backslash B_0(\delta_1)}(v_\alpha^1)^{2^\star}dv_{g_\alpha} \ge C
\end{equation}
for all $\alpha$, where $C > 0$ is independent of $\alpha$, and the $v_\alpha^1$'s are given by 
(\ref{ProofStep7Sec7Eqt3}). We let $\delta > 0$ small as in Step \ref{step6sec7}. For $x \in B_0(\delta)\backslash\{0\}$, 
$y \in B_0(\delta_2)\backslash B_0(\delta_1)$, and $\alpha$ sufficiently large such that $x \not= \sqrt{\mu_\alpha}y$, 
we let also $\hat G_\alpha(x,y)$ be given by
\begin{equation}\label{ProofStep7Sec7Eqt8}
\hat G_\alpha(x,y) = G_\alpha\left(\exp_{x_\alpha}(\sqrt{\mu_\alpha}x),\exp_{x_\alpha}(\mu_\alpha y)\right)
\hskip.1cm .
\end{equation}
By (\ref{GenericEqtSec7}),
$L_g^\alpha u_\alpha^1 = (u_\alpha^1)^{2^\star-1} - \sum_{j\not= 1}A^\alpha_{ij}u_\alpha^j$.
We fix $x$. Then,
\begin{equation}\label{EqtStep7ForStep8Sec7}
\begin{split}
\tilde u_\alpha^1(x)
& =  \int_MG_\alpha\left(\exp_{x_\alpha}(\sqrt{\mu_\alpha}x),y\right)(u_\alpha^1)^{2^\star-1}(y)dv_g(y)\\
&\hskip.4cm - \sum_{j\not= 1}\int_MG_\alpha\left(\exp_{x_\alpha}(\sqrt{\mu_\alpha}x),y\right)A^\alpha_{1j}(y)u_\alpha^j(y)dv_g(y)
\hskip.1cm ,
\end{split}
\end{equation}
and since we assumed that the $-A(\alpha)$'s are cooperative for all $\alpha$, we can write that
\begin{equation}\label{ProofStep7Sec7Eqt9}
\tilde u_\alpha^1(x) \ge 
\int_MG_\alpha\left(\exp_{x_\alpha}(\sqrt{\mu_\alpha}x),y\right)(u_\alpha^1)^{2^\star-1}(y)dv_g(y)
\end{equation}
for all $\alpha$ sufficiently large. By (\ref{ProofStep7Sec7Eqt9}) we also have that
\begin{equation}\label{ProofStep7Sec7Eqt10}
\begin{split}
\tilde u_\alpha^1(x)
&\ge \int_{B_{x_\alpha}(\delta_2\mu_\alpha)\backslash B_{x_\alpha}(\delta_1\mu_\alpha)}
G_\alpha\left(\exp_{x_\alpha}(\sqrt{\mu_\alpha}x),y\right)(u_\alpha^1)^{2^\star-1}(y)dv_g(y)\\
&\ge \mu_\alpha^{\frac{n-2}{2}} \int_{B_0(\delta_2)\backslash B_0(\delta_1)}
\hat G_\alpha(x,y)(v_\alpha^1)^{2^\star-1}(y)dv_{g_\alpha}(y)\hskip.1cm ,
\end{split}
\end{equation}
while, by (\ref{ProofStep7Sec7Eqt1}), 
\begin{equation}\label{ProofStep7Sec7Eqt11}
\begin{split}
\mu_\alpha^{\frac{n-2}{2}}\hat G_\alpha(x,y)
&\ge \frac{C\mu_\alpha^{\frac{n-2}{2}}}{\left\vert\sqrt{\mu_\alpha}x - \mu_\alpha y\right\vert^{n-2}}\\
&\ge \frac{C}{\left\vert x - \sqrt{\mu_\alpha}y\right\vert^{n-2}}\hskip.1cm .
\end{split}
\end{equation}
for all $y \in B_0(\delta_2)\backslash B_0(\delta_1)$, and all $\alpha$ sufficiently large. Combining 
(\ref{ProofStep7Sec7Eqt4}), (\ref{ProofStep7Sec7Eqt7}), (\ref{ProofStep7Sec7Eqt10}), and 
(\ref{ProofStep7Sec7Eqt11}), we get that there exists $C > 0$ such that
\begin{equation}\label{ProofStep7Sec7ConclEqt1}
\tilde u_1(x) \ge \frac{C}{\vert x\vert^{n-2}}
\end{equation}
for all $x \in B_0(\delta)\backslash\{0\}$, where $\tilde u_1$ is the $C^1_{loc}$-limit of the $\tilde u_\alpha^1$'s 
(which exists by Step \ref{step6sec7}). By (\ref{SharpPtseEstSec7}),
\begin{equation}\label{ProofStep7Sec7ConclEqt2}
\tilde u_1(x) = \frac{A_1}{\vert x\vert^{n-2}} + \varphi_1(x)\hskip.1cm ,
\end{equation}
for all $x \in B_0(\delta)\backslash\{0\}$, where $\varphi_1$ is harmonic in $B_0(\delta)$. In particular, it follows from 
(\ref{ProofStep7Sec7ConclEqt1}) and (\ref{ProofStep7Sec7ConclEqt2}) that $A_1 > 0$. This proves 
Step \ref{step7sec7}.
\end{proof}

The last step in the proof of Theorem \ref{SharpPtwseAsymptThm} 
is as follows.

\begin{Step}\label{step8sec7} Let $({\mathcal U}_\alpha)_\alpha$ 
be a bounded sequence in $H_{1,p}^2(M)$ of nonnegative solutions of (\ref{GenericEqtSec7}) 
which blows up. If the $-A(\alpha)$'s are cooperative for all $\alpha$, 
the limit system (\ref{LimitSystSec7}) is fully coupled, and ${\mathcal U}^0 \not\equiv 0$, where 
${\mathcal U}^0$ is as in Theorem \ref{SobTheorThm}, then the $\varphi_i$'s in (\ref{SharpPtseEstSec7}) are positive functions 
in $B_0(\delta)$ for all $i$.
\end{Step}

\begin{proof}[Proof of Step \ref{step8sec7}] By (\ref{CondSec7}), the limit map $-A$ is cooperative 
when the $-A(\alpha)$'s are assumed to be cooperative. Assuming that ${\mathcal U}^0 \not\equiv 0$, since 
the limit system (\ref{LimitSystSec7}) is also fully coupled and ${\mathcal U}^0$ is a solution 
of (\ref{LimitSystSec7}), we get, 
see Lemma \ref{PosSolSec1} in Section \ref{Prel}, that $u^0_i > 0$ for all $i$, where the $u^0_i$'s 
are the components of ${\mathcal U}^0$. We fix $i = 1,\dots,p$. Then we let $L_g^\alpha$ be the operator 
given by $L_g^\alpha u = \Delta_gu + A_{ii}^\alpha u$, and let $G_\alpha$ be the Green's function of $L_g^\alpha$. 
As already mentionned, see (\ref{Sec7RefinStep6Eqt}), 
there exists a finite subset $\hat{\mathcal S}$ of ${\mathbb R}^n$, 
such that, up to a subsequence, 
\begin{equation}\label{ProofStep8Sec7ConvEqt}
\tilde u_\alpha^i(x) \to \sum_{j=1}^m\frac{A_j^i}{\vert x-x_j\vert^{n-2}} + K_i
\end{equation}
in $C^1_{loc}({\mathbb R}^n\backslash\hat{\mathcal S})$ as $\alpha \to +\infty$, where the $A^i_j$'s and $K_i$'s are nonnegative 
constants, and the $x_j$'s, $j = 1,\dots,m$, are the points in $\hat{\mathcal S}$. 
We let $x \in {\mathbb R}^n\backslash\{0\}$. Like in 
(\ref{EqtStep7ForStep8Sec7}), we can write that
\begin{equation}\label{ProofStep8Sec7Eqt1}
\begin{split}
\tilde u_\alpha^i(x)
& =  \int_MG_\alpha\left(\exp_{x_\alpha}(\sqrt{\mu_\alpha}x),y\right)(u_\alpha^i)^{2^\star-1}(y)dv_g(y)\\
&\hskip.4cm - \sum_{j\not= i}\int_MG_\alpha\left(\exp_{x_\alpha}(\sqrt{\mu_\alpha}x),y\right)A^\alpha_{ij}(y)u_\alpha^j(y)dv_g(y)
\hskip.1cm ,
\end{split}
\end{equation}
for all $\alpha$, where the $x_\alpha$'s and $\mu_\alpha$'s are given by 
(\ref{HighestWeightSec7}) and (\ref{HighestWeightBisSec7}). Then, 
by (\ref{ProofStep8Sec7Eqt1}), 
\begin{eqnarray*}
\tilde u_\alpha^i(x)
& = & \int_MG_\alpha\left(\exp_{x_\alpha}(\sqrt{\mu_\alpha}x),y\right)(u_\alpha^i)^{2^\star-1}(y)dv_g(y)\\
&&\hskip.4cm - \sum_{j\not= i}
\int_{M\backslash B_{x_0}(r)}G_\alpha\left(\exp_{x_\alpha}(\sqrt{\mu_\alpha}x),y\right)A^\alpha_{ij}(y)u_\alpha^j(y)dv_g(y)\\
&&\hskip.4cm - \sum_{j\not= i}
\int_{B_{x_0}(r)}G_\alpha\left(\exp_{x_\alpha}(\sqrt{\mu_\alpha}x),y\right)A^\alpha_{ij}(y)u_\alpha^j(y)dv_g(y)
\end{eqnarray*}
for all $\alpha$ and  $r > 0$ small, where $x_0$ is the limit of the $x_\alpha$'s, 
and since the $-A(\alpha)$'s are cooperative, we can write that
\begin{equation}\label{ProofStep8Sec7Eqt2}
\begin{split}
\tilde u_\alpha^i(x)
& \ge \int_{M\backslash B_{x_0}(r)}G_\alpha\left(\exp_{x_\alpha}(\sqrt{\mu_\alpha}x),y\right)(u_\alpha^i)^{2^\star-1}(y)dv_g(y)\\
&\hskip.4cm - \sum_{j\not= i}
\int_{M\backslash B_{x_0}(r)}G_\alpha\left(\exp_{x_\alpha}(\sqrt{\mu_\alpha}x),y\right)A^\alpha_{ij}(y)u_\alpha^j(y)dv_g(y)
\end{split}
\end{equation}
for all $\alpha$ and  $r > 0$ small. Letting $\alpha \to +\infty$, it follows from 
(\ref{ProofStep8Sec7Eqt2}) that
\begin{eqnarray*}
\tilde u_i(x)
& \ge & \int_{M\backslash B_{x_0}(r)}G\left(x_0,y\right)(u_i^0)^{2^\star-1}(y)dv_g(y)\\
&&\hskip.4cm - \sum_{j\not= i}
\int_{M\backslash B_{x_0}(r)}G\left(x_0,y\right)A_{ij}(y)u_j^0(y)dv_g(y)
\end{eqnarray*}
where $\tilde u_i$ is the $C^1_{loc}$-limit of the $\tilde u_\alpha^i$'s in ${\mathbb R}^n\backslash\hat{\mathcal S}$, 
and $G$ is the Green's function of $\Delta_g + A_{ii}$. 
Then, letting $r \to 0$, we get that
\begin{eqnarray*}
\tilde u_i(x)
& \ge & \int_MG\left(x_0,y\right)(u_i^0)^{2^\star-1}(y)dv_g(y)\\
&&\hskip.4cm - \sum_{j\not= i}\int_MG\left(x_0,y\right)A_{ij}(y)u_j^0(y)dv_g(y)\\
& \ge & u_i^0(x_0)
\end{eqnarray*}
since ${\mathcal U}^0$ is a solution of the limit system (\ref{LimitSystSec7}). By (\ref{ProofStep8Sec7ConvEqt}), 
letting $\vert x\vert \to +\infty$, we then get that $K_i \ge u^0_i(x_0)$, and hence that $K_i > 0$ for all $i$. 
Since $\varphi_i$ in (\ref{SharpPtseEstSec7}) is given by
$$\varphi_i(x) = \sum_{x_j \in \hat{\mathcal S}\backslash\{0\}}\frac{A_j^i}{\vert x-x_j\vert^{n-2}} + K_i\hskip.1cm ,$$
and since the $A_j^i$'s are nonnegative, it follows that the $\varphi_i $'s are positive functions 
in $B_0(\delta)$ for all $i$. This proves Step \ref{step8sec7}.
\end{proof}

Clearly Theorem \ref{SharpPtwseAsymptThm} 
follows from Steps \ref{step1sec7} to \ref{step8sec7}. 
The proof of Theorem \ref{SharpPtwseAsymptThm} proceeds as follows.

\begin{proof}[Proof of Theorem \ref{SharpPtwseAsymptThm}] 
We let the $A(\alpha)$'s be such that (\ref{CondSec7}) is satisfied, 
and let $({\mathcal U}_\alpha)_\alpha$ 
be a bounded sequence in $H_{1,p}^2(M)$ of nonnegative solutions of (\ref{GenericEqtSec7}) 
which blows up. By Step \ref{step6sec7},
there exist $\delta > 0$, nonnegative real numbers $A_i$, and harmonic functions 
$\varphi_i: B_0(\delta) \to {\mathbb R}$, $i = 1,\dots,p$, such that, up to a subsequence, for any $i$,
$$\tilde u_\alpha^i(x) \rightarrow \frac{A_i}{\vert x\vert^{n-2}} + \varphi_i(x)$$
in $C^1_{loc}\left(B_0(\delta)\backslash\{0\}\right)$ as $\alpha \to +\infty$, where the $\tilde u_\alpha^i$'s 
are given by (\ref{DefRescaledFctSec7}). By Step \ref{step7sec7}, 
if we assume that the $-A(\alpha)$'s are cooperative for all $\alpha$, then 
at least one of the $A_i$'s is positive. By Step \ref{step8sec7}, if we also assume 
that the limit system (\ref{LimitSystSec7}) is fully coupled, and that ${\mathcal U}^0 \not\equiv 0$, where 
${\mathcal U}^0$ is as in Theorem \ref{SobTheorThm}, then the $\varphi_i$'s are positive functions 
in $B_0(\delta)$ for all $i$. This ends the proof of Theorem \ref{SharpPtwseAsymptThm}.
\end{proof}

Let ${\mathcal R} = \sum_{i=1}^pA_i\varphi_i(0)$. By Theorem 
\ref{SharpPtwseAsymptThm}, ${\mathcal R} > 0$ if we assume 
that the $-A(\alpha)$'s are cooperative for all $\alpha$, that ${\mathcal U}^0 \not\equiv 0$, and that 
the limit system (\ref{LimitSystSec7}) is fully coupled. 
If we drop this assumption that the limit system (\ref{LimitSystSec7}) should be 
fully coupled, and assume only that the $-A(\alpha)$'s are cooperative for all $\alpha$ 
and that ${\mathcal U}^0 \not\equiv 0$, we might get that 
${\mathcal R} = 0$. For instance, if $u_0$ is a positive solution of the Yamabe equation 
(\ref{YamEqtSphereEx}) on the unit sphere, and $(u_\alpha)_\alpha$ 
is a blowing-up sequence of positive solutions of (\ref{YamEqtSphereEx}), then 
${\mathcal U}_\alpha = (u_0,u_\alpha)$ is a solution of (\ref{GenericEqtSec7}) on the unit 
sphere, where $p = 2$ and $A(\alpha)$ is such that $A^\alpha_{11} = \frac{n(n-2)}{4}$, 
$A_{12}^\alpha = 0$, and $A_{22}^\alpha =\frac{n(n-2)}{4}$ for all $\alpha$. 
Here ${\mathcal U}^0 = (u_0,0)$ is not zero. However, ${\mathcal R} = 0$. 

\section{Standard rescaling}\label{StdRescaling}

In what follows we let $(M,g)$ be a smooth compact Riemannian manifold of dimension $n \ge 3$,  $p \ge 1$ integer, 
and $\left(A(\alpha)\right)_\alpha$, $\alpha\in{\mathbb N}$, be a sequence of smooth maps 
$A(\alpha): M \to M_p^s({\mathbb R})$. We let also $A(\alpha) = (A_{ij}^\alpha)$, and consider systems like
\begin{equation}\label{GenericEqtSec8}
\Delta_gu_i + \sum_{j=1}^pA_{ij}^\alpha(x)u_j = \vert u_i\vert^{2^\star-2}u_i
\end{equation}
in $M$, for all $i = 1,\dots p$. We assume that the $A(\alpha)$'s satisfy that 
there exists a $C^{0,\theta}$-map $A: M \to M_p^s({\mathbb R})$, 
$A = (A_{ij})$ and $0 < \theta < 1$, such that
\begin{equation}\label{CondSec8}
A^\alpha_{ij} \to A_{ij}\hskip.1cm\hbox{in}\hskip.1cm C^{0,\theta}(M)
\end{equation}
for all $i, j$ as $\alpha \to +\infty$. 
We let $({\mathcal U}_\alpha)_\alpha$, ${\mathcal U}_\alpha \not\equiv 0$, 
be a bounded sequence in $H_{1,p}^2(M)$ 
of nonnegative solutions of (\ref{GenericEqtSec8}). In particular, the sequence $({\mathcal U}_\alpha)_\alpha$ is a 
Palais-Smale sequence for (\ref{GenericEqtSec8}) and we can apply 
Theorem \ref{SobTheorThm}. Up to a subsequence, (\ref{EqtThmSobTherSec4}) 
and (\ref{EqtThmSobTherSec4Bis}) hold. In what follows, we assume that the sequence blows up and 
let the $x_{j,\alpha}^i$'s and $\mu_{j,\alpha}^i$'s be the centers and weights 
of the $1$-bubbles $(B_{j,\alpha}^i)_\alpha$ in 
(\ref{EqtThmSobTherSec4Bis}) from which the $p$-bubbles $({\mathcal B}_{j,\alpha})_\alpha$'s in 
(\ref{EqtThmSobTherSec4}) are defined, where $i = 1,\dots,p$, $j = 1,\dots,k_i$, 
and $\sum_i k_i = k$.  We fix $i_0 = 1,\dots,p$ and $j_0 = 1,\dots,k_i$ arbitrary. Then, 
for $x_0 \in {\mathbb R}^n$, we define $\hat x_\alpha \in {\mathbb R}^n$ and 
$\hat\mu_\alpha > 0$ by
\begin{equation}\label{PertCenterEqtSec8FirstEqt}
\hat x_\alpha = \exp_{x^{i_0}_{j_0,\alpha}}\left(\hat\mu_\alpha x_0\right)
\hskip.1cm\hbox{and}\hskip.1cm \hat\mu_\alpha = \mu^{i_0}_{j_0,\alpha}\hskip.1cm ,
\end{equation}
for all $\alpha$, 
where the notation $\exp_x$ stands for the exponential map at $x$. 
We also define the 
standard rescaled functions $\hat u_\alpha^i$ 
of the $u_\alpha^i$'s with respect to the $\hat x_\alpha$'s and $\hat\mu_\alpha$'s, 
$i = 1,\dots,p$, by
\begin{equation}\label{StdRescalFctEqtSec8FirstEqt}
\hat u_\alpha^i(x) = \hat\mu_\alpha^{\frac{n-2}{2}} 
u_\alpha^i\left(\exp_{\hat x_\alpha}(\hat\mu_\alpha x)\right)\hskip.1cm ,
\end{equation}
where $x \in B_0(1)$, the $\hat x_\alpha$'s 
and $\hat\mu_\alpha$'s are as in (\ref{PertCenterEqtSec8FirstEqt}), and 
${\mathcal U}_\alpha = (u_\alpha^1,\dots,u_\alpha^p)$ for all $\alpha$.
We claim here that the following convergence result holds.

\begin{thm}[Standard Rescaling]\label{StdRescThm} Let $(M,g)$ be a smooth compact Riemannian 
manifold of dimension $n \ge 3$, $p \ge 1$, and $\left(A(\alpha)\right)_\alpha$ a sequence of smooth maps 
$A(\alpha): M \to M_p^s({\mathbb R})$ satisfying (\ref{CondSec8}).
Let $({\mathcal U}_\alpha)_\alpha$ 
be a bounded sequence in $H_{1,p}^2(M)$ of nonnegative solutions of (\ref{GenericEqtSec8}) 
which blows up. 
Then there exist $\delta > 0$, and $x_0 \in {\mathbb R}^n$ such that, up to a subsequence, for any $i$,
\begin{equation}\label{EqtTheStdRescThmSec8}
\hat u_\alpha^i \rightarrow u_i\hskip.2cm\hbox{in}\hskip.1cm C^2\left(B_0(\delta)\right)
\end{equation}
as $\alpha \to +\infty$, where $\hat u_\alpha^i$ is given by 
(\ref{PertCenterEqtSec8FirstEqt})--(\ref{StdRescalFctEqtSec8FirstEqt}),  
and the $u_i$'s are nonnegative (not all trivial) solutions of the Euclidean equation 
$\Delta u = u^{2^\star-1}$ in ${\mathbb R}^n$. 
\end{thm}

The $u_i$'s, $i = 1,\dots,p$, see the discussion at the beginning of Section \ref{SobTheory}, can be written as
\begin{equation}\label{LimitStdRescSec8}
u_i(x) = \left(\frac{\lambda_i}{\lambda_i^2 + \frac{\vert x-x_i\vert^2}{n(n-2)}}\right)^{\frac{n-2}{2}}
\end{equation}
where $x \in {\mathbb R}^n$, $\lambda_i \ge 0$, and $x_i \in {\mathbb R}^n$. 
An equivalent statement to the statement that the $u_i$'s are not all trivial is that 
at least one of the $\lambda_i$'s is positive. 

\begin{proof} First we claim that we can choose $x_0 \in {\mathbb R}^n$, $\vert x_0\vert > 0$ small, 
such that, up to a subsequence,
\begin{equation}\label{ProofThSec8Eqt1}
d_g\left(\hat x_\alpha,x^i_{j,\alpha}\right) \ge R\hat\mu_\alpha
\end{equation}
for all $\alpha$ and all $i, j$, where $R > 0$ small is independent of $\alpha$, $i$, and $j$. 
In order to prove (\ref{ProofThSec8Eqt1}), we 
let $\tilde x_\alpha = x^i_{j,\alpha}$ when $i = i_0$ and $j = j_0$. Up to passing to a subsequence, 
we can assume that the ratios $d_g(\tilde x_\alpha,x^i_{j,\alpha})/\hat\mu_\alpha$ 
have limits (possibly $+\infty$) as $\alpha \to +\infty$, for all $i, j$. We write that
$$\lambda^i_j = \lim_{\alpha\to+\infty} \frac{d_g(\tilde x_\alpha,x^i_{j,\alpha})}{\hat\mu_\alpha}\hskip.1cm ,$$
and let ${\mathcal H}$ be the set consisting of the $i, j$'s such that 
$\lambda^i_j > 0$. Let $\lambda > 0$ be such that 
$\lambda < \lambda^i_j$ for all $i, j \in {\mathcal H}$. We choose $x_0 \in {\mathbb R}^n$ such that 
$\vert x_0\vert < \lambda/2$. If $i, j \in {\mathcal H}$, then
\begin{eqnarray*}
\frac{d_g\left(\hat x_\alpha,x^i_{j,\alpha}\right)}{\hat\mu_\alpha}
& \ge & \frac{d_g\left(\tilde x_\alpha,x^i_{j,\alpha}\right)}{\hat\mu_\alpha} 
- \frac{d_g\left(\hat x_\alpha,\tilde x_\alpha\right)}{\hat\mu_\alpha}\\
& = &  \frac{d_g\left(\tilde x_\alpha,x^i_{j,\alpha}\right)}{\hat\mu_\alpha} 
- \vert x_0\vert\\
& \ge & \vert x_0\vert
\end{eqnarray*}
for $\alpha$ large, and we get that (\ref{ProofThSec8Eqt1}) holds true for any such $i, j$'s. If, on the contrary, 
$i, j \not\in {\mathcal H}$, then $\lambda^i_j = 0$, and
\begin{eqnarray*}
\frac{d_g\left(\hat x_\alpha,x^i_{j,\alpha}\right)}{\hat\mu_\alpha}
& \ge & \frac{d_g\left(\hat x_\alpha,\tilde x_\alpha\right)}{\hat\mu_\alpha}
- \frac{d_g\left(\tilde x_\alpha,x^i_{j,\alpha}\right)}{\hat\mu_\alpha}\\
& = & \vert x_0\vert - \frac{d_g\left(\tilde x_\alpha,x^i_{j,\alpha}\right)}{\hat\mu_\alpha}\\
& \ge & \frac{\vert x_0\vert}{2}
\end{eqnarray*}
for $\alpha$ large and, here again, we get that (\ref{ProofThSec8Eqt1}) holds true. Clearly, this proves that 
we can choose $x_0 \in {\mathbb R}^n$, $\vert x_0\vert > 0$ small, 
such that, up to a subsequence, (\ref{ProofThSec8Eqt1}) is true. Let $\delta_0 > 0$ small 
such that $2\delta_0 < R$, 
where $R > 0$ is as in (\ref{ProofThSec8Eqt1}). By (\ref{ProofThSec8Eqt1}), we can write that for any $x \in B_0(\delta_0)$, 
any $i, j$, and any $\alpha$,
\begin{equation}\label{ProofThSec8Eqt2}
\begin{split}
d_g\left(x^i_{j,\alpha},\exp_{\hat x_\alpha}(\hat\mu_\alpha x)\right)
& \ge d_g\left(x^i_{j,\alpha},\hat x_\alpha\right) 
- d_g\left(\hat x_\alpha,\exp_{\hat x_\alpha}(\hat\mu_\alpha x)\right)\\
& \ge R\hat\mu_\alpha - \delta_0\hat\mu_\alpha\\
& \ge \frac{R}{2} \hat\mu_\alpha\hskip.1cm .
\end{split}
\end{equation}
By Theorem \ref{PointEstThm}, it follows from (\ref{ProofThSec8Eqt2}) that the $\hat u_\alpha^i$'s are 
uniformly bounded in $B_0(\delta_0)$. They are solutions of the system
\begin{equation}\label{ProofThSec8Eqt3}
\Delta_{g_\alpha}\hat u_\alpha^i + \hat\mu_\alpha^2\sum_{j=1}^p\hat A^\alpha_{ij}\hat u_\alpha^j = 
(\hat u_\alpha^i)^{2^\star-1}
\end{equation}
for all $i$, where 
$$g_\alpha(x) = \left(\exp_{\hat x_\alpha}^\star g\right)(\hat\mu_\alpha x)
\hskip.2cm\hbox{and}\hskip.2cm
\hat A^\alpha_{ij}(x) = A^\alpha_{ij}\left(\exp_{\hat x_\alpha}(\hat\mu_\alpha x)\right)
\hskip.1cm .$$
Clearly, 
$g_\alpha \to \xi$ in $C^2(K)$ for any $K \subset\subset {\mathbb R}^n$, where $\xi$ is the Euclidean 
metric. Then, since the $\hat u_\alpha^i$'s are 
uniformly bounded in $B_0(\delta_0)$, it follows from (\ref{CondSec8}), 
(\ref{ProofThSec8Eqt3}), and standard elliptic theory, that the $\hat u_\alpha^i$'s are 
uniformly bounded in $C^{2,\theta}\left(B_0(\delta_1)\right)$ for $\delta_1 < \delta_0$, 
where $0 < \theta < 1$. Up to a subsequence, we may then assume that
$\hat u_\alpha^i \rightarrow u_i$ in $C^2\left(B_0(\delta)\right)$ 
as $\alpha \to +\infty$, for some $\delta > 0$ small, and all $i$. By rescaling invariance, the 
$\hat u_\alpha^i$'s are bounded in $H_1^2\left(B_0(R^\prime)\right)$ for all $R^\prime > 0$. 
Up to passing to another subsequence, we may therefore 
assume that we also have that $\hat u_\alpha^i \rightharpoonup v_i$ weakly in $H_{1,loc}^2({\mathbb R}^n)$, 
that  $\hat u_\alpha^i \to v_i$ strongly in $L^2_{loc}({\mathbb R}^n)$, and that $\hat u_\alpha^i \to v_i$ 
almost everywhere, for all $i$. 
We can also write that 
$$(\hat u_\alpha^i)^{2^\star-1} \rightharpoonup (v_i)^{2^\star-1}$$
weakly in 
$L^s_{loc}({\mathbb R}^n)$ for all $i$, where $s = 2^\star/(2^\star-1)$. 
Then, we get with (\ref{CondSec8}) and (\ref{ProofThSec8Eqt3}) that the $v_i$'s are nonnegative solutions 
of the Euclidean equation $\Delta u = u^{2^\star-1}$ in ${\mathbb R}^n$.
Clearly, we also have that $v_i = u_i$ in $B_0(\delta)$, for all $i$. It follows that 
there exist $\delta > 0$, and $x_0 \in {\mathbb R}^n$ such that, up to a subsequence, and for any $i$,
\begin{equation}\label{ProofThSec8Eqt4}
\hat u_\alpha^i \rightarrow u_i\hskip.2cm\hbox{in}\hskip.1cm C^2\left(B_0(\delta)\right)
\end{equation}
as $\alpha \to +\infty$, where 
the $u_i$'s are nonnegative solutions of the Euclidean equation 
$\Delta u = u^{2^\star-1}$ in ${\mathbb R}^n$. Now it remains to prove that the $u_i$'s are not all trivial. 
We prove in what follows that $u_{i_0} \not\equiv 0$, and hence that $u_{i_0}$ is everywhere positive 
in ${\mathbb R}^n$. By Theorem \ref{SobTheorThm}, given $0 < r < \delta$, we can write that
\begin{equation}\label{ProofThSec8Eqt5}
\int_{B_{\hat x_\alpha}(r\hat\mu_\alpha)}(u_\alpha^{i_0})^{2^\star}dv_g 
\ge \int_{B_{\hat x_\alpha}(r\hat\mu_\alpha)}(B_{j_0,\alpha}^{i_0})^{2^\star}dv_g + o(1)
\hskip.1cm ,
\end{equation}
where $o(1) \to 0$ as $\alpha \to +\infty$, and $(B_{j_0,\alpha}^{i_0})_\alpha$ is the $1$-bubble of centers the 
$\tilde x_\alpha$'s and weights the $\hat\mu_\alpha$'s. 
By rescaling invariance, and (\ref{ProofThSec8Eqt4}), we also have that
\begin{equation}\label{ProofThSec8Eqt6}
\begin{split}
\int_{B_{\hat x_\alpha}(r\hat\mu_\alpha)}(u_\alpha^{i_0})^{2^\star}dv_g
& = \int_{B_0(r)}(\hat u_\alpha^{i_0})^{2^\star}dv_{g_\alpha}\\
& = \int_{B_0(r)}u_{i_0}^{2^\star}dx + o(1)\hskip.1cm ,
\end{split}
\end{equation}
where $o(1) \to 0$ as $\alpha \to +\infty$, since $g_\alpha \to \xi$ as $\alpha \to +\infty$. 
Independently,
\begin{equation}\label{ProofThSec8Eqt7}
\int_{B_{\hat x_\alpha}(r\hat\mu_\alpha)}(B_{j_0,\alpha}^{i_0})^{2^\star}dv_g 
= \int_{\Omega_\alpha}u^{2^\star}dv_{\hat g_\alpha}\hskip.1cm ,
\end{equation}
where $u$ is given by (\ref{PosSolCritEuclEqt}), $\hat g_\alpha$ is the metric given by 
$\hat g_\alpha(x) = \left(\exp_{\tilde x_\alpha}^\star g\right)(\hat\mu_\alpha x)$, and 
$$\Omega_\alpha = \frac{1}{\mu_\alpha} \exp_{\tilde x_\alpha}^{-1}\left(B_{\hat x_\alpha}(r\hat\mu_\alpha)\right)
\hskip.1cm .$$
As is easily checked, there exists $t_0 > 0$ such that $B_{x_0}(t_0) \subset \Omega_\alpha$ for all $\alpha$. Hence, 
\begin{equation}\label{ProofThSec8Eqt8}
\begin{split}
\int_{\Omega_\alpha}u^{2^\star}dv_{\hat g_\alpha}
&\ge \int_{B_{x_0}(t_0)}u^{2^\star}dv_{\hat g_\alpha}\\
&\ge \int_{B_{x_0}(t_0)}u^{2^\star}dx + o(1)\hskip.1cm ,
\end{split}
\end{equation}
where $o(1) \to 0$ as $\alpha \to +\infty$, since we also have that $\hat g_\alpha \to \xi$ as $\alpha \to +\infty$.
Combining (\ref{ProofThSec8Eqt5})--(\ref{ProofThSec8Eqt8}), letting $\alpha \to +\infty$, we get that
$$\int_{B_0(r)}u_{i_0}^{2^\star}dx \ge \int_{B_{x_0}(t_0)}u^{2^\star}dx > 0\hskip.1cm ,$$
and it follows that $u_{i_0} \not\equiv 0$. In particular, 
the $u_i$'s are not all trivial. Together with (\ref{ProofThSec8Eqt4}), this ends the 
proof of Theorem \ref{StdRescThm}.
\end{proof}

An easy but important corollary to Theorem \ref{StdRescThm} is as follows. 

\begin{cor}\label{CorSec8} Let $(M,g)$ be a smooth compact Riemannian 
manifold of dimension $n \ge 3$, $p \ge 1$, and $\left(A(\alpha)\right)_\alpha$ a sequence of smooth maps 
$A(\alpha): M \to M_p^s({\mathbb R})$ satisfying (\ref{CondSec8}).
Let $({\mathcal U}_\alpha)_\alpha$ 
be a bounded sequence in $H_{1,p}^2(M)$ of nonnegative solutions of (\ref{GenericEqtSec8}) 
which blows up. For any 
$\delta > 0$ there 
exists $C > 0$ such that, up to a subsequence,
\begin{equation}\label{L2boundSec8}
\sum_{i=1}^p\int_{B_{x_\alpha}(\delta\sqrt{\mu_\alpha})}(u_\alpha^i)^2dv_g \ge C\mu_\alpha^2
\hskip.1cm ,
\end{equation}
where ${\mathcal U}_\alpha = (u_\alpha^1,\dots,u_\alpha^p)$ for all $\alpha$, and the $x_\alpha$'s and $\mu_\alpha$'s 
are given by (\ref{HighestWeightSec7}) and (\ref{HighestWeightBisSec7}) in Section  \ref{SharpPtAsympt}. 
\end{cor}

\begin{proof} Let $x_\alpha$ and $\mu_\alpha$ be given by (\ref{HighestWeightSec7}) 
and (\ref{HighestWeightBisSec7}). By Theorem \ref{StdRescThm}, there exists $x_0 \in {\mathbb R}^n$ and 
$\delta_0 > 0$ such that, up to a subsequence, for any $i$,
$\hat u_\alpha^i \rightarrow u_i$ in $C^2\left(B_0(\delta_0)\right)$ 
as $\alpha \to +\infty$, where $\hat u_\alpha^i$ is given by
$$\hat u_\alpha^i(x) = \mu_\alpha^{\frac{n-2}{2}} 
u_\alpha^i\left(\exp_{\hat x_\alpha}(\mu_\alpha x)\right)\hskip.1cm ,$$
where $\hat x_\alpha = \exp_{x_\alpha}(\mu_\alpha x_0)$, 
and the $u_i$'s are nonnegative (not all trivial) solutions of the Euclidean equation 
$\Delta u = u^{2^\star-1}$ in ${\mathbb R}^n$. Let $\delta > 0$ be given, and $\delta_1 = \delta_0/2$. Clearly, 
$$B_{\hat x_\alpha}(\delta_1\mu_\alpha) \subset B_{x_\alpha}(\delta\sqrt{\mu_\alpha})$$
for $\alpha$ large. In particular, we can write that, up to a subsequence,
\begin{equation}\label{ProofCorSec8Eqt1}
\begin{split}
\sum_{i=1}^p\int_{B_{x_\alpha}(\delta\sqrt{\mu_\alpha})}(u_\alpha^i)^2dv_g
& \ge \sum_{i=1}^p\int_{B_{\hat x_\alpha}(\delta_1\mu_\alpha)}(u_\alpha^i)^2dv_g\\
& =  \mu_\alpha^2\sum_{i=1}^p\int_{B_0(\delta_1)}(\hat u_\alpha^i)^2dv_{g_\alpha}\hskip.1cm ,
\end{split}
\end{equation}
where $g_\alpha$ is the metric given by 
$$g_\alpha(x) = \left(\exp_{x_\alpha}^\star g\right)(\mu_\alpha x)\hskip.1cm .$$
Clearly, $g_\alpha \to \xi$ in $C^2(K)$ for all $K \subset\subset {\mathbb R}^n$ 
as $\alpha \to +\infty$. We also have that 
$\hat u_\alpha^i \rightarrow u_i$ in $C^0\left(B_0(\delta_0)\right)$ as $\alpha \to +\infty$, 
for all $i$. Since the $u_i$'s are not all trivial, so that at least one of the $u_i$'s is a positive 
function, we can write that, up to a subsequence,
\begin{equation}\label{ProofCorSec8Eqt2}
\sum_{i=1}^p\int_{B_0(\delta_1)}(\hat u_\alpha^i)^2dv_{g_\alpha}
\ge \frac{1}{2} \sum_{i=1}^p\int_{B_0(\delta_1)}u_i^2dx > 0
\end{equation}
for all $\alpha$. Combining (\ref{ProofCorSec8Eqt1}) 
and (\ref{ProofCorSec8Eqt2}), we get that (\ref{L2boundSec8}) holds true. This proves Corollary 
\ref{CorSec8}.
\end{proof}

By changing the $x_\alpha$'s in 
the estimates (\ref{SharpPtseEstSec7ThmStat}) of Theorem \ref{SharpPtwseAsymptThm}, as in (\ref{PertCenterEqtSec8FirstEqt}), 
we can ask to have both the estimates (\ref{EqtTheStdRescThmSec8})  
of Theorem \ref{StdRescThm} and the estimates (\ref{SharpPtseEstSec7ThmStat}) of Theorem \ref{SharpPtwseAsymptThm}. 
In order to see this, 
we assume in the sequel that (\ref{CondSec7}) 
holds. We let also $\mu_\alpha$ stand for 
the largest weight among all the possible weights $\mu_{j,\alpha}^i$, and 
$x_\alpha$ stand for the corresponding 
$x_{j,\alpha}^i$. 
Up to renumbering, and up to a subsequence, we can assume that the $x_\alpha$'s and $\mu_\alpha$'s 
correspond to the choice $i = j = 1$. Then 
(\ref{HighestWeightSec7}) and (\ref{HighestWeightBisSec7}) hold.
For $x_0 \in {\mathbb R}^n$, we let $\hat x_\alpha \in {\mathbb R}^n$ be given by 
the first equation in (\ref{PertCenterEqtSec8FirstEqt})  when $i_0 = j_0 = 1$. 
Namely,
\begin{equation}\label{PertCenterEqtSec8SecEqt}
\hat x_\alpha = \exp_{x_\alpha}\left(\mu_\alpha x_0\right)
\hskip.1cm ,
\end{equation}
and, as in (\ref{StdRescalFctEqtSec8FirstEqt}), we define the 
standard rescaled function $\hat u_\alpha^i$ 
of $u_\alpha^i$ with respect to $\hat x_\alpha$ and $\hat\mu_\alpha = \mu_\alpha$ by
\begin{equation}\label{StdRescalFctEqtSec8}
\hat u_\alpha^i(x) = \mu_\alpha^{\frac{n-2}{2}} 
u_\alpha^i\left(\exp_{\hat x_\alpha}(\mu_\alpha x)\right)\hskip.1cm ,
\end{equation}
where $x \in B_0(1)$, $i = 1,\dots,p$, and 
${\mathcal U}_\alpha = (u_\alpha^1,\dots,u_\alpha^p)$.
Given $R > 0$, we let  
$B_0(R)$ be the Euclidean ball centered at $0$ of radius $R$, 
and let $f_\alpha:B_0(R) \to {\mathbb R}^n$ 
be given by
\begin{equation}\label{ChgeExpSec8}
f_\alpha(x) = \frac{1}{\sqrt{\mu_\alpha}} \exp_{x_\alpha}^{-1}
\left(\exp_{\hat x_\alpha}\left(\sqrt{\mu_\alpha} x\right)\right)
\hskip.1cm ,
\end{equation}
where $x \in B_0(R)$, and $\hat x_\alpha$ is as in 
(\ref{PertCenterEqtSec8SecEqt}). It is easily seen that there exists $C > 0$ 
such that
$$\left\vert f_\alpha(x) - x\right\vert \le \frac{C}{\sqrt{\mu_\alpha}} d_g(x_\alpha,\hat x_\alpha)$$
for all $x \in B_0(R)$ and all $\alpha$. In particular, 
\begin{equation}\label{ConvPropPertCenterBisSec8}
\lim_{\alpha\to +\infty}\left(\sup_{x \in B_0(R)}\left\vert f_\alpha(x) - x\right\vert\right) = 0
\hskip.1cm .
\end{equation}
Choosing $R \gg 1$, we then get with Theorem 
\ref{SharpPtwseAsymptThm} and (\ref{ConvPropPertCenterBisSec8}) that 
there exist $\delta > 0$, nonnegative real numbers $A_i$, and harmonic functions 
$\varphi_i: B_0(\delta) \to {\mathbb R}$, $i = 1,\dots,p$, such that up to a subsequence, for any $i$,
\begin{equation}\label{SharpPtseEstSec8}
\tilde u_\alpha^i(x) \rightarrow \frac{A_i}{\vert x\vert^{n-2}} + \varphi_i(x)
\end{equation}
in $C^0_{loc}\left(B_0(\delta)\backslash\{0\}\right)$ as $\alpha \to +\infty$, where 
${\mathcal U}_\alpha = (u_\alpha^1,\dots,u_\alpha^p)$, the $\tilde u_\alpha^i$'s 
are the rescaled functions given by
$$\tilde u_\alpha^i(x) = u_\alpha^i\left(\exp_{\hat x_\alpha}\left(\sqrt{\mu_\alpha}x\right)\right)\hskip.1cm ,$$
$\hat x_\alpha$ is as in 
(\ref{PertCenterEqtSec8SecEqt}), and $x_0 \in {\mathbb R}^n$ is arbitrary. 
Independently, by Theorem \ref{StdRescThm}, 
we also have that for a suitable choice of $x_0$, and for $\delta > 0$ small, 
$\hat u_\alpha^i \to u_i$ in $C^2\left(B_0(\delta)\right)$ as $\alpha \to +\infty$, where the $\hat u_\alpha^i$'s are given by 
(\ref{StdRescalFctEqtSec8}),  
and the $u_i$'s are nonnegative (non all trivial) solutions of the Euclidean equation 
$\Delta u = u^{2^\star-1}$ in ${\mathbb R}^n$. 
In particular, 
by slightly changing the $x_\alpha$'s in the estimates of Theorem \ref{SharpPtwseAsymptThm},  
as in (\ref{PertCenterEqtSec8SecEqt}), 
we can ask to have both the estimates (\ref{EqtTheStdRescThmSec8})  
of Theorem \ref{StdRescThm} and the estimates (\ref{SharpPtseEstSec8}) of Theorem \ref{SharpPtwseAsymptThm}. 

\section{Pseudo-compactness and compactness}\label{Cptness}

In what follows we let $(M,g)$ be a smooth compact Riemannian manifold of dimension $n \ge 3$,  $p \ge 1$ integer, 
and $\left(A(\alpha)\right)_\alpha$, $\alpha\in{\mathbb N}$, be a sequence of smooth maps 
$A(\alpha): M \to M_p^s({\mathbb R})$. We let also $A(\alpha) = (A_{ij}^\alpha)$, and consider systems like
\begin{equation}\label{GenericEqtSec9}
\Delta_gu_i + \sum_{j=1}^pA_{ij}^\alpha(x)u_j = \vert u_i\vert^{2^\star-2}u_i
\end{equation}
in $M$, for all $i = 1,\dots p$. We assume that the $A(\alpha)$'s satisfy that 
there exists a $C^1$-map $A: M \to M_p^s({\mathbb R})$, 
$A = (A_{ij})$, such that
\begin{equation}\label{CondSec9}
\begin{split}
&\Delta_g^p+A\hskip.1cm\hbox{is coercive, and}\\
&A^\alpha_{ij} \to A_{ij}\hskip.1cm\hbox{in}\hskip.1cm C^1(M)
\end{split}
\end{equation}
as $\alpha \to +\infty$, where the second equation in (\ref{CondSec9}) 
should hold for all $i, j$. The limit system we get by combining (\ref{GenericEqtSec9}) and (\ref{CondSec9}) reads as
\begin{equation}\label{LimitSystSec9}
\Delta_gu_i + \sum_{j=1}^pA_{ij}(x)u_j = \vert u_i\vert^{2^\star-2}u_i
\end{equation}
in $M$, for all $i = 1,\dots p$. The goal in this section is to prove 
compactness results for sequences of nonnegative nontrivial solutions of (\ref{GenericEqtSec9}). 
Compactness for second order scalar equations goes back to the remarkable work of Schoen 
\cite{Sch2,Sch3,Sch4} on the 
Yamabe equation
\begin{equation}\label{YamEqtSec9}
\Delta_gu + \frac{n-2}{4(n-1)}S_gu = u^{2^\star-1}\hskip.1cm ,
\end{equation}
where $S_g$ is the scalar curvature of $g$. Compactness 
for scalar curvature type equations was also discussed in 
Devillanova and Solimini \cite{DevSol}, Druet \cite{Dru1,Dru2}, Druet and Hebey \cite{DruHeb1},  
Khuri and Schoen \cite{KhuSch}, Li and Zhu \cite{LiZhu}, Li and Zhang \cite{LiZha,LiZhabis}, 
and Marques \cite{Mar}. A possible survey reference 
is Druet and Hebey \cite{DruHeb2}. We refer also to Druet, Hebey and Robert \cite{DruHebRob1,DruHebRob2} where the blow-up analysis used in \cite{Dru1} is developed. 
The dynamical viewpoint, in the sense of the terminology introduced in Hebey 
\cite{HebBBConf}, where we discuss sequences 
of equations and not only one equation, is considered in 
Druet \cite{Dru1}, Druet and Hebey \cite{DruHeb1,DruHeb2}, 
and Druet, Hebey and Robert \cite{DruHebRob1,DruHebRob2}. With this viewpoint, that we also adopt here,  
several blow-up phenomena can be 
made concrete in specific examples. In particular, we can construct examples 
of sequences of solutions with no a priori 
bound on the energy, or which blow up with an 
arbitrary number of bubbles in their Sobolev decomposition 
and an arbitrary number of geometrical blow-up points. In the case of systems, 
by combining the examples in Druet and Hebey \cite{DruHeb1} with the construction at the beginning 
of Section \ref{HighEn}, we get that the following result holds.

\begin{prop}\label{PropExBlowUpSec9} Let $(M,g)$ be a space form of positive curvature $+1$, and 
dimension $n \ge 6$. Let also $s = \pm 1$. 
There exist sequences $(h_\alpha)_\alpha$, $(\tilde h_\alpha)_\alpha$ of smooth functions, 
converging in $C^1(M)$ to $\frac{n-2}{4(n-1)}S_g$, and sequences 
$(\varepsilon_\alpha)_\alpha$, $(\tilde\varepsilon_\alpha)_\alpha$, and $(\beta_\alpha)_\alpha$ 
of smooth positive functions, 
converging to $0$ in $C^1(M)$, such that 
the $2$-systems (\ref{GenericEqtSec9}), where $p = 2$ and
\begin{equation}\label{ExBlowUpSec9}
A(\alpha) = \left(
\begin{matrix}
h_\alpha - s\varepsilon_\alpha & s\beta_\alpha\\
s\beta_\alpha & \tilde h_\alpha - s\tilde\varepsilon_\alpha
\end{matrix}
\right)\hskip.1cm ,
\end{equation}
are fully coupled, and possess sequences of strictly positive solutions with either no bound on their 
$H_{1,2}^2$-norm, or which blow up with an 
arbitrary number of bubbles in their Sobolev decomposition 
and an arbitrary number of geometrical blow-up points.
\end{prop}

By a sequence of solutions which blows up with an 
arbitrary number of bubbles in its Sobolev decomposition 
and an arbitrary number of geometrical blow-up points, we mean 
here that the components 
of the solutions have an arbitrary number $m_i$ of bubbles in their $H_1^2$-decomposition 
(\ref{EqtThmSobTherSec4Bis}), and that 
the set of geometrical blow-up points can be chosen with any number $m \le \sum m_i$ of points. 
In particular, when $m < \sum m_i$, there are distinct bubbles which 
accumulate on the same point. Refinements where we 
also control the number of bubbles which accumulate on a given 
point are possible. Depending on whether $s = -1$ or 
$s = +1$, $-A(\alpha)$ 
or $A(\alpha)$ in (\ref{ExBlowUpSec9}) is cooperative.

\begin{proof} (i) We prove that there are $2$-systems like (\ref{GenericEqtSec9}), where $A(\alpha)$ 
is as in (\ref{ExBlowUpSec9}), which are fully coupled, and possess sequences of strictly positive solutions 
which blow up in any possible configuration. We fix $m_1$, $m_2$, $m$, and 
$s = \pm 1$. Let $q_1, q_2$ integers such that 
$m = q_1 + q_2$, $q_1 \le m_1$, and $q_2 \le m_2$. By the examples in Druet and Hebey \cite{DruHeb1}, see also \cite{DruHeb2}, 
there exist sequences $(h_\alpha)_\alpha$, $(\tilde h_\alpha)_\alpha$ of smooth functions, 
converging in $C^1(M)$ to $\frac{n-2}{4(n-1)}S_g$, and sequences $(u_\alpha)_\alpha$ and 
$(\tilde u_\alpha)_\alpha$ of positive solutions of the scalar equations
\begin{equation}\label{ProofPropSec9Eqt1}
\begin{split}
&\Delta_g u_\alpha + h_\alpha u_\alpha = u_\alpha^{2^\star-1}
\hskip.1cm ,\hskip.1cm\hbox{and}\\
&\Delta_g\tilde u_\alpha + \tilde h_\alpha \tilde u_\alpha = \tilde u_\alpha^{2^\star-1}
\end{split}
\end{equation}
such that the $u_\alpha$'s have $m_1$ bubbles in their $H_1^2$-decomposition, and 
$q_1$ geometrical blow-up points, and such that the $\tilde u_\alpha$'s have $m_2$ 
bubbles in their $H_1^2$-decomposition, and 
$q_2$ geometrical blow-up points. We let $(\beta_\alpha)_\alpha$ be a sequence of positive functions 
converging to $0$ in $C^1(M)$. We choose the $\beta_\alpha$'s sufficiently small in the $C^1$-norm such that 
$\beta_\alpha v_\alpha \to 0$ and $\beta_\alpha w_\alpha \to 0$  in $C^1(M)$ as $\alpha \to +\infty$, where 
$v_\alpha = u_\alpha/\tilde u_\alpha$ and $w_\alpha = \tilde u_\alpha / u_\alpha$.
Letting the $\varepsilon_\alpha$'s and $\tilde\varepsilon_\alpha$'s be given by
\begin{equation}\label{ProofPropSec9Eqt2}
\varepsilon_\alpha \frac{u_\alpha}{\tilde u_\alpha} 
= \tilde\varepsilon_\alpha \frac{\tilde u_\alpha}{u_\alpha} 
= \beta_\alpha
\end{equation}
for all $\alpha$, the sequences $(\varepsilon_\alpha)_\alpha$ and $(\tilde\varepsilon_\alpha)_\alpha$ 
consist of positive functions and converge to $0$ in $C^1(M)$. Moreover, by combining 
(\ref{ProofPropSec9Eqt1}) and (\ref{ProofPropSec9Eqt2}), we get that the ${\mathcal U}_\alpha$'s given by 
${\mathcal U}_\alpha = (u_\alpha,\tilde u_\alpha)$ for all $\alpha$, are solutions of the $2$-systems 
(\ref{GenericEqtSec9}), where the $A(\alpha)$'s are as in (\ref{ExBlowUpSec9}), and the sequences 
$(h_\alpha)_\alpha$, $(\tilde h_\alpha)_\alpha$, 
$(\varepsilon_\alpha)_\alpha$, $(\tilde\varepsilon_\alpha)_\alpha$, and $(\beta_\alpha)_\alpha$ are as above. 
The systems are fully coupled, the components of the ${\mathcal U}_\alpha$'s have $m_i$ bubbles 
in their Sobolev decomposition, and the set of geometrical blow-up points consists of $m$ points.
In particular, there 
are $2$-systems like (\ref{GenericEqtSec9}), where $A(\alpha)$ 
is as in (\ref{ExBlowUpSec9}), which are fully coupled, and possess sequences of strictly positive solutions 
which blow up with an 
arbitrary number of bubbles in their Sobolev decomposition 
and an arbitrary number of geometrical blow-up points..

\medskip\noindent (ii) We prove that there are $2$-systems like (\ref{GenericEqtSec9}), where $A(\alpha)$ 
is as in (\ref{ExBlowUpSec9}), which are fully coupled, and possess sequences of strictly positive solutions 
with no bound on their $H_{1,2}^2$-norm. 
By the examples in Druet and Hebey \cite{DruHeb1}, see also \cite{DruHeb2}, 
there are sequences $(h_\alpha)_\alpha$ of smooth functions, 
converging in $C^1(M)$ to $\frac{n-2}{4(n-1)}S_g$, and sequences $(u_\alpha)_\alpha$ 
of positive solutions of the scalar equations
\begin{equation}\label{ProofPropSec9Eqt3}
\Delta_g u_\alpha + h_\alpha u_\alpha = u_\alpha^{2^\star-1}
\hskip.1cm ,
\end{equation}
such that
$$\lim_{\alpha\to+\infty}\Vert u_\alpha\Vert_{H_1^2} = +\infty\hskip.1cm .$$
We let 
$(\varepsilon_\alpha)_\alpha$ be a sequence of positive functions 
converging to $0$ in $C^1(M)$. Letting 
$\tilde\varepsilon_\alpha = \beta_\alpha = \varepsilon_\alpha$, and 
$\tilde h_\alpha = h_\alpha$ for all $\alpha$, 
the ${\mathcal U}_\alpha$'s given by 
${\mathcal U}_\alpha = (u_\alpha,u_\alpha)$ for all $\alpha$, are solutions of the $2$-systems 
(\ref{GenericEqtSec9}), where the $A(\alpha)$'s are as in (\ref{ExBlowUpSec9}). 
The systems are fully coupled, and the components of the ${\mathcal U}_\alpha$'s have no bound 
on their $H_1^2$-norm. 
In particular, there are $2$-systems like (\ref{GenericEqtSec9}), where $A(\alpha)$ 
is as in (\ref{ExBlowUpSec9}), which are fully coupled, and possess sequences of strictly positive solutions 
with no bound on their $H_{1,2}^2$-norm. Together with (i), this ends the proof of the Proposition.
 \end{proof}

There are several notions of compactness in the literature. We distinguish two notions in this paper: {\it pseudo-compactness}, 
and {\it compactness}. The family (\ref{GenericEqtSec9}) is said to be {\it pseudo-compact} if 
for any 
bounded sequence $({\mathcal U}_\alpha)_\alpha$ in $H_{1,p}^2(M)$ of nontrivial nonnegative 
solutions of 
(\ref{GenericEqtSec9}) which converges 
weakly in $H_{1,p}^2(M)$, the weak limit ${\mathcal U}^0$ of the ${\mathcal U}_\alpha$'s is not zero, 
and thus is a nonnegative nontrivial solution of 
the limit system (\ref{LimitSystSec9}). 
In contrast, we say that the family (\ref{GenericEqtSec9}) is {\it compact} if any 
bounded sequence 
$({\mathcal U}_\alpha)_\alpha$ in $H_{1,p}^2(M)$ of nonnegative solutions 
of (\ref{GenericEqtSec9}) is actually bounded 
in $C^{2,\theta}(M)$, $0 < \theta < 1$, and thus converges, up to a subsequence, in $C^2(M)$ to some 
$p$-map ${\mathcal U}^0$, where the bound in $C^{2,\theta}(M)$ and the convergence in $C^2(M)$ 
have to be understood for the components of the ${\mathcal U}_\alpha$'s and ${\mathcal U}^0$. 
Because of the examples in Proposition \ref{PropExBlowUpSec9}, we need to assume 
a bound on the energy in these definitions.
 In terms of the $H_1^2$-decomposition
$${\mathcal U}_\alpha = {\mathcal U}^0 + \sum_{j=1}^k{\mathcal B}_{j,\alpha} + {\mathcal R}_\alpha$$
of the ${\mathcal U}_\alpha$'s, given by Theorem \ref{SobTheorThm}, pseudo-compactness reduces to  
${\mathcal U}^0 \not\equiv 0$, and, by the pointwise estimates 
in Section \ref{PointEst} and elliptic theory, compactness reduces to $k = 0$. 
Compactness 
is clearly a stronger notion than pseudo-compactness. 
Indeed, by (\ref{CondSec9}), we easily get that if $k = 0$, then ${\mathcal U}^0 \not\equiv 0$. In order to see this we 
can write, as in Section \ref{L2Conc}, that
\begin{eqnarray*} \int_M\vert{\mathcal U}_\alpha\vert^{2^\star}dv_g
&\ge& C_1 I_A({\mathcal U}_\alpha)\\
&\ge& C_2 \left(\int_M\vert{\mathcal U}_\alpha\vert^{2^\star}dv_g\right)^{2/2^\star}\hskip.1cm ,
\end{eqnarray*}
where $C_1, C_2 > 0$ are independent of $\alpha$, and $I_A$ is given by 
(\ref{FonctionalDef}). 
In particular, if the ${\mathcal U}_\alpha$'s are nonzero, then 
$$\Vert{\mathcal U}_\alpha\Vert_{2^\star} \ge C_3$$
for all $\alpha$ and 
some $C_3 > 0$ independent of $\alpha$, so that ${\mathcal U}^0 \not\equiv 0$ if $k = 0$. 
Pseudo-compactness provides the existence of a nontrivial 
solution of the limit system (\ref{LimitSystSec9}). Compactness is closely related 
to the a priori estimates one needs to prove to develop 
a Morse theory for systems like (\ref{GenericEqtSec9}).  For minimal energy sequences, 
as studied in Hebey \cite{Heb2}, pseudo-compactness 
and compactness 
are one and only one notion.
Compactness and a priori estimates for nonlinear (in general subcritical) systems in 
bounded domains of the Euclidean space, or for 
the critical Toda system on Riemann surfaces, are discussed in 
Angenent and Van der Vorst \cite{AngVan1,AngVan2}, Cl\'ement, Man\'asevich and Mitidieri \cite{CleManMit}, 
De Figueiredo \cite{DeF}, Jost, Lin and Wang \cite{JosLinWan}, 
Jost and Wang \cite{JosWan}, Montenegro \cite{Mon}, and Qing \cite{Qin}.

\medskip We know from the work of Schoen \cite{Sch2,Sch3,Sch4} 
that the Yamabe equation 
on manifolds 
conformally distinct to the sphere is compact. Related references are Druet \cite{Dru2}, 
Khuri and Schoen \cite{KhuSch}, Li and Zhu \cite{LiZhu}, Li and Zhang \cite{LiZha,LiZhabis}, 
and Marques \cite{Mar}. A very elegant 
proof of this compactness 
result is in Schoen \cite{Sch4} when the manifold is assumed to be conformally flat. 
In particular, it follows from this result that for manifolds 
distinct to the sphere, 
the trivially coupled system consisting of $p$-copies of the Yamabe equation is compact. 
We prove here that 
our systems (\ref{GenericEqtSec9}) are also compact when their coupling stands far from this trivial coupling. 
By far we mean here that one of the three following 
situations occur:

\medskip (A1) $A - \frac{n-2}{4(n-1)}S_gId_p$ has a sign, or

\medskip (A2) $A_{ii} < \frac{n-2}{4(n-1)}S_g$ for all $i$, or

\medskip (A3) $A_{ii} > \frac{n-2}{4(n-1)}S_g$ for all $i$,

\medskip\noindent where $A$ is the limit map in (\ref{CondSec9}), $S_g$ is the scalar curvature 
of $g$, $Id_p$ is the $p\times p$-identity matrix, and, for $B: M \to M_p^s({\mathbb R})$, we say that 
$B$ has a sign if either $B(x).(X,X) < 0$
for any $x \in M$ and any $X \in {\mathbb R}^p\backslash\{0\}$, or $B(x).(X,X) > 0$
for any $x \in M$ and any $X \in {\mathbb R}^p\backslash\{0\}$. 
Because of the dynamical 
viewpoint we adopted here, and the examples in Proposition \ref{PropExBlowUpSec9}, we need to assume 
conditions like (A1), (A2), or (A3) in order to get pseudo-compactness or compactness. Actually, it is easily seen 
that we also need to 
assume more than (A2) or (A3), as indicated in the remark below. 

\medskip\noindent{\bf Remark 9.1:}  
Conditions like (A2) or (A3) alone are not 
sufficient to guarantee compactness 
(or pseudo-compactness). Let, for instance, $(u_\alpha)_\alpha$ be a blowing-up sequence of solutions of the 
Yamabe equation on the sphere, and 
${\mathcal U}_\alpha = (u_\alpha,u_\alpha)$ for all $\alpha$. Let also $a$, $b$, $c$ be 
real numbers such that
$$a + b = b + c = \frac{n(n-2)}{4}\hskip.1cm ,$$
and $a, c > 0$.
Then $({\mathcal U}_\alpha)_\alpha$ is a blowing-up sequence of solutions of the $2$-system
\begin{eqnarray*}
&&\Delta_{g_0}u + au + bv = u^{2^\star-1}\\
&&\Delta_{g_0}v + bu + cv = v^{2^\star-1}\hskip.1cm ,
\end{eqnarray*}
and we can choose $a$, $b$, and $c$ such that the matrix in the system is coercive and 
either $a, c < \frac{n(n-2)}{4}$, or $a, c > \frac{n(n-2)}{4}$. In the first case $b > 0$, 
and $A$ is cooperative. In the second case, $b < 0$ and $-A$ is cooperative. 

\medskip Our compactness result states as follows. With respect to the examples 
we just discussed in Remark 9.1, the assertions in Theorem \ref{ThmCptnessSec9} that 
the $-A(\alpha)$'s should be cooperative when we assume 
(A2), or that the $A(\alpha)$'s should be cooperative when we assume (A3), are sharp. 

\begin{thm}[Compactness]\label{ThmCptnessSec9} Let $(M,g)$ be a smooth compact conformally flat Riemannian 
manifold of dimension 
$n \ge 4$, $p \ge 1$, and $\left(A(\alpha)\right)_\alpha$ a sequence of smooth maps 
$A(\alpha): M \to M_p^s({\mathbb R})$ satisfying (\ref{CondSec9}). If (A1) 
holds, then the family (\ref{GenericEqtSec9}) is pseudo-compact when $n \ge 4$, and compact when 
$n \ge 7$. If (A2) holds, and the $-A(\alpha)$'s are cooperative for all $\alpha$, 
then the family (\ref{GenericEqtSec9}) is pseudo-compact when $n \ge 4$, and compact 
when $n \ge 7$. It is 
also compact when $n \ge 4$ if we assume in addition 
that the limit system (\ref{LimitSystSec9}) is fully coupled. At last, if (A3) holds, and the $A(\alpha)$'s are cooperative 
for all $\alpha$, then the family (\ref{GenericEqtSec9}) is pseudo-compact when $n \ge 4$, and 
compact when $n \ge 7$.
\end{thm}

We prove Theorem \ref{ThmCptnessSec9} by proving first the pseudo-compactness part, 
and then the compactness part of the theorem. For both these parts, the proof is 
by contradiction. We use in the process the blow-up theory developed 
in Sections \ref{SobTheory} to \ref{StdRescaling}. We also use conformal invariance, and an 
independent relation given 
by the Euclidean Pohozaev identity (\ref{PohozaevSec9}) below. The key estimate for 
pseudo-compactness is $L^2$-concentration in Section \ref{L2Conc}. The key estimate for compactness 
is the sharp asymptotic in Section \ref{SharpPtAsympt}. For $\Omega$ a smooth bounded domain 
in the Euclidean space, and $u$ a smooth function in $\overline{\Omega}$, the Pohozaev 
identity we use in the sequel reads as 
\begin{equation}\label{PohozaevSec9}
\begin{split}
&\int_\Omega(x^k\partial_ku)\Delta u dx + \frac{n-2}{2} \int_\Omega u(\Delta u)dx\\
&= - \int_{\partial\Omega}(x^k\partial_ku)\partial_\nu ud\sigma 
+ \frac{1}{2} \int_{\partial\Omega}(x,\nu)\vert\nabla u\vert^2d\sigma \\
&\hskip.4cm 
- \frac{n-2}{2}\int_{\partial\Omega}u\partial_\nu ud\sigma\hskip.1cm ,
\end{split}
\end{equation}
where $\Delta$ is the Euclidean Laplacian, $\nu$ is the outward unit normal to $\partial\Omega$,  
$d\sigma$ is the Euclidean volume element on $\partial\Omega$, and 
there is a sum over $k$ from $1$ to $n$. We easily get 
(\ref{PohozaevSec9}) by integrating by parts the first term in the left hand side of (\ref{PohozaevSec9}).
As one can check from the proof of the theorem, 
see below, the pseudo-compactness part in Theorem \ref{ThmCptnessSec9}, and the compactness 
part when $n \ge 7$, are also true if we assume that (A2) holds and $-A$ is cooperative, or 
that (A3) holds and $A$ is cooperative, where $A$ is the limit of the $A(\alpha)$'s in 
(\ref{CondSec9}). 

\begin{proof}[Proof of Theorem \ref{ThmCptnessSec9}] We prove Theorem \ref{ThmCptnessSec9} 
by contradiction, and, as already mentionned, we proceed in two steps 
by proving first the pseudo-compactness part of the theorem, 
and then the compactness part of the theorem. In what follows, we let 
$({\mathcal U}_\alpha)_\alpha$ be a bounded sequence in $H_{1,p}^2(M)$ of nontrivial nonnegative 
solutions of (\ref{GenericEqtSec9}). We assume that the ${\mathcal U}_\alpha$'s blow up. 
Up to a subsequence, we may then write that 
\begin{equation}\label{SobDecSec9}
{\mathcal U}_\alpha = {\mathcal U}^0 + \sum_{j=1}^k{\mathcal B}_{j,\alpha} + {\mathcal R}_\alpha\hskip.1cm ,
\end{equation}
where ${\mathcal U}^0$ is a weak nonnegative solution of the limit system 
(\ref{LimitSystSec9}), the $({\mathcal B}_{j,\alpha})_\alpha$'s are $p$-bubbles, and the 
${\mathcal R}_\alpha$'s converge strongly to $0$ in $H_{1,p}^2(M)$ as $\alpha \to +\infty$. 
 First  we prove that, under the assumptions concerning 
pseudo-compactness in Theorem \ref{ThmCptnessSec9}, we have that ${\mathcal U}^0 \not\equiv 0$ 
in (\ref{SobDecSec9}). Then we prove that, under the assumptions concerning 
compactness in Theorem \ref{ThmCptnessSec9}, we also get that $k = 0$ in (\ref{SobDecSec9}). We let 
${\mathcal S}_{geom}$ be the set consisting of the geometrical blow-up points of 
the sequence $({\mathcal U}_\alpha)_\alpha$. We let also the 
$u_\alpha^i$'s be the components of the ${\mathcal U}_\alpha$'s. In the first case, when proving 
pseudo-compactness, we apply (\ref{PohozaevSec9}) to the $u_\alpha^i$'s 
on balls $B_{x_i}(\delta)$, $\delta > 0$, where $x_i \in {\mathcal S}_{geom}$. 
In the second case, when proving 
compactness, we apply (\ref{PohozaevSec9}) on the smaller balls $B_{x_\alpha}(\delta\sqrt{\mu_\alpha})$, 
$\delta > 0$, where $\mu_\alpha$ stand for 
the largest weight among all the possible weights of the bubbles in (\ref{SobDecSec9}), and 
$x_\alpha$ stand for the corresponding center in (\ref{SobDecSec9}).

\begin{proof}[Proof of pseudo-compactness] We assume that ${\mathcal U}^0 \equiv 0$ in 
(\ref{SobDecSec9}). 
Let $x_0 \in {\mathcal S}_{geom}$. Since $g$ is conformally flat, 
there exist $\delta > 0$ and a conformal metric $\hat g$ to $g$ such that $\hat g$ is flat in $B_{x_0}(4\delta)$. 
Let $\hat g = \varphi^{4/(n-2)}g$, where $\varphi$ is smooth and positive, and 
$\hat u_\alpha^i = u_\alpha^i\varphi^{-1}$ for all $\alpha$ and $i$. By conformal invariance of the 
conformal Laplacian, see, for instance, Lee and Parker \cite{LeePar}, and by (\ref{GenericEqtSec9}), 
$\hat{\mathcal U}_\alpha = (\hat u_\alpha^1,\dots,\hat u_\alpha^p)$ is a solution of the system
\begin{equation}\label{ProofPseudoCptEqt1}
\Delta\hat u_\alpha^i + \sum_{j=1}^p\hat A^\alpha_{ij}\hat u_\alpha^j = (\hat u_\alpha^i)^{2^\star-1}
\end{equation}
in $B_{x_0}(4\delta)$ for all $i$ and all $\alpha$, where $\Delta = \Delta_{\hat g}$ is the Euclidean Laplacian, 
\begin{equation}\label{ProofPseudoCptEqt2}
\hat A^\alpha_{ij} = \frac{1}{\varphi^{2^\star-2}} 
\left(A^\alpha_{ij} - \frac{n-2}{4(n-1)}S_g\delta_{ij}\right)\hskip.1cm ,
\end{equation}
and $S_g$ is the scalar curvature of $g$. We choose $\delta > 0$ sufficiently small such that 
${\mathcal S}_{geom}\cap B_{x_0}(4\delta) = \{x_0\}$. We regard the $\hat u_\alpha^i$'s, $\varphi$, and the 
$\hat A^\alpha_{ij}$'s as defined in the Euclidean space. Also we assimilate $x_0$ with $0$. 
The $u_\alpha^i$'s, see Section \ref{PointEst}, 
are uniformly bounded in $C^0_{loc}(M\backslash{\mathcal S}_{geom})$. 
By the De Giorgi-Nash-Moser iterative scheme 
in Section \ref{L2Conc}, it follows that  the $C^0$-norm of the $\hat u_\alpha^i$'s in small 
neighbourhood of 
$\partial B_0(\delta)$ are controlled by the $L^2$-norm of the 
$\vert\hat{\mathcal U}_\alpha\vert$'s in 
annuli like $B_0(2\delta)\backslash B_0(\delta/2)$. By standard elliptic theory, as developed 
in Gilbarg and Trudinger \cite{GilTru}, we then get that 
for $T_\delta$ a sufficiently small neighbourhood of $\partial B_0(\delta)$,
\begin{equation}\label{ProofPseudoCptEqt3}
\left\Vert\hat u_\alpha^i\right\Vert_{C^{1,\theta}(T_\delta)}^2 \le 
C\int_{B_0(2\delta)\backslash B_0(\delta/2)}\vert\hat {\mathcal U}_\alpha\vert^2dx
\end{equation}
for all $\alpha$ and $i$, where $C > 0$ does not depend on $\alpha$, 
$\vert\hat {\mathcal U}_\alpha\vert^2 = \sum_i(\hat u_\alpha^i)^2$, and $\theta \in (0,1)$. 
Let
\begin{equation}\label{ProofPseudoCptEqt4}
{\mathcal R}_\alpha^\delta = \int_{B_0(2\delta)\backslash B_0(\delta/2)}\vert\hat {\mathcal U}_\alpha\vert^2dx
\hskip.1cm .
\end{equation}
By (\ref{ProofPseudoCptEqt3}), 
plugging the $\hat u_\alpha^i$'s in the Pohozaev identity (\ref{PohozaevSec9}), with $\Omega = B_0(\delta)$, 
we can write that
\begin{equation}\label{ProofPseudoCptEqt5}
\int_{B_0(\delta)}(x^k\partial_k\hat u_\alpha^i)\Delta\hat u_\alpha^i dx 
+ \frac{n-2}{2} \int_{B_0(\delta)}\hat u_\alpha^i(\Delta\hat u_\alpha^i)dx 
= O\left({\mathcal R}_\alpha^\delta\right)
\end{equation}
for all $\alpha$ and $i$, 
where ${\mathcal R}_\alpha^\delta$ is as in (\ref{ProofPseudoCptEqt4}). 
Combining (\ref{ProofPseudoCptEqt1}) and (\ref{ProofPseudoCptEqt5}), summing over 
$i=1,\dots,p$, we then get that
\begin{equation}\label{ProofPseudoCptEqt6}
\begin{split}
&\sum_{i=1}^p\int_{B_0(\delta)}(x^k\partial_k\hat u_\alpha^i)(\hat u_\alpha^i)^{2^\star-1}dx
-\sum_{i,j=1}^p\int_{B_0(\delta)}(x^k\partial_k\hat u_\alpha^i)\hat A^\alpha_{ij}\hat u_\alpha^j dx\\
&+\frac{n-2}{2} \sum_{i=1}^p\int_{B_0(\delta)}(\hat u_\alpha^i)^{2^\star}dx 
-\frac{n-2}{2}  \sum_{i,j=1}^p\int_{B_0(\delta)}\hat A^\alpha_{ij}\hat u_\alpha^i\hat u_\alpha^j dx\\
&= O\left({\mathcal R}_\alpha^\delta\right)
\end{split}
\end{equation}
for all $\alpha$, where the $\hat A^\alpha_{ij}$'s are as in (\ref{ProofPseudoCptEqt2}), and 
${\mathcal R}_\alpha^\delta$ is as in (\ref{ProofPseudoCptEqt4}). Integrating by parts, we can write 
with (\ref{ProofPseudoCptEqt3}) that 
\begin{equation}\label{ProofPseudoCptEqt7}
\begin{split}
&\sum_{i,j=1}^p\int_{B_0(\delta)}(x^k\partial_k\hat u_\alpha^i)\hat A^\alpha_{ij}\hat u_\alpha^j dx\\
&= -\frac{n}{2}\sum_{i,j=1}^p\int_{B_0(\delta)}\hat A^\alpha_{ij}\hat u_\alpha^i\hat u_\alpha^j dx\\
&\hskip.4cm - \frac{1}{2} 
\sum_{i,j=1}^p\int_{B_0(\delta)}(x^k\partial_k\hat A^\alpha_{ij})\hat u_\alpha^i\hat u_\alpha^j dx
+ O\left({\mathcal R}_\alpha^\delta\right)
\end{split}
\end{equation}
for all $\alpha$, where ${\mathcal R}_\alpha^\delta$ is as in (\ref{ProofPseudoCptEqt4}). We can also write that 
\begin{equation}\label{ProofPseudoCptEqt8}
\begin{split}
&\sum_{i=1}^p\int_{B_0(\delta)}(x^k\partial_k\hat u_\alpha^i)(\hat u_\alpha^i)^{2^\star-1}dx\\
&= - \frac{n-2}{2}\sum_{i=1}^p\int_{B_0(\delta)}(\hat u_\alpha^i)^{2^\star}dx + O\left({\mathcal R}_\alpha^\delta\right)
\end{split}
\end{equation}
for all $\alpha$, where ${\mathcal R}_\alpha^\delta$ is as in (\ref{ProofPseudoCptEqt4}). In particular, 
plugging (\ref{ProofPseudoCptEqt7}) and (\ref{ProofPseudoCptEqt8}) in (\ref{ProofPseudoCptEqt6}), it 
follows that 
\begin{equation}\label{ProofPseudoCptEqt9}
\begin{split}
&\sum_{i,j=1}^p\int_{B_0(\delta)}\hat A^\alpha_{ij}\hat u_\alpha^i\hat u_\alpha^j dx
+ \frac{1}{2} \sum_{i,j=1}^p\int_{B_0(\delta)}(x^k\partial_k\hat A^\alpha_{ij})\hat u_\alpha^i\hat u_\alpha^j dx\\
&= O\left({\mathcal R}_\alpha^\delta\right)
\end{split}
\end{equation}
for all $\alpha$, where ${\mathcal R}_\alpha^\delta$ is as in (\ref{ProofPseudoCptEqt4}). By the $C^1$-convergence 
in (\ref{CondSec9}), the $\partial_k\hat A^\alpha_{ij}$'s are uniformly bounded. Coming back to the manifold, 
noting that $dv_{\hat g} = \varphi^{2^\star}dv_g$, and 
summing over the $x_0$ in ${\mathcal S}_{geom}$, we get with (\ref{CondSec9}), (\ref{ProofPseudoCptEqt9}), 
and $L^2$-concentration in Section \ref{L2Conc}, that
\begin{equation}\label{ProofPseudoCptEqt10}
\begin{split}
&\sum_{i,j=1}^p\int_{{\mathcal B}_\delta}
\left(A_{ij}-\frac{n-2}{4(n-1)}S_g\delta_{ij}\right)
u_\alpha^iu_\alpha^jdv_g\\
&= \varepsilon_\delta O\left(\sum_{i=1}^p\int_M(u_\alpha^i)^2dv_g\right)
+ o\left(\sum_{i=1}^p\int_M(u_\alpha^i)^2dv_g\right)
\end{split}
\end{equation}
for all $\alpha$, where $\varepsilon_\delta \to 0$ as $\delta \to 0$, the first term in the right hand side of 
(\ref{ProofPseudoCptEqt10}) depends on $\delta$ only by $\varepsilon_\delta$, 
the $A_{ij}$'s are the components of the limit map $A$ given by 
(\ref{CondSec9}), and 
$${\mathcal B}_\delta = \bigcup_{x \in {\mathcal S}_{geom}}B_x(\delta)\hskip.1cm .$$
If we assume that (A1) holds, then either 
$A_{ij} > \frac{n-2}{4(n-1)}S_g\delta_{ij}$ or 
$A_{ij} < \frac{n-2}{4(n-1)}S_g\delta_{ij}$ in the sense of bilinear forms. In both cases, it follows from 
(\ref{ProofPseudoCptEqt10}) and $L^2$-concentration that there exists $C > 0$, independent of $\alpha$ and $\delta$, such that 
\begin{equation}\label{ProofPseudoCptConclusEqt}
C\sum_{i=1}^p\int_M(u_\alpha^i)^2dv_g
\le  \varepsilon_\delta O\left(\sum_{i=1}^p\int_M(u_\alpha^i)^2dv_g\right)
+ o\left(\sum_{i=1}^p\int_M(u_\alpha^i)^2dv_g\right)
\end{equation}
for all $\alpha$, where $\varepsilon_\delta \to 0$ as $\delta \to 0$. If we assume that the $A(\alpha)$'s are cooperative, 
then, by (\ref{CondSec9}), $A$ is also cooperative. In particular, $A_{ij}u_\alpha^iu_\alpha^j \ge 0$ for 
$i\not= j$, and if we assume in addition that (A3) holds, then we get here again 
that there exists $C > 0$, independent of $\alpha$ and $\delta$, such that (\ref{ProofPseudoCptConclusEqt}) 
holds for all $\alpha$. The same conclusion (\ref{ProofPseudoCptConclusEqt}) 
holds if we assume (A2) and that $-A$ is, or that the $-A(\alpha)$'s are, cooperative. 
Taking $\delta > 0$ sufficiently small, the contradiction follows from (\ref{ProofPseudoCptConclusEqt}). In particular, ${\mathcal U}^0 \not\equiv 0$ in (\ref{SobDecSec9}), and this 
ends the proof of the 
pseudo-compactness part of Theorem \ref{ThmCptnessSec9}.
\end{proof}

\begin{proof}[Proof of compactness] We let the $x_\alpha$'s and $\mu_\alpha$'s be given by (\ref{HighestWeightSec7}) 
and (\ref{HighestWeightBisSec7}) in Section \ref{SharpPtAsympt}. We let also $\delta > 0$ small be less than 
the $\delta$ given by 
Theorem \ref{SharpPtwseAsymptThm}.
Since $g$ is conformally flat, 
there exists a conformal metric $\hat g$ to $g$ such that $\hat g$ is flat in $B_{x_0}(4\delta)$, where 
$x_0$ is the limit of the $x_\alpha$'s. 
We let $\hat g = \varphi^{4/(n-2)}g$, where $\varphi$ is smooth, positive, and such that $\varphi(x_0) = 1$, and 
let $\hat u_\alpha^i = u_\alpha^i\varphi^{-1}$ for all $\alpha$ and $i$. 
By conformal invariance of the conformal Laplacian, and by (\ref{GenericEqtSec9}), 
equation (\ref{ProofPseudoCptEqt1}) holds in $B_{x_0}(4\delta)$ for all $\alpha$ and all $i$. Namely,
$$\Delta\hat u_\alpha^i + \sum_{j=1}^p\hat A^\alpha_{ij}\hat u_\alpha^j = (\hat u_\alpha^i)^{2^\star-1}$$
in $B_{x_0}(4\delta)$ for all $i$ and all $\alpha$, where $\Delta = \Delta_{\hat g}$ is the Euclidean Laplacian, 
and the $\hat A^\alpha_{ij}$'s are given by (\ref{ProofPseudoCptEqt2}). 
In what follows 
we assimilate $x_\alpha$ with $0 \in {\mathbb R}^n$ (thanks to the exponential map 
at $x_\alpha$). We also regard the $\hat u_\alpha^i$'s, $\varphi$, and the 
$\hat A^\alpha_{ij}$'s in (\ref{ProofPseudoCptEqt1}) as defined in the Euclidean space. 
Plugging the $\hat u_\alpha^i$'s in the Pohozaev identity (\ref{PohozaevSec9}), with 
$\Omega = B_0(\delta\mu_\alpha)$, 
we get that
\begin{equation}\label{ProofCptEqt1}
\begin{split}
&\int_{B_0(\delta\sqrt{\mu_\alpha})}(x^k\partial_k\hat u_\alpha^i)\Delta\hat u_\alpha^i dx 
+ \frac{n-2}{2} \int_{B_0(\delta\sqrt{\mu_\alpha})}\hat u_\alpha^i(\Delta\hat u_\alpha^i)dx\\
&= - \int_{\partial B_0(\delta\sqrt{\mu_\alpha})}(x^k\partial_k\hat u_\alpha^i)\partial_\nu\hat u_\alpha^id\sigma 
+ \frac{1}{2} \int_{\partial B_0(\delta\sqrt{\mu_\alpha})}(x,\nu)\vert\nabla\hat u_\alpha^i\vert^2d\sigma \\
&\hskip.4cm 
- \frac{n-2}{2}\int_{\partial B_0(\delta\sqrt{\mu_\alpha})}\hat u_\alpha^i\partial_\nu\hat u_\alpha^id\sigma
\end{split}
\end{equation}
for all $\alpha$ and $i$. 
Combining (\ref{ProofPseudoCptEqt1}) and (\ref{ProofCptEqt1}), summing over $i = 1,\dots,p$, 
it follows that
\begin{equation}\label{ProofCptEqt2}
\begin{split}
&\sum_{i=1}^p\int_{B_0(\delta\sqrt{\mu_\alpha})}(x^k\partial_k\hat u_\alpha^i)(\hat u_\alpha^i)^{2^\star-1}dx 
- \sum_{i,j=1}^p\int_{B_0(\delta\sqrt{\mu_\alpha})}(x^k\partial_k\hat u_\alpha^i)\hat A^\alpha_{ij}\hat u_\alpha^jdx\\
&+ \frac{n-2}{2} \sum_{i=1}^p\int_{B_0(\delta\sqrt{\mu_\alpha})}(\hat u_\alpha^i)^{2^\star}dx
- \frac{n-2}{2}\sum_{i,j=1}^p\int_{B_0(\delta\sqrt{\mu_\alpha})}\hat A^\alpha_{ij}\hat u_\alpha^i\hat u_\alpha^jdx\\
&= - \sum_{i=1}^p\int_{\partial B_0(\delta\sqrt{\mu_\alpha})}(x^k\partial_k\hat u_\alpha^i)\partial_\nu\hat u_\alpha^id\sigma 
+ \frac{1}{2} \sum_{i=1}^p \int_{\partial B_0(\delta\sqrt{\mu_\alpha})}(x,\nu)\vert\nabla\hat u_\alpha^i\vert^2d\sigma \\
&\hskip.4cm 
- \frac{n-2}{2} \sum_{i=1}^p \int_{\partial B_0(\delta\sqrt{\mu_\alpha})}\hat u_\alpha^i\partial_\nu\hat u_\alpha^id\sigma
\end{split}
\end{equation}
for all $\alpha$. 
Integrating by parts,
\begin{equation}\label{ProofCptEqt3}
\begin{split}
&\sum_{i=1}^p\int_{B_0(\delta\sqrt{\mu_\alpha})}(x^k\partial_k\hat u_\alpha^i)(\hat u_\alpha^i)^{2^\star-1}dx\\
&= - \frac{n-2}{2} \sum_{i=1}^p\int_{B_0(\delta\sqrt{\mu_\alpha})}(\hat u_\alpha^i)^{2^\star}dx
+ \frac{n-2}{2n} \sum_{i=1}^p \int_{\partial B_0(\delta\sqrt{\mu_\alpha})}(x,\nu)(\hat u_\alpha^i)^{2^\star}d\sigma
\hskip.1cm ,
\end{split}
\end{equation}
and
\begin{equation}\label{ProofCptEqt4}
\begin{split}
&\sum_{i,j=1}^p\int_{B_0(\delta\sqrt{\mu_\alpha})}(x^k\partial_k\hat u_\alpha^i)\hat A^\alpha_{ij}\hat u_\alpha^jdx\\
&= - \frac{n}{2}\sum_{i,j=1}^p\int_{B_0(\delta\sqrt{\mu_\alpha})} \hat A^\alpha_{ij}\hat u_\alpha^i\hat u_\alpha^jdx 
- \frac{1}{2} \sum_{i,j=1}^p\int_{B_0(\delta\sqrt{\mu_\alpha})}(x^k\partial_k\hat A^\alpha_{ij})\hat u_\alpha^i\hat u_\alpha^jdx\\
&\hskip.4cm + \frac{1}{2}
\sum_{i,j=1}^p\int_{\partial B_0(\delta\sqrt{\mu_\alpha})} (x,\nu) \hat A^\alpha_{ij}\hat u_\alpha^i\hat u_\alpha^jd\sigma 
\end{split}
\end{equation}
for all $\alpha$. By (\ref{CondSec9}), we can write that
\begin{eqnarray*}
&&\int_{\partial B_0(\delta\sqrt{\mu_\alpha})}(x,\nu)(\hat u_\alpha^i)^{2^\star}d\sigma = 
o\left(\int_{\partial B_0(\delta\sqrt{\mu_\alpha})}\vert\hat{\mathcal U}_\alpha\vert^{2^\star}d\sigma\right)
\hskip.1cm ,\hskip.1cm\hbox{and}\\
&&\int_{\partial B_0(\delta\sqrt{\mu_\alpha})} (x,\nu) \hat A^\alpha_{ij}\hat u_\alpha^i\hat u_\alpha^jdx = 
o\left(\int_{\partial B_0(\delta\sqrt{\mu_\alpha})}\vert\hat{\mathcal U}_\alpha\vert^2d\sigma\right)
\end{eqnarray*}
for all $\alpha$, $i$, and $j$, where $\vert\hat{\mathcal U}_\alpha\vert^q = \sum_i\vert\hat u_\alpha^i\vert^q$. By Theorem 
\ref{SharpPtwseAsymptThm}, and the change of variables $x = \sqrt{\mu_\alpha}y$, we also have that
$$\int_{\partial B_0(\delta\sqrt{\mu_\alpha})}\vert\hat{\mathcal U}_\alpha\vert^qd\sigma = 
O\left(\mu_\alpha^{\frac{n-1}{2}}\right)$$
for all $\alpha$, and all $q = 2, 2^\star$. Hence,
\begin{equation}\label{ProofCptEqt5}
\begin{split}
&\int_{\partial B_0(\delta\sqrt{\mu_\alpha})}(x,\nu)(\hat u_\alpha^i)^{2^\star}d\sigma = 
o\left(\mu_\alpha^{\frac{n-1}{2}}\right)
\hskip.1cm ,\hskip.1cm\hbox{and}\\
&\int_{\partial B_0(\delta\sqrt{\mu_\alpha})} (x,\nu) \hat A^\alpha_{ij}\hat u_\alpha^i\hat u_\alpha^jdx = 
o\left(\mu_\alpha^{\frac{n-1}{2}}\right)
\end{split}
\end{equation}
for all $\alpha$, $i$, and $j$. By the $C^1$-convergence in (\ref{CondSec9}),
\begin{equation}\label{ProofCptEqt6}
\int_{B_0(\delta\sqrt{\mu_\alpha})}(x^k\partial_k\hat A^\alpha_{ij})\hat u_\alpha^i\hat u_\alpha^jdx
= o\left(\int_{B_0(\delta\sqrt{\mu_\alpha})}\vert\hat{\mathcal U}_\alpha\vert^2dx\right)
\end{equation}
for all $\alpha$, $i$, and $j$. By (\ref{ProofCptEqt2})--(\ref{ProofCptEqt6}), we then get that 
\begin{equation}\label{ProofCptEqt7}
\begin{split}
&\sum_{i,j=1}^p\int_{B_0(\delta\sqrt{\mu_\alpha})} \hat A^\alpha_{ij}\hat u_\alpha^i\hat u_\alpha^jdx 
+ o\left(\int_{B_0(\delta\sqrt{\mu_\alpha})}\vert\hat{\mathcal U}_\alpha\vert^2dx\right)
+ o\left(\mu_\alpha^{\frac{n-1}{2}}\right)\\
&= - \sum_{i=1}^p\int_{\partial B_0(\delta\sqrt{\mu_\alpha})}(x^k\partial_k\hat u_\alpha^i)\partial_\nu\hat u_\alpha^id\sigma 
+ \frac{1}{2} \sum_{i=1}^p \int_{\partial B_0(\delta\sqrt{\mu_\alpha})}(x,\nu)\vert\nabla\hat u_\alpha^i\vert^2d\sigma \\
&\hskip.4cm 
- \frac{n-2}{2} \sum_{i=1}^p \int_{\partial B_0(\delta\sqrt{\mu_\alpha})}\hat u_\alpha^i\partial_\nu\hat u_\alpha^id\sigma
\end{split}
\end{equation}
for all $\alpha$. In what follows we let ${\mathcal R}_{HS}(\alpha)$ 
be the right hand side in (\ref{ProofCptEqt7}). By the change of variables $x = \sqrt{\mu_\alpha}y$, by 
Theorem \ref{SharpPtwseAsymptThm}, and since we assumed that $\varphi(x_0) = 1$, we can write that
\begin{equation}\label{ProofCptEqt8}
\lim_{\alpha\to+\infty}\mu_\alpha^{-\frac{n-2}{2}}{\mathcal R}_{HS}(\alpha)  = {\mathcal R}\hskip.1cm ,
\end{equation}
where
\begin{equation}\label{ProofCptEqt9}
\begin{split}
{\mathcal R}
&= - \sum_{i=1}^p\int_{\partial B_0(\delta)}(x^k\partial_k\tilde u_i)\partial_\nu\tilde u_id\sigma 
+ \frac{1}{2} \sum_{i=1}^p \int_{\partial B_0(\delta)}(x,\nu)\vert\nabla\tilde u_i\vert^2d\sigma \\
&\hskip.4cm 
- \frac{n-2}{2} \sum_{i=1}^p \int_{\partial B_0(\delta)}\tilde u_i\partial_\nu\tilde u_id\sigma\hskip.1cm ,
\end{split}
\end{equation}
and where, for $i = 1,\dots,p$, the $\tilde u_i$'s are the limits of the 
$\tilde u_\alpha^i$'s given by Theorem \ref{SharpPtwseAsymptThm}. In particular,
\begin{equation}\label{ProofCptEqt10}
\tilde u_i(x) = \frac{A_i}{\vert x\vert^{n-2}} + \varphi_i(x)
\end{equation}
for all $i$, where $A_i \ge 0$, and $\varphi_i$ is harmonic in 
an open ball $B_0(\delta^\prime)$, $\delta^\prime > \delta$. We clearly have that $\Delta\tilde u_i = 0$ 
in $B_0(\delta^\prime)\backslash\{0\}$. Applying the Pohozaev identity (\ref{PohozaevSec9}) to 
the $\tilde u_i$'s in annuli like $B_0(\delta)\backslash B_0(r)$, where $0 < r < \delta$, 
and letting $r \to 0$, it follows that
\begin{equation}\label{ProofCptEqt11}
\begin{split}
{\mathcal R}
&= \lim_{r \to 0}\Biggl[- \sum_{i=1}^p\int_{\partial B_0(r)}(x^k\partial_k\tilde u_i)\partial_\nu\tilde u_id\sigma 
+ \frac{1}{2} \sum_{i=1}^p \int_{\partial B_0(r)}(x,\nu)\vert\nabla\tilde u_i\vert^2d\sigma \\
&\hskip 1cm 
- \frac{n-2}{2} \sum_{i=1}^p \int_{\partial B_0(r)}\tilde u_i\partial_\nu\tilde u_id\sigma\Biggr]\hskip.1cm ,
\end{split}
\end{equation}
and by (\ref{ProofCptEqt10}), we can write with (\ref{ProofCptEqt11}) that
${\mathcal R} = \frac{(n-2)^2}{2} \omega_{n-1}\sum_{i=1}^pA_i\varphi_i(0)$, 
where $\omega_{n-1}$ is the volume of the unit $(n-1)$-sphere.
We have that $dv_{\hat g} = \varphi^{2^\star}dv_g$. By (\ref{ProofCptEqt7})--(\ref{ProofCptEqt8}), and also by (\ref{CondSec9}), we then get, 
by coming back to the manifold, that
\begin{equation}\label{ProofCptConcludEqt1}
\begin{split}
&\sum_{i,j=1}^p\int_{B_{x_\alpha}(\delta\sqrt{\mu_\alpha})} \left(A_{ij} - \frac{n-2}{4(n-1)}S_g\delta_{ij}\right) u_\alpha^i u_\alpha^jdv_g\\
&+ o\left(\int_{B_{x_\alpha}(\delta\sqrt{\mu_\alpha})}\vert{\mathcal U}_\alpha\vert^2dv_g\right) 
= \frac{(n-2)^2}{2} \left(\sum_{i=1}^pA_i\varphi_i(0) + o(1)\right) \omega_{n-1} \mu_\alpha^{\frac{n-2}{2}} 
\end{split}
\end{equation}
for all $\alpha$, 
where the $A_i$'s and $\varphi_i$'s are given by Theorem \ref{SharpPtwseAsymptThm}, 
the $A_{ij}$'s are the components of the limit map $A$ given by (\ref{CondSec9}), 
$o(1) \to 0$ as $\alpha \to +\infty$, and $\omega_{n-1}$ is the volume of the unit $(n-1)$-sphere. 
Now, if we assume that (A1) holds, then either 
$A_{ij} > \frac{n-2}{4(n-1)}S_g\delta_{ij}$ or 
$A_{ij} < \frac{n-2}{4(n-1)}S_g\delta_{ij}$ in the sense of bilinear forms. In both cases, it follows from 
(\ref{ProofCptConcludEqt1}) that there exists $C > 0$, independent of $\alpha$, such that 
\begin{equation}\label{ProofCptConcludEqt2}
C\sum_{i=1}^p\int_{B_{x_\alpha}(\delta\sqrt{\mu_\alpha})}(u_\alpha^i)^2dv_g
=  O\left(\mu_\alpha^{\frac{n-2}{2}}\right)
\end{equation}
for all $\alpha$. When $n \ge 7$, $\frac{n-2}{2} > 2$ so that $\mu_\alpha^{\frac{n-2}{2}} = o\left(\mu_\alpha^2\right)$, 
and the contradiction follows 
from (\ref{ProofCptConcludEqt2}) and Corollary \ref{CorSec8}. 
If we assume that the $A(\alpha)$'s are cooperative, 
then, by (\ref{CondSec9}), $A$ is also cooperative. In particular, $A_{ij}u_\alpha^iu_\alpha^j \ge 0$ for 
$i\not= j$, and if we assume in addition that (A3) holds, then we get here again 
that there exists $C > 0$, independent of $\alpha$ and $\delta$, such that (\ref{ProofCptConcludEqt2}) 
holds for all $\alpha$. The same conclusion (\ref{ProofCptConcludEqt2}) 
holds if we assume (A2) and that the $-A(\alpha)$'s are cooperative. 
Here again, the contradiction follows 
from (\ref{ProofCptConcludEqt2}) and Corollary \ref{CorSec8} 
when the dimension $n \ge 7$. In particular, ${\mathcal U}^0 \not\equiv 0$ and $k = 0$ in (\ref{SobDecSec9}). 
Now we assume that $n \ge 4$, that (A2) holds, that the $-A(\alpha)$'s are cooperative for all $\alpha$, and 
that the limit system (\ref{LimitSystSec9}) is fully coupled. By (\ref{ProofCptConcludEqt2}), since we assumed 
(A2) and that the $-A(\alpha)$'s are cooperative for all $\alpha$, there exists $C > 0$, independent of $\alpha$, 
such that
\begin{equation}\label{ProofCptConcludEqt3}
C\int_{B_{x_\alpha}(\delta\sqrt{\mu_\alpha})}\vert{\mathcal U}_\alpha\vert^2dv_g 
+ \frac{(n-2)^2}{2} \left(\sum_{i=1}^pA_i\varphi_i(0) + o(1)\right) \omega_{n-1} \mu_\alpha^{\frac{n-2}{2}} \le 0
\end{equation}
for all $\alpha$. By the pseudo-compactness we proved above, we also have that ${\mathcal U}^0$ in (\ref{SobDecSec9}) 
is nonzero when $n \ge 4$. Then, 
by Theorem \ref{SharpPtwseAsymptThm}, since we assumed that the $-A(\alpha)$'s are 
cooperative for all $\alpha$ and that the limit system (\ref{LimitSystSec9}) is fully coupled, 
we can write that $\sum_{i=1}^pA_i\varphi_i(0)  > 0$, and the contradiction follows 
from (\ref{ProofCptConcludEqt3}). This proves the compactness part of Theorem \ref{ThmCptnessSec9}.
\end{proof}

Summarizing, we let 
$({\mathcal U}_\alpha)_\alpha$ be a bounded sequence in $H_{1,p}^2(M)$ of nontrivial nonnegative 
solutions of (\ref{GenericEqtSec9}). We assume that the ${\mathcal U}_\alpha$'s blow up and, up 
to a subsequence, that the decomposition 
(\ref{SobDecSec9}) is true. 
By the first part of the proof, the pseudo-compactness part, if (A1) 
holds, or (A2) holds and the $-A(\alpha)$'s are cooperative for all $\alpha$, or (A3) holds 
and the $A(\alpha)$'s are cooperative 
for all $\alpha$, then ${\mathcal U}^0$ in (\ref{SobDecSec9}) is nonzero when $n \ge 4$. By the second part of the proof, 
the compactness part, 
we also get that $k = 0$ when the dimension $n \ge 7$ if (A1) holds, 
or (A2) holds and the $-A(\alpha)$'s are cooperative for all $\alpha$, or (A3) holds 
and the $A(\alpha)$'s are cooperative 
for all $\alpha$. Moreover, when $n \ge 4$, still by the second part of the proof, 
if we assume that (A2) holds, that the $-A(\alpha)$'s are cooperative for all $\alpha$, and that the limit system 
(\ref{LimitSystSec9}) is 
fully coupled, then, we again get that $k = 0$ in (\ref{SobDecSec9}).
Since $({\mathcal U}_\alpha)_\alpha$ 
is arbitrary, this ends the proof of Theorem \ref{ThmCptnessSec9}.
\end{proof}

One more remark with respect to Theorem \ref{ThmCptnessSec9} is as follows.

\medskip\noindent{\bf Remark 9.2:} We have already mentionned, in Remark 9.1 just before 
stating Theorem  \ref{ThmCptnessSec9}, 
that we cannot get compactness or pseudo-compactness 
for systems like (\ref{GenericEqtSec9}) if we only assume (A2) or (A3). 
In the same spirit, we mention that, without further assumptions, we 
cannot hope as well to get compactness or pseudo-compactness with a 
mix of conditions like (A2) and (A3) on the diagonal 
of the limit matrix $A$. Let, for instance, $(u_\alpha)_\alpha$ be a blowing-up sequence of solutions of the 
Yamabe equation on the sphere, and 
${\mathcal U}_\alpha = (u_\alpha,u_\alpha,u_\alpha)$ for all $\alpha$. Let also $A$ be the matrix
\begin{equation}\label{ExSphereCptnessSect}
A = \left(
\begin{matrix}
a & b & 0\\
b & c & -d\\
0 & -d & e
\end{matrix}
\right)\hskip.1cm ,
\end{equation}
where $a$, $b$, $c$, $d$, $e$ are positive real numbers. If $\lambda_n = \frac{n(n-2)}{4}$, 
hence $\lambda_n = \frac{n-2}{4(n-1)}S_{g_0}$, we choose 
$a$, $b$, $c$, $d$, $e$ such that $a + b = \lambda_n$, $b + c = d + \lambda_n$, and 
$e = d + \lambda_n$. For any $\alpha$, the map ${\mathcal U}_\alpha$ 
is a positive solution of the $3$-system
$$\Delta_{g_0}u^i + \sum_{j=1}^pA_{ij}u^j  = (u^i)^{2^\star-1}$$
in $S^n$, for all $i = 1,2,3$, where the $A_{ij}$'s are the components of the matrix $A$ in 
(\ref{ExSphereCptnessSect}). The $u_\alpha$'s blow up with a 
zero (pointwise) limit. It follows that the system is not pseudo-compact (and thus, not compact as well). 
However, by choosing $d < b$ sufficiently small, the operator $\Delta_{g_0}^3 + A$ is coercive, 
$A_{ii} < \lambda_n$ for $i = 1,2$, and $A_{33} > \lambda_n$. Compactness when we mix 
conditions like (A2) and (A3) 
on the diagonal is false in general. We get compactness type behaviour under such assumptions 
when the energy of the sequence is 
of minimal type, see Hebey \cite{Heb2}, or, of course, if we 
assume that $A^\alpha_{ij} \equiv 0$ for all $i \not= j$ and all $\alpha$.

\medskip The following easy corollary of Theorem \ref{ThmCptnessSec9} shows that the 
theorem is sharp when regarded on the unit sphere. Corollary \ref{corSphereSec9} below states that 
the trivially coupled system on the sphere, consisting of $p$ 
copies of the Yamabe equation, is, in some sense, the only system on the sphere which is not compact. 
Concerning notations, we let $(S^n,g_0)$ be the unit $n$-sphere, and 
$A: S^n \to M_p^s(S^n)$ be a smooth map. We assume that $A(x)$ is positive in the sense of bilinear forms 
for all $x \in S^n$, and for $t \in {\mathbb R}$, we define 
\begin{equation}\label{CorExSec9}
A_t = \frac{n(n-2)}{4}Id_p + tA\hskip.1cm ,
\end{equation}
where $Id_p$ is the $p\times p$-identity matrix in $M^s_p({\mathbb R})$. 
We let $\Lambda_0$ be the maximum over $i = 1,\dots,p$ and $x \in S^n$ 
of the eigenvalues $\lambda_i(x)$ of $A(x)$. 
The operator $\Delta_g^p+A_t$ is coercive for all $t$ in the interval $I_0 = (-\Lambda,+\infty)$, 
where $4\Lambda_0\Lambda = - n(n-2)$. For $t \in I_0 = (-\Lambda,+\infty)$, 
where $4\Lambda_0\Lambda = - n(n-2)$, we consider the systems
\begin{equation}\label{CorExSec9EqtSyst}
\Delta_gu_i + \sum_{j=1}^pA_{ij}^t(x)u_j = \vert u_i\vert^{2^\star-2}u_i
\end{equation}
in $M$, for all $i = 1,\dots p$, where the $A_{ij}^t$'s are the components of $A_t$, and 
$A_t$ is given by (\ref{CorExSec9}). Corollary \ref{corSphereSec9} isolates 
the trivially coupled system on the sphere as the only system 
for which pseudo-compactness or compactness fails in families of systems like (\ref{CorExSec9EqtSyst}), 
providing an illustration that 
blow-up phenomena and noncompactness occur only with the geometric equation (in the case of the sphere), 
or perturbations of the geometric equation (by Proposition \ref{PropExBlowUpSec9}). 

\begin{cor}\label{corSphereSec9} Let 
$(S^n,g_0)$ be the unit $n$-sphere, and 
$A: S^n \to M_p^s(S^n)$ be a smooth map such that $A(x)$ is positive in the sense of bilinear forms 
for all $x \in S^n$. For any $t \in I_0\backslash\{0\}$, 
where $I_0 = (-\Lambda,+\infty)$ is as above, the system (\ref{CorExSec9EqtSyst}) is pseudo-compact 
when $n \ge 4$, and compact when $n \ge 7$. On the other hand, when $t = 0$, 
(\ref{CorExSec9EqtSyst}) is neither  
compact nor pseudo-compact. In particular,  
the trivially coupled system (\ref{CorExSec9EqtSyst}) when $t=0$, 
corresponding to $p$-copies of the Yamabe equation on the sphere,  
is the only system 
in the family (\ref{CorExSec9EqtSyst}), $t  > - \Lambda$, which is not compact, respectively not pseudo-compact (depending on whether 
$n \ge 4$ or $n \ge 7$).
\end{cor}

Pseudo-compactness and compactness in the corollary follow from 
Theorem \ref{ThmCptnessSec9}. Since we assumed that 
$A(x)$ is positive in the sense of bilinear forms 
for all $x \in S^n$, 
(A1) holds for the 
$A_t$'s when $t \not= 0$ and we can apply 
the theorem. On the other hand, that (\ref{CorExSec9EqtSyst}) is neither  
compact nor pseudo-compact when $t=0$ follows from the observation 
that (\ref{CorExSec9EqtSyst}) when $t=0$ 
consists of $p$-copies of the Yamabe equation on the sphere 
which, as is well known, possesses sequences of solutions which blow up with one bubble and a 
zero (pointwise) limit in their Sobolev decomposition. 

\medskip\noindent{\bf Acknowledgements:} The author is indebted to 
Olivier Druet and Fr\'ed\'eric Robert 
for their valuable comments on the manuscript.

\end{document}